\documentclass[a4paper,oneside,reqno]{amsart}
\usepackage[T1]{fontenc}
\usepackage{amsmath}
\usepackage{amsthm}
\usepackage{amssymb}
\usepackage{mathabx}
\usepackage{times}
\usepackage{bm}
\usepackage{amscd}
\usepackage[latin2]{inputenc}
\usepackage{t1enc}
\usepackage[mathscr]{eucal}
\usepackage{indentfirst}
\usepackage{graphicx}
\usepackage{graphics}
\usepackage{pict2e}
\usepackage{epic}

\usepackage{epstopdf} 
\usepackage[colorlinks,linkcolor=blue]{hyperref}
\usepackage[title]{appendix}
\usepackage[capitalise,noabbrev]{cleveref}
\usepackage{todonotes}
\usepackage{orcidlink}

 \pdfoutput=1
\numberwithin{equation}{section}
\crefformat{equation}{(#2#1#3)}
\crefmultiformat{equation}{(#2#1#3)}{ and~(#2#1#3)}{, (#2#1#3)}{ and~(#2#1#3)}
\crefrangeformat{equation}{(#3#1#4) to~(#5#2#6)}

\usepackage[normalem]{ulem}

\allowdisplaybreaks

\theoremstyle{plain}
\newtheorem{Thm}{Theorem}[section]
\newtheorem*{Thm*}{Theorem}
\newtheorem{Lem}[Thm]{Lemma}

\newtheorem{Prop}[Thm]{Proposition}

\theoremstyle{definition}

\newtheorem{Rem}[Thm]{Remark}

\AddToHook{env/Lem/begin}{\crefalias{Thm}{Lem}}
\AddToHook{env/Prop/begin}{\crefalias{Thm}{Prop}}
\AddToHook{env/Cor/begin}{\crefalias{Thm}{Cor}}
\AddToHook{env/Def/begin}{\crefalias{Thm}{Def}}
\AddToHook{env/Conj/begin}{\crefalias{Thm}{Conj}}
\AddToHook{env/Rem/begin}{\crefalias{Thm}{Rem}}
\AddToHook{env/Ex/begin}{\crefalias{Thm}{Ex}}

\crefname{Lem}{Lemma}{Lemmas}
\crefname{Prop}{Proposition}{Propositions}
\crefname{Cor}{Corollary}{Corollaries}
\crefname{Def}{Definition}{Definitions}
\crefname{Conj}{Conjecture}{Conjectures}
\crefname{Rem}{Remark}{Remarks}
\crefname{Ex}{Example}{Examples}

\newcommand{\p}{\partial}
\newcommand{\R}{\mathbb{R}}
\newcommand{\e}{{\red{\varepsilon}}}

\newcommand{\E}{\mathbf{e}}

\newcommand{\Torus}{\mathbb{T}}

\newcommand{\dv}{\mathsf{div}\,}
\newcommand{\vs}{\mathsf{vs}}
\newcommand{\od}{\flat}
\newcommand{\md}{\sharp}
\newcommand{\lap}{\triangle}
\newcommand{\dnab}{\cdot \nabla}

\newcommand{\rhob}{\bar{\rho}}

\newcommand{\rhot}{\widetilde{\rho}}

\newcommand{\ub}{\bar{u}}
\newcommand{\ut}{\widetilde{u}}
\newcommand{\uv}{\mathbf{u}}
\newcommand{\mv}{\mathbf{m}}
\newcommand{\wv}{\mathbf{w}}
\newcommand{\wt}{\widetilde{w}}
\newcommand{\wvt}{\widetilde{\mathbf{w}}}

\newcommand{\uvb}{\bar{\mathbf{u}}}

\newcommand{\uvt}{\widetilde{\mathbf{u}}}
\newcommand{\uvs}{\mathbf{u}^\vs}
\newcommand{\mvt}{\widetilde{\mathbf{m}}}

\newcommand{\mt}{\widetilde{m}}

\newcommand{\mut}{\tilde{\mu}}
\newcommand{\xp}{x_{\perp}}

\newcommand{\rv}{\mathbf{r}}
\newcommand{\db}{\mathbf{d}}

\newcommand{\Fv}{\mathbf{F}}

\newcommand{\Lv}{\mathbf{L}}
\newcommand{\Qv}{\mathbf{Q}}
\newcommand{\Zv}{\mathbf{Z}}

\newcommand{\vt}{\widetilde{v}}

\newcommand{\zetab}{\bar{\zeta}}

\newcommand{\andd}{\quad \text{and} \quad }

\newcommand{\Ac}{\mathcal{A}}
\newcommand{\Bc}{\mathcal{B}}
\newcommand{\Ec}{\mathcal{E}}
\newcommand{\Nc}{\mathcal{N}}

\newcommand{\Rc}{\mathfrak{E}}
\newcommand{\vv}{\mathbf{v}}

\newcommand{\red}{\textcolor{red}}

\newcommand{\abs}[1]{\big\lvert#1\big\rvert}

\newcommand{\norm}[1]{\big\lVert#1\big\rVert}

\newcommand{\cdn}{\cdot\nabla}
\newcommand{\vvt}{\widetilde{\mathbf{v}}}

\begin{document}
	
	
\title[Vortex sheets for compressible Navier-Stokes]{Global nonlinear stability of vortex sheets for the Navier-Stokes equations with large data}
 
\author[Yuan]{Qian Yuan}
\address{State Key Laboratory of Mathematical Sciences, Academy of Mathematics and Systems Science, Chinese Academy of Sciences, Beijing, PRC.}
\email{qyuan@amss.ac.cn}

\author[Zhao]{Wenbin Zhao$^{\dagger}$
\orcidlink{0000-0003-0804-1934}
}
\address{School of Mathematics, Renmin University of China, Beijing, PRC.}
\email{wenbizhao2@ruc.edu.cn}

\thanks{$\dagger$ Corresponding author.}



\begin{abstract}

This paper concerns the global nonlinear stability of vortex sheets for the Navier-Stokes equations. When the Mach number is small, we allow both the amplitude and vorticity of the vortex sheets to be large. 
We introduce an auxiliary flow and reformulate the problem as a vortex sheet with small vorticity but subjected to a large perturbation. Based on the decomposition of frequency, the largeness of the perturbation is encoded in the zero modes of the tangential velocity. We discover an essential cancellation property that there are no nonlinear interactions among these large zero modes in the zero-mode perturbed system. 
This cancellation is owing to the shear structure inherent in the vortex sheets.  Furthermore, with the aid of the anti-derivative technique, we establish a faster decay rate for the large zero modes. 
These observations enable us to derive the global estimates for strong solutions that are uniform with respect to the Mach number. 
As a byproduct, we can justify the incompressible limit.

\end{abstract}

\maketitle



\section{Introduction}

The three-dimensional compressible isentropic Navier-Stokes equations read
\begin{equation}\label{NS-1}
	\begin{cases}
		\p_t \rho + \dv \mv = 0, \\
		\p_t \mv + \dv ( \rho \uv\otimes\uv ) + \nabla p = \mu \lap \uv + \big(\mu+\lambda \big) \nabla \dv \uv,
	\end{cases} 
	\quad x\in\R^3,\, t>0, 
\end{equation}
where $ \rho $ is the density, $ \uv $ is the velocity, and $ p=p(\rho)=\rho^{\gamma} $ is the pressure with $\gamma>1$. 
The viscosity coefficients $ \mu $ and $ \lambda $ are assumed to satisfy that
\begin{equation}\label{viscous}
	\mu > 0 \andd \lambda + \mu \geq 0,
\end{equation}
which covers the physical requirements.
The system \cref{NS-1} describes the motion for a viscous isentropic compressible fluid.
When $\mu=\lambda=0$ in \eqref{NS-1}, this is the compressible isentropic Euler equations. 

In an inviscid fluid, a vortex sheet is a phenomenon of an interface, across which the tangential velocity of the flow is discontinuous while the normal velocity and the pressure are continuous. 
A planar vortex sheet in $\R^3$ is a piece-wise constant solution to the Euler equations, and up to a Galilean transformation, it can be written as 
\begin{equation}\label{vs}
	(\rhob, \uvb^\vs)(x_3,t) = \begin{cases}
		(\rhob,\, -\uvb), \qquad x_3 <0, \\
		(\rhob,\, \uvb), \qquad \quad x_3 >0,
	\end{cases}
\end{equation}
where $\rhob>0$ and $ \uvb =  (\ub_1, \ub_2, 0) \in \R^3 $ are any given constants.

\subsection{Prior works on instability/stability of vortex sheets}

The inviscid vortex sheets \cref{vs} are usually subject to the Kelvin-Helmholtz instability. It is known that 3D vortex sheets are violently unstable while 2D vortex sheets (such as \cref{vs} with $\uvb = (\ub_2,0) \in \R^2$) are weakly stable under the supersonic condition (\cite{Coulomb2004,Coulomb2008,FM1963,Miles1958}). In the incompressible limit with the sound speed tending to infinity, the incompressible vortex sheets are always unstable. 
Mathematical analysis of the incompressible vortex sheets in the analytic class has been established in \cite{SSBF1981,Wu2006}. 

When viscosity is
present, the study of the vortex layers are of great importance for both mathematics and mechanics, e.g. the mixing of flows and the separation of boundary layers (\cite{MB2002,Ruban2018}). 
As for the regularization of the inviscid vortex sheets, Wu raised the following open problem in \cite{Wu2002}:
\begin{quotation}
	``...the vortex sheet in general fails to be
	a curve beyond the initial time for general data. Therefore it becomes interesting to study the vortex layers or considering the effects of viscosity. 
    ''
\end{quotation}
In \cite{CLS2020}, Caflisch et al. considered the 2D inviscid vortex layers with small thickness along a given curve, the vorticity of which essentially concentrated along this curve and decayed exponentially away from it. They showed that the center of the layer could be well approximated by the vortex sheet for a short time. When taking viscosity into account, this type of vorticity distribution would be more appropriate. By choosing the initial thickness of the 2D vortex layer proportional to the square root of the viscosity, 
\begin{equation*}
    \mathrm{thickness}\sim \mu^{1/2},
\end{equation*}
Caflisch et al. investigated the small viscosity limit in \cite{CS2006} and the roll-up process in \cite{CGSS2022}.  
The interested readers are advised to consult \cite{CGSS2022,CLS2020} and references within for more details.   

\subsection{Large time behavior of vortex sheets for Navier-Stokes equations}

On the large time scale, the vortex layers have the thickness $\sim (\mu t)^{1/2}. $ 
In fact, for the planar vortex sheet \cref{vs}, the associated viscous vortex layer is given explicitly by
\begin{equation}\label{profile}
    \big(\rhob,\, \uv^\vs\big)(x_3,t) = \Big(\rhob,\, \Theta\Big(\frac{x_3}{\sqrt{t}}\Big) \uvb\Big), \quad t\geq 0,
\end{equation}
with
\begin{equation}\label{Theta}
    \Theta(\xi) := \frac{2}{\sqrt{\pi}} \int_{0}^{\frac{1}{2} \sqrt{\frac{\rhob}{\mu}} \xi} e^{-\eta^2} d \eta.
\end{equation}
The vortex layer \cref{profile} is a solution to the Navier-Stokes equations \cref{NS-1}, which approaches the vortex sheet \cref{vs} as $t\to 0+,$ while moves away from it as $t\to +\infty$ due to the fact, 
\begin{equation*}
    \begin{aligned}
        \norm{(\rhob, \uv^\vs) - (\rhob, \uvb^\vs)}_{L^p(\R; dx_3)} & \sim \abs{\uvb} (\mu t)^{\frac{1}{p}} \qquad \forall p \in [1, +\infty).
    \end{aligned}
\end{equation*}
The vorticity of the vortex layer \eqref{profile} satisfies the Gaussian distribution along $x_3$:
\begin{equation}\label{vor-vs}
    \nabla\times\uvs = (\rhob/\pi)^{\frac{1}{2}} (\mu t)^{-\frac{1}{2}} e^{-\frac{\rhob}{4\mu t}|x_3|^2} \mathbf{e}_3 \times \uvb \quad \text{with } \mathbf{e}_3 =(0,0,1).
\end{equation}
Recently, \cite{HXXY2023} showed that the vortex layer with small initial vorticity, that is, the solution
\begin{equation}\label{profile-t0}
    \big(\rhob,\, \uv^\vs\big)(x_3,t+t_0) = \Big(\rhob,\, \Theta\Big(\frac{x_3}{\sqrt{t+t_0}}\Big) \uvb\Big),
\end{equation}
with a suitably large constant $t_0>0,$ is nonlinearly stable for the compressible Navier-Stokes equations, \cref{NS-1}. 
Although the largeness assumption, $t_0 \gg 1, $ in \cite{HXXY2023} is essential in the proof, the amplitude of the vortex layer, $\abs{\uvb}$,
can be large.
Compared to the results \cite{Coulomb2004,Coulomb2008,Ebin1988,Serre1999} for the Euler equations, the global nonlinear stability obtained in \cite{HXXY2023} indicates a strong stabilizing effect of the viscosity on the vortex sheets.

\subsection{Main Results} 
In this paper, we consider a general vortex layer \cref{profile-t0} with arbitrary $\rhob>0$, $\uvb =(\ub_1,\ub_2,0)\in \R^3$ and $t_0>0 $. 
It should be noted that the vorticity of the vortex layer can be large. 
The study of the vortex layers with large vorticity is an important step for the future study of the small viscosity limit problem. 
In addition to largeness of the background flow, the analysis of the paper also allows the initial perturbations to be large in some sense. 

Now we formulate the problem. 
Consider the compressible isentropic Navier-Stokes equations in a non-dimensional way,
\begin{equation}\label{NS}
	\begin{cases}
		\p_t \rho^\e + \dv \mv^\e = 0, \\
		\p_t \mv^\e + \dv ( \rho^\e \uv^\e\otimes\uv^\e ) + \frac{1}{\e^2} \nabla p(\rho^\e) = \mu \lap \uv^\e + \big(\mu+\lambda \big) \nabla \dv \uv^\e,
	\end{cases} 
\end{equation}
where $\e$ denotes the Mach number, namely the ratio of a characteristic velocity to the sound speed in the fluid. The system \cref{NS} can be obtained by scaling the variables 
\begin{equation}\label{scaling}
    t \to \e^2 t, \quad x \to \e x, \quad \uv \to \e \uv.
\end{equation}
We refer to \cite{Alazard06} for other changes of variables instead of \cref{scaling}.

In the infinitely long nozzle domain 
\begin{equation*}
    \Omega := \big\{ x=(x_1,x_2,x_3): (x_1,x_2) \in \Torus^2, x_3 \in \R \big\},
\end{equation*}
we consider a Cauchy problem for \cref{NS}, in which the initial data is a general perturbation of the vortex layer \cref{profile-t0}, namely, 
\begin{equation}\label{ic}
	\begin{aligned}
		(\rho^\e, \uv^\e)(x,t=0) & = (\rho_0^\e, \uv_0^\e)(x)  \\ 
		& := (\rhob, \uv^\vs)(x_3,t=0) + (\e b_0, \vv_0)(x), \quad x \in \Omega, 
	\end{aligned}
\end{equation}
where $ (b_0, \vv_0) = (b_0, v_{01}, v_{02}, v_{03}) $ belongs to $ H^3(\Omega)$.

\vspace{.1cm}

Before stating the main theorem, we introduce some notations.
\begin{itemize}
    \item For any vector $\vv=(v_1,\,v_2,\,v_3)\in\R^3, $ denote $\vv_\perp = (v_1,\,v_2)$.

    \item For any $ f(x) \in L^\infty(\Omega) $, denote $ f^\od $ as its \textit{zero mode},
    \begin{equation}\label{od}
	f^\od(x_3) := \int_{\Torus^2} f(\xp, x_3) d\xp,
    \end{equation}
    and $ f^\md $ as its \textit{non-zero mode},
    \begin{equation}\label{md}
	f^\md(x) := f(x) - f^\od(x_3).
    \end{equation}

    \item Denote $\langle x_3\rangle := (1+x_3^2)^{\frac{1}{2}}.$ For $\alpha>0,$ define 
\begin{equation}
    \begin{aligned}
        H^s_\alpha(\Omega) & := \Big\{ f\in H^s(\Omega): \norm{f^\od}_{L^2_\alpha(\R)}:= \norm{\langle x_3\rangle^{\alpha} f^\od}_{L^2(\R)} < \infty \Big\},
    \end{aligned}
\end{equation}
and
\begin{equation}
	\norm{f}_{H^s_\alpha(\Omega)} := \norm{f^\od}_{L^2_\alpha(\R)} + \norm{f}_{H^s(\Omega)}.
\end{equation}
\end{itemize}

Now we are ready to state the main results of this paper, 

\begin{Thm}\label{Thm}
    Given a vortex layer \cref{profile-t0}, where $\rhob>0, t_0>0, \ub_1 $ and $ \ub_2$ in \cref{profile-t0} are arbitrary constants.
    Assume that the initial perturbations $(b_0, \vv_0) \in H^3_{3/2}$, and denote
    \begin{align}
        M_0 & := \norm{(b_0, \vv_0)}_{H^3_{3/2}(\Omega)}, \label{M-0} \\
        \chi & := \norm{(b_0, v_{03})^\od}_{H^1_{3/2}(\R)} + \norm{(b_0, \vv_0)^\md}_{H^1(\Omega)}. \label{chi}
    \end{align}
    Then there exist 
    \begin{itemize}
        \item $ \e_0>0 $ and $ \chi_0>0 $, depending on $\mu, \lambda, \rhob, |\uvb| $ and $t_0$,
        \item a generic constant $k_0>0,$ depending only on the space dimension,
    \end{itemize}
    such that given any $M_0 >0,$ if
    \begin{equation}\label{small}
        \begin{aligned}
            0< \e \leq \e_0, \quad \chi \leq \chi_0 \andd (\e + \chi) M_0^{k_0} \leq 1,
        \end{aligned}
    \end{equation}
    then the Cauchy problem \cref{NS}, \cref{ic} admits a unique strong solution, $ (\rho^\e,\uv^\e), $  globally in time. Moreover, the perturbations, 
	\begin{equation}\label{pert-bv}
		b^\e :=  \e^{-1}(\rho^\e - \rhob), \quad \vv^\e := \uv^\e - \uv^\vs,
	\end{equation}
	satisfy the uniform (with respect to $\e$) estimates,
    \begin{equation}\label{Thm-L2}
        \begin{aligned}
            & \sup_{t\geq 0} \norm{(b^\e, \vv^\e)}_{H^3(\Omega)}^2 \\
            & \qquad + \int_0^\infty
            \Big(\norm{\nabla b^\e}_{H^{2}(\Omega)}^2 + \norm{ \nabla \vv^\e }_{H^{3}(\Omega)}^2 \Big) dt \leq C,
        \end{aligned}
        \end{equation}
        and
        \begin{equation}\label{rate-Thm}
		\norm{(b^\e, \vv^\e)}_{L^\infty(\Omega)} \leq C (t+1)^{-\frac{1}{2}} \qquad \forall t\geq 0,
	\end{equation}
    where $ C>0 $ is a constant, independent of $ \e $ and $ t. $
\end{Thm}

\begin{Rem}\label{Rem-large}

We give two remarks on the \textit{largeness} of the initial data. 
\begin{itemize}
    \item For the background flow \cref{profile-t0}, both the amplitude and initial vorticity can be large.

\item Comparing \eqref{M-0} and \eqref{chi}, the tangential velocity contains a large perturbation around the background flow, since the associated zero mode, $\vv_{0\perp}^\od = (v_{01}, v_{02})^\od$, can be arbitrary in $H^3_{3/2}(\R).$

\end{itemize}

\begin{Rem}
The results in \cref{Thm} still hold true in the 2D domain $\Omega = \Torus \times \R$,
just by letting $\ub_1 = 0, v_{01} = 0$ and $(b_0, v_{02}, v_{03})$ be independent of $x_1$.
\end{Rem}

\begin{Rem}
The background vortex sheet can degenerate to a constant state, that is, the case $\uvb = 0 $ is included.
\end{Rem}

With the global stability which is uniform with respect to the Mach number $\e$, the incompressible limit of the solutions $(\rho^\e, \uv^\e)$ can follow from the classical arguments in \cite{Metivier2001,Alazard06}. In addition, one can obtain the global nonlinear stability of vortex layers for the incompressible Naiver-Stokes equations with large data. 
We include these results and the proofs in \cref{Sec-lim} for the sake of completeness.

\vspace{.2cm}

\end{Rem}

\subsection{Difficulties and novelties in the proof}\label{Sec-novelty}

We outline the main difficulties, new observations, and novel ingredients in the proof. To simplify the notations, we omit the upper index $\e.$
Denote the perturbation $U = (b, \vv) := (\e^{-1} (\rho- \rhob), \uv - \uv^\vs)$. By \cref{NS}, the system of $U$ is formally given by
\begin{equation}\label{system-formal}
    \begin{cases}
        (\p_t + \uvs\cdn) b +\frac{1}{\e} \rhob \, \dv \vv + \mathbb{N}_1(U, \nabla U) = 0, \\
        (\p_t + \uvs\cdn)\vv
        + \frac{1}{\e}p'(\rhob)\nabla b + \mathbb{L}(U) + \mathbb{N}_2(U, \nabla U) \\
        \qquad  = \frac{\mu}{\rhob} \lap \vv + \frac{\mu+\lambda}{\rhob} \nabla\dv \vv + O(\e) (\cdots).
    \end{cases}
\end{equation}
where
\begin{equation}\label{diffi-lin}
    \mathbb{L}(U) = \vv \cdn\uvs = \p_3\uvs \, v_3,
\end{equation}
and
\begin{equation}\label{nonlin}
    \begin{aligned}
        \mathbb{N}_1(U, \nabla U) = \dv(b\vv), \quad \mathbb{N}_2(U, \nabla U) = \vv\cdn\vv + 2 p''(\rhob) b \nabla b.
    \end{aligned}
\end{equation}

$\spadesuit$ Difficulties in linear level.
Due to the slow decay rate of the background flow,
$$\abs{\p_3 \uv^\vs} \lesssim \abs{\uvb} (t+t_0)^{-1/2},$$ 
the zero-order term \cref{diffi-lin} results in an unbounded $L^2$-estimate of $\vv.$
To overcome this difficulty, we decompose the perturbation into the \textit{zero modes} and \textit{non-zero modes}, and prove their estimates separately.

For the zero modes, following \cite{Yuan2023s,HXXY2023}, we apply the \textit{anti-derivative} method and use an \textit{effective momentum}.
Roughly speaking, if it holds that $ \vv^\od \sim \p_3 \mathbf{V} $ with some $ \mathbf{V} \in L^2(\R)$, then the $L^2$-estimate of $\vv^\od$ gives that
\begin{equation}\label{diff-od}
    \frac{d}{dt} \norm{\vv^\od}_{L^2}^2 + \norm{\p_3 \vv^\od}_{L^2}^2 \lesssim \norm{\p_3 \mathbf{V}}_{L^2}^2 + \cdots.
\end{equation}
Thanks to the inherent structures of the Navier-Stokes equations, the system of $\mathbf{V}$ has no zero order terms of $\mathbf{V}$.
Thus, the derivative $\norm{\p_3\mathbf{V}}_{L^2}$ on the right hand side of \cref{diff-od} can be bounded by the dissipation. 
On the other hand, we set the anti-derivative variable $\mathbf{V}$ as an effective one as in \cite{HXXY2023} to overcome the difficulty arising from the large amplitude of the vortex sheet.
We also refer to \cite{Ilin1960,MN1985,G1986,HMX} for the anti-derivative technique in the study on one-dimensional shocks and contact discontinuities. 

However, in the estimates of non-zero modes,
\begin{equation}\label{diff-md}
    \frac{d}{dt} \norm{\vv^\md}_{L^2}^2 + \norm{\nabla\vv^\md}_{L^2}^2 \lesssim \abs{\uvb} (t+t_0)^{-\frac{1}{2}} \norm{\vv^\md}_{L^2}^2 + \cdots,
\end{equation}
we cannot construct anti-derivatives in multiple dimensions and this anti-derivative method is not effective for the non-zero modes.
In \cite{HXXY2023} where $t_0$ is assumed to be large, the coefficient $\abs{\uvb} t_0^{-\frac{1}{2}}$ in \cref{diff-md} is small. The bad term on the right-hand side can be controlled by the dissipation since the non-zero mode satisfies the Poincar\'{e} inequality, $\norm{\vv^\md}_{L^2} \lesssim \norm{\nabla\vv^\md}_{L^2}$.

To cope with a general $t_0>0$ in this paper, we shall introduce an \textit{auxiliary flow} $\uvt^\vs(x_3,t):=\uv^\vs(x_3, t+\Lambda)$ with $\Lambda>0$ being a large constant. 
We consider $(\rhob, \uvt^\vs)$ as a new background solution, where the velocity field can be decomposed as 
\begin{equation}\label{pert-add}
    \begin{aligned}
	\uv=\uvs+\vv & = \uvt^\vs+(\uvs-\uvt^\vs)+\vv := \uvt^\vs+ \vvt.
\end{aligned}
\end{equation}
The system of $(b,\,\vvt)$ shares the similar structure as \eqref{system-formal}, while the convection term $\vvt\cdn\uvt^\vs$ has a small coefficient as in \cite{HXXY2023}.

However, the cost of the replacement is the extra perturbation of the tangential velocity field, $\uv_\perp^\vs - \uvt_\perp^\vs$, which makes the associated zero mode 
$\vvt_\perp^\od$ as a large perturbation.

\vspace{.2cm}

$\spadesuit$ Difficulties in nonlinear level.
To deal with the large nonlinearities in \cref{nonlin}, we introduce the \textit{key energy functional}, which excludes the zero mode associated with the tangential velocity,
\begin{equation}\label{intro-small}
\Ec^*(t) := \norm{(b,\tilde{v}_3)^\od}_{H^1(\R)}^2 + \norm{(b, \vvt)^\md}_{H^1(\Omega)}^2.
\end{equation}
It is noted that the extra perturbation, $\uv_\perp^\vs-\uvt_\perp^\vs$, in \cref{pert-add} does not affect the initial value of $\Ec^*.$ Thus, by the initial assumption \cref{chi}, the energy $\Ec^*$ is small at $t=0.$ 
One of our main efforts is to prove the smallness of $\Ec^*(t)$ globally in time. 

In the estimates of the zero modes, we find that the nonlinearities, $(\mathbb{N}_i(U, \nabla U))^\od$ for $i=1,2$, excluding the nonlinear interactions among the large modes $\vvt_\perp^\od$, are all interactions of $U$ and the small perturbations in \cref{intro-small}. 
This cancellation is owing to the shear structure inherent in the vortex sheets. Thus, the smallness of \cref{intro-small} is sufficient to control these nonlinearities for the zero-mode estimates.

However, this is not the situation for the non-zero modes.
In the system of $ \vvt^\md,$ the nonlinear convection $(\vvt \cdot \nabla \vvt)^\md$ contains a nonlinear interaction of the large zero mode and non-zero mode, $ \p_3 \vvt^\od \tilde{v}_3^\md, $ resulting in 
\begin{equation}\label{diff-md-n}
    \frac{d}{dt} \norm{\vvt^\md}_{L^2}^2 + \norm{\nabla\vvt^\md}_{L^2}^2 \lesssim \norm{\p_3 \vvt_\perp^\od}_{L^\infty} \norm{\vvt^\md}_{L^2}^2 + \cdots,
\end{equation}
where $\p_3 \vvt_\perp^\od$ has no smallness and its decay rate $(t+1)^{-3/4}$, which is obtained through the energy method, is not sufficient to achieve the global boundedness of $\norm{\vvt^\md}_{L^2}^2.$

As one novelty of this paper, we find that the large part of the zero mode $\p_3 \vvt_\perp^\od$ decays at a faster rate $(t+1)^{-\frac{5}{4}} $ in $L^\infty$, namely, it holds that
\begin{equation}\label{intro-key}
    \norm{\p_3 \vvt_\perp^\od}_{L^\infty(\R)} \lesssim (t+1)^{-\frac{5}{4}} + \text{(small parameters)}.
\end{equation}
The refined estimate plays a key role in the proof and it is derived from two observations. 
First, owing to the cancellation property of the vortex sheets, we find that $\vvt_\perp^\od$ satisfies the parabolic-type equation with small source terms,
\begin{equation}\label{equ-vi}
\begin{cases}
    \p_t \vvt_\perp^\od = \frac{\mu}{\rhob} \p_3^2 \vvt_\perp^\od + \text{(small source)}, \\
    \vvt_\perp^\od|_{t=0} = (\uvs_\perp -\uvt_\perp^\vs)|_{t=0} + \vv_{0\perp}^\od, \ \text{which is large.}
\end{cases}
\end{equation}
This means that the largeness of $\vvt_\perp^\od$ is totally from the initial data.
The second observation is the existence of the $L^2$-integrable anti-derivative variables such that $\vvt_\perp^\od \sim \p_3 \mathbf{V}_\perp. $ Formally, $\mathbf{V}_\perp$ satisfies the integrated system and data of \cref{equ-vi}, and the use of the Green's function yields \cref{intro-key}.

\textit{Role of small Mach number.} 
Throughout this paper, we shall use the smallness of the Mach number $\e$ to deal with many complicated nonlinear interactions in the energy estimates.
In particular, the small Mach number plays a key role in achieving the global smallness of the energy functional, \cref{intro-small}. Indeed, the nonlinear pressure results in a bad second-order term in the estimate of $\Ec^*,$ which gives
\begin{equation}\label{diff-E*}
	\begin{aligned}
		\frac{d}{dt} \Ec^*(t) + \Nc^*(t) & \lesssim \norm{b^\md \nabla b^\md}_{L^2}^2 + \cdots \\
		& \lesssim \Ec^*(t) \norm{\nabla^2 b^\md}_{L^2}^2 + \cdots,
	\end{aligned}
\end{equation}
where $\Nc^* \sim \norm{\nabla b}_{L^2}^2 + \norm{\p_3 \vt_3^\od}_{H^1}^2 + \norm{\nabla \vvt^\md}_{H^1}^2 $ is the sum of low-order dissipation terms associated with $\Ec^*$, and $\norm{\nabla^2 b^\md}_{L^2}^2$ on the right-hand side is a large energy which does NOT belong to $\Ec^*.$ 
One key point in the proof is the second-order estimate of the density,
\begin{equation}\label{intro-b2}
    \frac{d}{dt} \big(\e^2 \norm{\nabla^2 b}_{L^2}^2 \big) + \norm{\nabla^2 b}_{L^2}^2 \lesssim \Nc^*(t) + O(\e).
\end{equation}
Plugging \cref{intro-b2} into \cref{diff-E*} and use the smallness of $\e,$ the large energy on the right-hand side of \cref{diff-E*} can be controlled, and the global smallness of $\Ec^*$ can be achieved.

\vspace{.2cm}

\subsection{Outline of the paper}
In \cref{Sec-ansatz}, we introduce the auxiliary flow and show the construction of the ansatz such that the anti-derivative variables of the perturbations exist in $L^2$. In addition, we reformulate the perturbed system and the theorem. 
In \cref{Sec-apriori}, we present the a priori estimates, and outline the bounds and decay rates of the perturbations.
In \cref{Sec-pf-apriori}, we outline the main steps of the a priori estimates, and postpone the detailed proof of each step to Sections \ref{Sec-od}-\ref{Sec-orn}, respectively.
The last section is devoted to the incompressible limit result.

\vspace{.2cm}

\underline{Notations}: 
\begin{itemize}

\item Since $t_0>0$ in \cref{profile-t0} is a fixed constant throughout the paper, we assume that $ t_0=1 $ for readers' convenience. 

\item For the pressure law $p=p(\rho)$, let $\varpi(\cdot, \cdot\cdot): \R^2 \to \R $ denote that
	\begin{equation}\label{varpi}
		\begin{aligned}
			\varpi(\rho_1, \rho_2) & := p(\rho_1) - p(\rho_2) - p'(\rho_2) (\rho_1 - \rho_2) \\
			& = \sigma(\rho_1,\rho_2) \abs{\rho_1-\rho_2}^2,
		\end{aligned}
	\end{equation}
where $\sigma(\rho_1,\rho_2) := \int_0^1 \int_0^1 p''(\rho_2 + r_1 r_2 (\rho_1-\rho_2)) dr_2 r_1 dr_1.$

\item We shall use the conventions
\begin{equation}\label{rel-leq}
    A \lesssim B, \quad A \gtrsim B, \quad A\sim B \andd A = O(1) B,
\end{equation}
which mean that
\begin{equation}
    A \leq CB, \quad A \geq C^{-1} B, \quad C^{-1} B \leq A \leq CB \andd \abs{A} \leq C \abs{B},
\end{equation}
respectively, where $ C \geq 1 $ is a constant of the form $$ C = C_0(\mu, \lambda, \rhob) (|\uvb|^{k} + M_0^l), $$ with $C_0$ depending on $\mu, \lambda $ and $ \rhob$, and $k$ and $l$ being some given positive integers.

\item For any $ \ p\in [1,+\infty]$ and $s \geq 0,$ we use the notations,  
$$ \norm{\cdot}_{L^p} := \norm{\cdot}_{L^p(\Omega)} \andd \norm{\cdot}_{H^s} := \norm{\cdot}_{H^s(\Omega)}. $$

\end{itemize}


\vspace{0.3cm}

\section{Ansatz and reformulation of problem}\label{Sec-ansatz}

\subsection{An auxiliary flow}
As stated in \cref{Sec-novelty}, we require the smallness of the vorticity of the background flow in the estimate of the linearized system for the perturbations.
We choose an \textit{auxiliary flow}, 
\begin{equation}\label{profile-t}
\uvt^\vs(x_3,t) := \uv^\vs(x_3, t-t_0+\Lambda) = \Theta\Big(\frac{x_3}{\sqrt{t+\Lambda}}\Big)\uvb,
\end{equation}
to replace $\uv^\vs$ in \cref{profile-t0}, where $\Lambda\geq 1$ is a constant to be determined.
In fact, $\Lambda$ will be chosen as
\begin{equation}\label{Lambda}
\Lambda =  \max\{ C_1(\mu,\rhob) \big(\abs{\uvb}^2 + M_0\big),\, 1 \},
\end{equation}
for some suitably large constant $C_1>0$ that depends only on $\mu$ and $\rhob$. 
Owing to the choice \cref{Lambda}, if the following constants depend on $\Lambda$, we still use the conventions in \cref{rel-leq}.  

\textit{Additional initial perturbation.} A cost for replacing the background flow is that the difference, $\uv^{\vs}-\uvt^{\vs}, $ provides an extra large perturbation to the zero mode of the tangential velocity. More precisely, we can rewrite the initial condition \cref{ic} as
\begin{equation}\label{ic-t}
	(\rho_0, \uv_0)(x) = (\rhob, \uvt^\vs)(x_3,t=0) + (\e b_0, \vvt_0)(x),
\end{equation}
where the new perturbation,
\begin{equation}
    \vvt_0(x) := (\uv^{\vs}-\uvt^{\vs})(x_3,t=0) + \vv_0(x) \in H^3_{3/4},
\end{equation}
satisfies that 
\begin{equation}\label{rel-v0}
	\begin{aligned}
		& \vvt_{0\perp}^\od(x_3) = \Big[ \Theta(x_3) - \Theta\Big(\frac{x_3}{\sqrt{\Lambda}}\Big) \Big] \uvb_\perp + \vv_{0\perp}^\od(x_3), \\
		& \vvt_0^\md = \vv^\md_0, \quad \vt_{03} = v_{03}.
	\end{aligned}
\end{equation}
We also define the new initial perturbation of momentum as  
\begin{equation}\label{w-0}
	\begin{aligned}
		\wvt_0(x) & := (\rho_0 \uv_0)(x) - (\rhob \uvt^\vs)(x_3,0) = \rhob \vvt_0(x) + \e b_0(x) \uv_0(x).
	\end{aligned}
\end{equation}
The following lemma shows that the replacement of the background flow does not affect the smallness assumption \cref{chi}.

\begin{Lem}\label{Lem-vw0}
    If $\e (M_0+\abs{\uvb}) \leq 1,$ then it holds that
    \begin{equation}\label{bdd-vwt-0}
	\begin{aligned}
            & \norm{(\vvt_0, \wvt_0)}_{H_{3/2}^3}  \lesssim 1, \\
            & \norm{(\vt_{03}, \wt_{03})^\od}_{H^1_{3/2}} + \norm{(\vvt_0, \wvt_0)^\md}_{H^1} \lesssim \chi.
	\end{aligned}
    \end{equation}
\end{Lem}
\begin{proof}
It follows from \cref{Theta} that
	\begin{equation}\label{diff-Theta}
		\begin{aligned}
			& \abs{\Theta\Big(\frac{x_3}{\sqrt{t+1}}\Big) - \Theta\Big(\frac{x_3}{\sqrt{t+\Lambda}}\Big)} \\
			& \qquad \lesssim \Big(\frac{1}{\sqrt{t+1}} - \frac{1}{\sqrt{t+\Lambda}} \Big)  \int_0^1 \exp \Big\{ -\frac{\rhob}{4\mu} \abs{h (r,t) x_3}^2 \Big\} x_3 dr,
		\end{aligned}
	\end{equation}
	where $ h(r,t) := \frac{1}{\sqrt{t+\Lambda}} + r \big(\frac{1}{\sqrt{t+1}} - \frac{1}{\sqrt{t+\Lambda}} \big). $ Then for any $\alpha\geq 0$, it holds that
    \begin{equation}\label{bdd-vt-0}
	\begin{aligned}
		\norm{\vvt_{0\perp}^\od}_{L^2_{\alpha}} & \leq C(\mu, \rhob)  \abs{\uvb} \Lambda^{\frac{2\alpha+1}{4}} + \norm{ \vv_{0\perp}^\od}_{L^2_{\alpha}}.
	\end{aligned}
    \end{equation}
	Note that $\vvt_0^\md = \vv_0^\md$ and
	\begin{align*}
		 \wvt_0^\md = \rhob \vv_0^\md + \e \big[ b_0^\od \uv_0^\md + b_0^\md \uv_0^\od + \big(b_0^\md \uv^\md_0\big)^\md \big].
	\end{align*}
	It follows from the Sobolev inequality that
	$$ \norm{(b_0, \uv_0)}_{W^{1,\infty}} \lesssim 1. $$
	Then collecting \cref{rel-v0,bdd-vt-0,w-0}, one can get \cref{bdd-vwt-0}.
\end{proof}

\vspace{.2cm}

\subsection{An ansatz and the anti-derivative variables}

To overcome the difficulty in the linear estimates of zero-modes, we use the anti-derivative technique which was initiated by \cite{Ilin1960} and further developed by \cite{Yuan2023s} to the multi-dimensional problems.
The strategy is to construct a suitable ansatz $(\rhot, \mvt)$ for the zero mode $(\rho^\od, \mv^\od), $ such that the perturbation satisfies that $ (\rho^\od-\rhot, \mv^\od-\mvt) = \p_3 (\Phi, \Psi) $ with $(\Phi, \Psi) \in H^s(\R)$ for some $s\geq 1$.
A necessary condition for the existence of $(\Phi,\Psi)$ is the zero-mass property, 
\begin{equation}\label{zero-mass}
    \int_\R (\rho^\od -\rhot, \mv^\od-\mvt)(x_3,t) dx_3 = 0 \quad \text{for all } \ t\geq 0.
\end{equation}
To achieve this, we follow \cite{Liu-vshock,HXXY2023} to introduce some linear diffusion waves, which propagate along each characteristic of the system \cref{NS}, to carry the excessive mass of the perturbation around the background flow, 
$\int_\R (\rho^\od -\rhob, \mv^\od-\rhob \uvt^\vs)(x_3,t) dx_3.$

\textit{Characteristics in the normal direction}.
The ansatz is devised to be a planar wave, that is, it depends only on $(x_3,t).$
Regardless of the tangential derivatives and the viscosity in \cref{NS}, we arrive at the hyperbolic system in the normal direction,
\begin{equation}
	\p_t \Big(\begin{array}{c}
		\rho\\
		\mv\end{array}
	\Big) + A(\rho,\uv) \p_3 \Big(\begin{array}{c}
		\rho\\
		\mv\end{array}\Big) = 0,
\end{equation}
where
\begin{equation}
	A(\rho,\uv) = \left(\begin{array}{cccc}
		0 & 0 & 0 & 1 \\
		-u_1 u_3 & u_3 & 0 & u_1 \\
		-u_2 u_3 & 0 & u_3 & u_2 \\
		\abs{\e^{-1} a(\rho)}^2 - \abs{u_3}^2 & 0 & 0 & 2u_3
	\end{array} \right),
\end{equation}
and $a$ is the sound speed $a(\rho) := \sqrt{p'(\rho)}$. 
The matrix $A(\rho,\uv)$ has four real eigenvalues
\begin{equation}
	\lambda_0 = u_3 - \e^{-1} a(\rho), \quad \lambda_1 = u_3, \quad \lambda_2 = u_3, \quad \lambda_3 = u_3 + \e^{-1} a(\rho),
\end{equation}
with the associated linearly-independent right eigenvectors
\begin{equation}\label{eigen-vec}
	\rv_0(\rho,\uv) = \left(\begin{array}{c}
		1 \\
		u_1 \\
		u_2 \\
		\lambda_0
	\end{array}\right), \quad \rv_1 = \left( \begin{array}{c}
		0 \\
		1 \\
		0\\
		0
	\end{array}\right), \quad \rv_2 = \left( \begin{array}{c}
		0 \\
		0 \\
		1 \\
		0
	\end{array}\right), \quad \rv_3(\rho,\uv) = \left(\begin{array}{c}
		1 \\
		u_1 \\
		u_2 \\
		\lambda_3
	\end{array}\right).
\end{equation}
With the constant states $(\rhob, \pm \uvb)$, we denote 
\begin{equation}
    \begin{aligned}
        & \bar{a}:= a(\rhob), \quad
        \lambda_0^-  := \lambda_0(\rhob, -\uvb) = -\e^{-1} \bar{a}, \quad \lambda_3^+ := \lambda_3(\rhob, \uvb) = \e^{-1} \bar{a}, \\
        & \rv_0^- := \rv_0(\rhob, -\uvb) = (1, -\uvb_\perp, -\e^{-1} \bar{a})^t, \qquad \rv_3^+ := \rv_3(\rhob, \uvb) = (1, \uvb_\perp, \e^{-1} \bar{a})^t.
    \end{aligned}
\end{equation}

\textit{Diffusion waves}.
Define 
\begin{equation}\label{kernel}
	\begin{aligned}
		\vartheta(x_3,t) :=
		\frac{\sqrt{\rhob}}{2 \sqrt{\pi\mu(t+\Lambda)}}
		\exp\Big(- \frac{ \rhob x_3^2}{4\mu (t+\Lambda)}\Big),
	\end{aligned}
\end{equation}
and
\begin{equation}\label{ker-0-3}
	\vartheta_\pm(x_3,t) := \vartheta(x_3 \mp \e^{-1} \bar{a} t, t), 
\end{equation}
which satisfy the transport-diffusion equations
\begin{equation*}
	\begin{aligned}
		\p_t \vartheta = \frac{\mu}{\rhob} \p_3^2 \vartheta, \qquad \p_t \vartheta_\pm \pm \e^{-1} \bar{a} \p_3 \vartheta_\pm = \frac{\mu}{\rhob} \p_3^2 \vartheta_\pm, 
	\end{aligned}
\end{equation*}
and 
\begin{equation*}
	\int_\R \vartheta(x_3,t) dx_3 = \int_\R \vartheta_\pm(x_3,t) dx_3 \equiv 1.
\end{equation*}

\textit{Ansatz.}
The ansatz $(\rhot,\mvt)$ is constructed as the form,
\begin{equation}\label{ansatz}
	\begin{aligned}
		(\rhot,\mvt)(x_3,t) & := (\rhob, \rhob \uvt^\vs)(x_3,t)+  (\alpha_1 \rv_1 + \alpha_2 \rv_2) \vartheta(x_3,t)  \\
		& \qquad + \e \big[\alpha_0  \rv_0^- \vartheta_-(x_3,t) + \alpha_3\rv_3^+  \vartheta_+(x_3,t)\big], \\
		\uvt &:= \mvt/\rhot,
	\end{aligned}
\end{equation}
where $\alpha_i \in \R$ for $i=0,\cdots,3$ are some constants to be chosen such that the zero-mass condition \cref{zero-mass} holds true. Equivalent to \cref{ansatz}, the ansatz satisfies that
\begin{equation}\label{ansatz-1}
	\begin{aligned}
		& \rhot = \rhob + \e \big(\alpha_0 \vartheta_- + \alpha_3 \vartheta_+\big), \quad \mt_3 = \bar{a} \big(  \alpha_3 \vartheta_+-\alpha_0 \vartheta_- \big). \\
		& \mt_i = \rhob\ut_i^\vs + \alpha_i \vartheta + \e \ub_i (\alpha_3 \vartheta_+ - \alpha_0 \vartheta_-), \qquad i=1,2.
	\end{aligned}
\end{equation}
Plugging \cref{ansatz-1} into \cref{NS} yields that
\begin{equation}\label{equ-ansatz}
	\begin{cases}
		\p_t \rhot + \dv \mvt  = \e \p_3 F_0, \\
		\p_t \mvt + \dv \big(\uvt \otimes \mvt\big) + \frac{1}{\e^2}\nabla p(\rhot) = \mu \lap \uvt + (\mu+\lambda) \nabla\dv \uvt + \p_3 \Fv,
	\end{cases}
\end{equation}
where the error terms, $ F_0=F_0(x_3,t)\in\R $ and $ \Fv = (F_1, F_2, F_3)(x_3,t) \in \R^3, $ are given by
\begin{equation}\label{F}
	\begin{aligned}
		F_0 & = \frac{\mu}{\rhob}  \big(\alpha_0 \p_3 \vartheta_- + \alpha_3 \p_3 \vartheta_+\big),\\
		F_i & = \mu \Big[ \frac{ \alpha_i}{\rhob} \p_3 \vartheta - \p_3 \big(\ut_i - \ut_i^\vs \big) \Big] + \Big[ \frac{\mt_3 \mt_i}{\rhot} - \bar{a} \ub_i \big( \alpha_0 \vartheta_- + \alpha_3 \vartheta_+ \big) \Big] \\
		& \quad + \frac{\e \mu \ub_i }{\rhob} \big( \alpha_3 \p_3 \vartheta_+ - \alpha_0 \p_3 \vartheta_- \big), \qquad i=1,2, \\
		F_3 & = \frac{\mt_3^2}{\rhot} + \frac{1}{\e^2} \varpi(\rhot,\rhob) + \frac{\mu\bar{a}}{\rhob}\big( \alpha_3 \p_3 \vartheta_+ - \alpha_0 \p_3 \vartheta_- \big) - \mut \p_3 \ut_3.
	\end{aligned}
\end{equation}
Recall here that $\varpi$ is defined by \cref{varpi} and $\mut:=2\mu+\lambda$. 
Owing to the conservative forms of \cref{NS,equ-ansatz}, it suffices to fulfill \cref{zero-mass} at $t=0.$ 
This determines the unique choice of $\{\alpha_i: i=0,\cdots,3\}$, satisfying that 
\begin{equation}\label{alpha-i}
	\int_\R (\e b_0, \wvt_0)^\od dx_3 = \alpha_1 \rv_1 + \alpha_2 \rv_2 + \e \big(\alpha_0 \rv_0^- + \alpha_3 \rv_3^+\big),
\end{equation}
which gives
\begin{equation}\label{alpha}
	\begin{aligned}
		& \alpha_0 = \frac{1}{2} \int_\R \big( b_0 - \bar{a}^{-1} \wt_{03} \big)^\od dx_3, \qquad\quad \alpha_3 = \frac{1}{2} \int_\R \big( b_0 + \bar{a}^{-1} \wt_{03} \big)^\od dx_3, \\ 
		& \alpha_i = \int_\R \big(\wt_{0i} - \e \bar{a}^{-1} \ub_i \wt_{03} \big)^\od dx_3 \qquad \text{for }\, i=1,2.
	\end{aligned}
\end{equation}
One can apply \cref{Lem-vw0} to obtain the bounds of the constants in \cref{alpha} directly.

\begin{Lem}\label{Lem-alpha}
    If $\e (M_0+\abs{\uvb}) \leq 1,$ then the constant coefficients satisfy that
    \begin{equation}\label{small-alpha}
		\abs{\alpha_0} + \abs{\alpha_3} \lesssim \chi, \qquad \abs{\alpha_1} + \abs{\alpha_2} \lesssim 1.
    \end{equation}
\end{Lem}


Using \cref{Lem-alpha}, one can also get the point-wise estimates of the errors associated with the ansatz.

\begin{Lem}\label{Lem-F}
	For $ j=0,1,2, \cdots, $ the difference between the ansatz \cref{ansatz-1} and the background vortex layer $(\rhob, \uvt^\vs)$ satisfies that
	\begin{equation}\label{bdd-diff}
		\begin{aligned}
			\e^{-1} \abs{\p_3^j (\rhot-\rhob)} + \abs{\p_3^j \mt_3}  + \abs{\p_3^j \ut_3} &  \lesssim \chi \, (t+\Lambda)^{-\frac{1+j}{2}} \Rc, \\
			\abs{\p_3^j (\mvt_\perp - \rhob \uvt_\perp^\vs)} + \abs{\p_3^j (\uvt_\perp - \uvt_\perp^\vs)} & \lesssim (t+\Lambda)^{-\frac{1+j}{2}} \Rc, \\
			\abs{\p_3^{j+1} \mvt_\perp} + \abs{\p_3^{j+1} \uvt_\perp} & \lesssim (t+\Lambda)^{-\frac{1+j}{2}} \Rc,
		\end{aligned}
	\end{equation}
	and the error terms in \cref{equ-ansatz} satisfy that
	\begin{equation}\label{bdd-F}
		\begin{aligned}
			\abs{\p_3^j F_0} + \abs{\p_3^j \Fv} & \lesssim \chi \, (t+\Lambda)^{-\frac{2+j}{2}} \Rc.
		\end{aligned}
	\end{equation}
	Here $\Rc$ denotes that
	\begin{equation*}
		\Rc = e^{- \frac{c \lvert x_3\rvert^2}{t+\Lambda} } + e^{- \frac{c \lvert x_3-\lambda_0^-(t+\Lambda)\rvert^2}{t+\Lambda}} + e^{ - \frac{c \lvert x_3-\lambda_3^+(t+\Lambda)\rvert^2}{t+\Lambda}}
	\end{equation*} 
	for some generic constant $ c>0, $ depending only on $\mu$ and $\rhob.$
\end{Lem}
\begin{proof}
We only show the estimate of $\Fv_\perp$ in \cref{bdd-F}, since the others can be derived from \cref{ansatz-1,F} directly.
Let $i\in \{1,2\}$ be fixed. It is noted that
$$
\ut_i - \ut_i^\vs = \frac{1}{\rhob} \big( \mt_i - \rhob \ut_i^\vs \big) + \big(\frac{1}{\rhot} - \frac{1}{\rhob}\big) \mt_i.
$$
Then using \cref{ansatz-1}, the first square bracket in \cref{F}$_2$ satisfies that
\begin{equation}\label{equality-1}
\begin{aligned}
    &  \frac{\alpha_i}{\rhob} \p_3\vartheta - \p_3 \big(\ut_i - \ut_i^\vs \big) \\
    & \qquad = \frac{\e}{\rhob} \ub_i ( \alpha_0 \p_3 \vartheta_- - \alpha_3 \p_3 \vartheta_+) + \frac{\e}{\rhob} (\alpha_0 \p_3 \vartheta_- + \alpha_3 \p_3 \vartheta_+) \ut_i.
\end{aligned}
\end{equation}
Moreover, the corresponding second square bracket satisfies that
\begin{equation}\label{equality-2}
	\begin{aligned}
    & \frac{\mt_3 \mt_i}{\rhot} - \bar{a} \ub_i \big( \alpha_0 \vartheta_- + \alpha_3 \vartheta_+ \big) \\
    & \quad = \bar{a} \alpha_3 \vartheta_+ (\ut_i^\vs - \ub_i) - \bar{a} \alpha_0 \vartheta_-  (\ut_i^\vs + \ub_i) - \frac{\e \bar{a}}{\rhot}  (\alpha_3^2  \vartheta_+^2 - \alpha_0^2 \vartheta_-^2) \ut_i^\vs \\
	& \qquad + \frac{\bar{a}}{\rhot} (\alpha_3 \vartheta_+ - \alpha_0 \vartheta_-) \big[ \alpha_i \vartheta + \e \ub_i (\alpha_3 \vartheta_+ - \alpha_0 \vartheta_-) \big].
\end{aligned}
\end{equation}
To estimate \cref{equality-2}, first note (by \cref{Theta}) that
\begin{align*}
	\abs{\ut_i^\vs \mp \ub_i} \leq C(\mu, \rhob) \abs{\uvb} e^{ -\frac{c_0 x_3^2}{t+\Lambda} } \qquad \text{for } \, \pm x_3 \geq 0.
\end{align*}
For any given positive constants $c_0, c_1,$ if $\e \leq \e_0$ with $\e_0$ suitably small, then there exists $c>0$ that is independent of $\e,$ such that
\begin{align}
	e^{ -\frac{c_0 x_3^2}{t+\Lambda} - \frac{c_1 (x_3 \pm \e^{-1} \bar{a} t)^2}{t+\Lambda}} \leq e^{- c (t+ |x_3|) } \qquad \forall x_3\in \R, t\geq 0.
\end{align}
Then the remaining proof is through direct calculations. 
\end{proof}

\textit{Anti-derivative variables.}
Owing to \cref{zero-mass}, we define the anti-derivative variables,
\begin{equation}\label{def: anti-derivative Phi Psi}
	\quad (\Phi, \Psi)(x_3,t) := \int_{-\infty}^{x_3} (\phi, \psi)^\od(y_3,t) dy_3 = - \int_{x_3}^{+\infty} (\phi, \psi)^\od(y_3,t) dy_3,
\end{equation}
where $(\Phi, \Psi)(x_3,0) = (\Phi_0, \Psi_0)(x_3) := \int_{-\infty}^{x_3} (\phi_0,\psi_0)^\od(y_3) dy_3.$ 

\begin{Lem}\label{Lem-ic}
	Under the assumptions in \cref{Thm}, the initial data in \cref{equ-phipsi}, $(\phi_0,\psi_0), $ belongs to $ H^3_{3/2}$; and the anti-derivative variable, $(\Phi_0, \Psi_0), $ belongs to $ H^4$. Moreover, it holds that
	\begin{equation*}
		\begin{aligned}
			\norm{(\Phi_0, \Psi_{03})}_{L^2} + \norm{(\phi_0, \psi_{03}, \zeta_{03})^\od}_{H^1} + \norm{(\phi_0, \psi_0, \zeta_0)^\md}_{H^1} & \lesssim \chi,
		\end{aligned}
	\end{equation*}
	and
	\begin{equation}\label{Mt}
		\norm{(\Phi_0, \Psi_0)}_{L^2} + \norm{(\phi_0, \psi_0, \zeta_0)}_{H^3} \lesssim 1,
	\end{equation}
	where $\zeta_0 := \uv_0 - \uvt(x_3,0) = \frac{1}{\rho_0} (\psi_0-\e \uvt(x_3,0) \phi_0). $
\end{Lem}
\begin{proof}
    We show only the estimates for the anti-derivatives. 
    With the use of the Minkowski inequality, it follows from \cref{kernel} that
    \begin{align*}
        \norm{\int_{\pm \infty}^{x_3} \vartheta(y_3,0) dy_3}_{L^2(\R_\pm; dx_3)} & \lesssim \Lambda^{\frac{1}{4}}.
    \end{align*}
    Note that $\Lambda$ is determined by \cref{Lambda}.
    Using \cref{small-alpha,IC-phipsi}, one has that
	\begin{align*}
		\norm{\Phi_0}_{L^2} & = \Big[ \int_{-\infty}^0 \Big(\int_{-\infty}^{x_3} \phi_0^\od(y_3) dy_3\Big)^2 dx_3 + \int_{0}^{+\infty} \Big(\int_{x_3}^{+\infty} \phi_0^\od(y_3) dy_3\Big)^2 dx_3 \Big]^{\frac{1}{2}} \\
		& \lesssim \norm{\langle x_3\rangle^\frac{3}{2} b_0}_{L^2}  \Big[ \int_0^{+\infty} \int_{x_3}^{+\infty} (y_3^2+1)^{-\frac{3}{2}} dy_3  dx_3 \Big]^{\frac{1}{2}} + \big(\abs{\alpha_0} + \abs{\alpha_3}\big) \Lambda^{\frac{1}{4}} \\
		& \lesssim \chi.
	\end{align*}
	Similarly, one can estimate $\Psi_{03}$. 
    Also, it follows from \cref{small-alpha,IC-phipsi,bdd-vwt-0} that
	\begin{align*}
		\norm{\Psi_{0\perp}}_{L^2}
		& \lesssim \norm{\langle x_3 \rangle^\frac{3}{2} \wvt_{0\perp}}_{L^2} + \big(\sum_{i=0}^{3} \abs{\alpha_i}\big) \Lambda^{\frac{1}{4}}
		\lesssim 1.
	\end{align*}
\end{proof} 

\vspace{.2cm}

\subsection{Reformulation of the Problem}

With the ansatz in \cref{ansatz}, we define the perturbations,
\begin{equation}\label{phipsi}
	\phi := \e^{-1} \big(\rho - \rhot\big), \quad \psi = (\psi_1, \psi_2, \psi_3) := \mv-\mvt,
\end{equation}
and
\begin{equation}\label{zeta}
    \zeta = (\zeta_1, \zeta_2, \zeta_3) := \uv - \uvt = \frac{1}{\rho} (\psi - \e \uvt \phi).
\end{equation}
Then using \cref{NS,equ-ansatz}, we arrive at the reformulation of the problem,
\begin{equation}\label{equ-phipsi}
	\begin{cases}
		\p_t \phi + \e^{-1} \dv \psi = - \p_3 F_0, \\
		\p_t \psi + \dv \big(\uv\otimes\mv - \uvt\otimes\mvt \big) + \e^{-1} \nabla \big( p'(\rhot) \phi \big) + \e^{-2} \nabla \big(\varpi(\rho,\rhot)\big)  \\
		\qquad - \mu \lap \zeta - (\mu+\lambda) \nabla \dv \zeta = - \p_3 \Fv, \\
		(\phi, \psi)(x,0) = (\phi_0, \psi_0)(x) := (\e^{-1}(\rho - \rhot), \mv-\mvt)(x,0),
	\end{cases}
\end{equation}
where $\varpi$ is defined by \cref{varpi}.
With \cref{ic,ansatz}, the initial data satisfies that 
\begin{equation*}
	\begin{aligned}
		(\e \phi_0,\psi_0)(x) & = (\e b_0, \wvt_0)(x) - (\alpha_1 \rv_1 + \alpha_2 \rv_2) \vartheta(x_3,0) \\
		& \quad - \e (\alpha_0 \rv_0^- \vartheta_- + \alpha_3 \rv_3^+ \vartheta_+)(x_3,0),
	\end{aligned}
\end{equation*}
which is equivalent to that 
\begin{equation}\label{IC-phipsi}
	\begin{aligned}
		\phi_0(x) & = b_0(x) - \big(\alpha_0 \vartheta_- + \alpha_3 \vartheta_+\big)(x_3,0), \\
		\psi_{03}(x) & = \wt_{03}(x) + \bar{a} \big( \alpha_0 \vartheta_- - \alpha_3 \vartheta_+ \big)(x_3,0), \\
		\psi_{0i}(x) &  = \wt_{0i}(x) - \alpha_i \vartheta(x_3,0) + \e \ub_i (\alpha_0 \vartheta_- - \alpha_3 \vartheta_+)(x_3,0) \quad \text{for } \ i=1,2.
	\end{aligned}
\end{equation}



\begin{Thm}\label{Thm-pert}
	Under the assumptions of \cref{Thm}, the Cauchy problem \cref{equ-phipsi} admits a unique strong solution $(\phi,\zeta)$ globally in time, satisfying that
    \begin{align*}
        \sup_{t>0} \norm{(\phi, \zeta)}_{H^3}^2 + \int_0^{+\infty} \Big(\norm{\nabla \phi}_{H^2}^2 + \norm{\nabla \zeta}_{H^3}^2\Big) dt \lesssim 1,
    \end{align*}
    and
	\begin{equation}\label{behavior-inf}
		\begin{aligned}
			\norm{(\phi, \zeta)(\cdot,t)}_{L^\infty} & \lesssim (t+1)^{-\frac{1}{2}} \qquad \forall t\geq 0.
		\end{aligned}
	\end{equation}
\end{Thm}
We show that the main result of this paper, \cref{Thm}, can follow from \cref{Thm-pert}.

\subsection{Proof of the main result}
Recall the notations \cref{pert-bv,phipsi,zeta}, which satisfy
\begin{equation}\label{rel-pert}
	\begin{aligned}
	b^\e & = \e^{-1} \big( \rhot - \rhob \big) + \phi, \\
	\vv^\e & = \big( \uvt^\vs - \uv^\vs \big) + \big( \uvt - \uvt^\vs \big) + \zeta.
	\end{aligned}
\end{equation}
Assume that \cref{Thm-pert} holds true, then it suffices to show \cref{Thm-L2,rate-Thm} to complete the proof of \cref{Thm}. In fact, it follows from \cref{bdd-diff} that for any $p\in [1,+\infty],$
\begin{equation}\label{Lp-diff1}
	\e^{-1} \norm{\p_3^j (\rhot -\rhob)}_{L^p} + \norm{\p_3^j (\uvt-\uvt^\vs)}_{L^p} \lesssim (t+\Lambda)^{-\frac{j+1}{2} + \frac{1}{2p}}, \quad j=0,1,2,\cdots.
\end{equation}
Besides, it follows from \cref{diff-Theta} that
\begin{align}
	& \norm{\uvt^\vs - \uv^\vs}_{L^p} \notag \\
		& \quad \lesssim \norm{\Theta\Big(\frac{x_3} {\sqrt{t+1}}\Big) - \Theta\Big(\frac{x_3}{\sqrt{t+\Lambda}}\Big)}_{L^p} \notag \\
		&\quad \lesssim \begin{cases}
		(t+\Lambda)^{\frac{1}{p}} - (t+1)^{\frac{1}{p}} \lesssim (\Lambda-1) (t+1)^{- 1 + \frac{1}{p}}, \quad p \in [1,+\infty), \\
		\ln\big(1+\frac{\Lambda}{t}\big) - \ln\big(1+\frac{1}{t}\big) \lesssim (\Lambda-1)\ t^{-1}, \quad p=+\infty.
		\end{cases} \label{Theta-Lp}
\end{align}
Then collecting \cref{behavior-inf,rel-pert,Lp-diff1,Theta-Lp}, one has that
	\begin{align*}
		\norm{(b^\e, \vv^\e)}_{L^2} \lesssim 1 \andd \norm{(b^\e, \vv^\e)}_{L^\infty} 
		\lesssim (t+1)^{-\frac{1}{2}}.
	\end{align*}
	The derivatives in \cref{Thm-L2} can be estimated similarly, and the remaining proof is omitted.

\vspace{.1cm}

\subsection{Useful lemmas}

At the end of this section, we list two useful lemmas.

\begin{Lem}[\cite{HY2020}, Theorem 1.4 \& Lemma 3.3]\label{Lem-GN}
	Assume that $ f(x) \in L^\infty(\Omega) $ is periodic in $ \xp = (x_1, x_2) $. Then there exists a decomposition $ f(x) = f^{(1)}(x_3) + \sum\limits_{k=2}^{3} f^{(k)}(x) $ such that 
	\begin{itemize}
		\item[i)] for the zero and non-zero modes defined in \cref{od,md}, it holds that
		\begin{equation}\label{decomp}
			f^{(1)} \equiv f^{\od}, \qquad f^{(2)} + f^{(3)} \equiv f^\md,
		\end{equation}
		
		\item[ii)] if $\nabla^l f $ belongs to $ L^p$ with an order $l\geq 0 $ and $p\in[1,+\infty]$, then each $f^{(k)}$ satisfies that
		\begin{equation}\label{ineq-cauchy}
			\norm{\nabla^l f^{(k)}}_{L^p(\Omega)} \leq C(p) \norm{\nabla^l f}_{L^p(\Omega)};
		\end{equation} 
		
		\item [iii)] each $ f^{(k)} $ satisfies the classical k-dimensional Gagliardo--Nirenberg inequality, that is,
		\begin{equation}\label{G-N-type-1}
			\norm{\nabla^j f^{(k)}}_{L^p(\Omega)} \leq C(p,r,q) \norm{\nabla^{m} f^{(k)}}_{L^{r}(\Omega)}^{\theta_k} \norm{f^{(k)}}_{L^{q}(\Omega)}^{1-\theta_k},
		\end{equation}
		where $ 0\leq j< m $ is any integer, $ 1\leq p \leq +\infty $ is any number and $\theta_k \in [\frac{j}{m}, 1)$ satisfies 
		$$ \frac{1}{p} = \frac{j}{k} + \Big(\frac{1}{r}-\frac{m}{k}\Big) \theta_k + \frac{1}{q}\big(1-\theta_k\big). $$
	\end{itemize}
\end{Lem}


Using the Poincar\'{e} inequality on the transverse domain $\Torus^2$, one can get that 

\begin{Lem}\label{Lem-poin}
	Suppose that $ f(x) $ belongs to $ W^{1,p}(\Omega) $ with $ p\in [1,+\infty). $ Then its non-zero mode $ f^\md $ satisfies that
	\begin{equation}\label{poincare}
		\norm{f^\md}_{L^p(\Omega)} \leq C(p) \norm{\nabla_{\xp} f^\md}_{L^p(\Omega)} = C(p) \norm{\nabla f}_{L^p(\Omega)},
	\end{equation}
	where the constant $C(p)>0$ depends only on $ p. $
\end{Lem}

\vspace{0.2cm}

\section{Local existence and a priori estimates}\label{Sec-apriori}

The proof of \cref{Thm-pert} is based on a bootstrap procedure using the local existence and the a priori estimates.

\subsection{Local existence}
Now we show that local existence for the problem \cref{equ-phipsi}. 
For any $ T>0, $ we define the solution space for \cref{equ-phipsi} as
\begin{align*}
	\mathbb{B}(0,T) := \Big\{ (\phi, \zeta): 
	& \ (\phi,\zeta) \in C(0,T; H^3(\Omega)), \\ 
    & \ (\Phi, \Psi) \in C(0,T; L^2(\R)), \\
	& \ \nabla\phi \in L^2(0,T; H^2(\Omega)), \ \nabla\zeta \in L^2(0,T; H^3(\Omega)) \Big\}.
\end{align*}

\begin{Thm}[Local existence]\label{Thm-local}
	Suppose that the hypotheses of \cref{Thm} hold true, and the initial data $ (\phi_0,\psi_0) $ satisfies that
	\begin{equation}\label{ic-loc}
		(\phi_0, \psi_0)\in {H^3(\Omega)}, \quad \inf_{x\in\Omega} \rho(x,0) := \underline{\rho} > 0.
	\end{equation}
	Then for any given $M_0 > 0,$ there exists a positive constant $T_0$ that depends on $M_0,$ such that if $\norm{(\phi_0, \psi_0)}_{H^3}\leq M_0,$ 
	the Cauchy problem \cref{equ-phipsi} has a unique solution $ (\phi, \zeta) \in \mathbb{B}(0,T_0) $, which satisfies that
	$\rho(x,t) \geq \frac{1}{2} \underline{\rho}>0,$ and
	\begin{equation*}
		\sup_{t\in [0,T_0]} \norm{(\phi, \zeta)}_{H^3}^2(t) + c_0 \int_0^{T_0} \big(\norm{\nabla\phi}_{H^2}^2 + \norm{\nabla \zeta}_{H^3}^2\big) d\tau \leq 2 \norm{(\phi_0,\zeta_0)}_{H^3}^2,
	\end{equation*}
	where $ c_0>0 $ is a constant, depending only on $\mu$ and $\underline{\rho}. $ If there holds in addition to \cref{ic-loc} that 
	$$ (\Phi_0, \Psi_0)(x_3):= \int_{-\infty}^{x_3} (\phi_0, \psi_0)^\od(y_3) dy_3 \in L^2(\R), $$
    then the anti-derivative variables, 
	$$ (\Phi, \Psi)(x_3,t) := \int_{-\infty}^{x_3} (\phi, \psi)^\od(y_3,t) dy_3, $$ 
	exist in $C(0,T_0; L^2(\R)), $ and satisfies that
    \begin{equation*}
        \sup_{t\in [0,T_0]} \norm{(\Phi, \Psi)}_{L^2}^2(t) \leq 2 \norm{(\Phi_0, \Psi_0)}_{L^2}^2 + 2 \norm{(\phi_0,\zeta_0)}_{H^3}^2.
    \end{equation*}
\end{Thm}
The local existence of $(\phi, \zeta)$ can be derived from a standard argument; see \cite{MN1980} for instance. The existence of the anti-derivatives, $(\Phi, \Psi), $ can follow from the standard parabolic theory, and we refer to \cite[Section 6]
{Yuan2023n} for the detailed proof.

\vspace{.1cm}

\subsection{A priori estimates}

Before stating the a priori estimates for the problem \cref{equ-phipsi}, we first introduce \textit{two effective variables} associated with the momentum.
\begin{itemize}
    \item As in \cite{HXXY2023}, we first introduce an effective variable in terms of the anti-derivatives,
	\begin{equation}\label{Zv}
		\Zv := \Psi - \e \uvt \Phi.
	\end{equation}
    The replacement of $\Psi$ by \cref{Zv} can help overcome the difficulty due to the large amplitude of the vortex sheet, $\abs{\uvb}.$
	
    \item Besides \cref{Zv}, in this paper we introduce a new variable in terms of the original perturbations,
	\begin{equation}\label{wv}
		\begin{aligned}
		  \wv & := \psi - \e \uvt \phi  \\
            & = \rho \zeta = \rhot \zeta + \e \phi \zeta.
		\end{aligned}
	\end{equation}
	The role of \cref{wv} plays an effective connection between the anti-derivatives, $(\Phi, \Zv),$ and the non-zero modes, $(\phi^\md, \zeta^\md)$, based on the facts
    \begin{equation}\label{wv-od}
        \begin{aligned}
            \wv^\od & = \p_3 \Zv + \e \p_3 \uvt \Phi \\
            & =  \rhot \zeta^\od + \e \phi^\od \zeta^\od + \e \big(\phi^\md \zeta^\md\big)^\od, 
        \end{aligned}
    \end{equation}
    and
    \begin{equation}\label{wv-md}
        \begin{aligned}
            \wv^\md & = \psi^\md - \e\uvt \phi^\md \\
        & = \rhot \zeta^\md + \e \big[ \phi^\od \zeta^\md + \phi^\md \zeta^\od + \big( \phi^\md \zeta^\md \big)^\md \big]. 
        \end{aligned}
    \end{equation}
    In fact, in the non-zero estimates (see \cref{Sec-md}), these relations assist in the estimates of
    the complex interactions of the zero modes and non-zero modes, especially for the ones arising from nonlinear convection, $\uv\otimes\uv.$

\end{itemize}


\vspace{.1cm}

For a fixed $ T>0 $, suppose that $ (\phi, \zeta) \in \mathbb{B}(0,T) $ is a solution to the problem \cref{equ-phipsi}.
Then it holds that
\begin{align*}
	& (\psi,\wv) \in C(0,T; H^3(\Omega)), \quad \nabla (\psi,\wv) \in L^2(0,T; H^3(\Omega)); \\
	& (\Phi, \Psi, \Zv) \in C(0,T; H^4(\R)).
\end{align*}
For $0\leq t\leq T, $ we define two weighted energy functionals,
\begin{equation}\label{Ec*}
	\Ec^*(t):= \sum_{j=0}^{2} (t+1)^j \norm{ \p_3^j \big(\Phi, Z_3\big)}_{L^2(\R)}^2+ \sum_{j=0}^{1} (t+1)^j \norm{\nabla^j \big(\phi^\md, \wv^\md\big)}_{L^2(\Omega)}^2,
\end{equation}
and
\begin{equation}\label{Ec}
	\begin{aligned}
		\Ec(t):= \Ec^*(t) + \sum_{j=0}^{2} (t+1)^{j} \norm{ \p_3^j \Zv_\perp}_{L^2(\R)}^2 + (t+1)^2 \norm{\nabla^2 (\phi,\wv)}_{H^1(\Omega)}^2.
	\end{aligned}
\end{equation}
It follows from \cref{Lem-ic} that
\begin{equation}\label{ic-Ec}
	\begin{aligned}
		\Ec^*(0) \lesssim \chi^2 << 1 \andd \Ec(0) \lesssim 1.
	\end{aligned}
\end{equation}
In the a priori estimates, only the components appearing in $\Ec^*(t)$ are assumed to be small. 

Similar to \cite{HXXY2023}, we will see that due to the nontrivial background flow, the energy functionals \cref{Ec,Ec*} grow at the rate $(t+1)^{1/2}.$ Thus, we define two constants
\begin{equation}\label{M&nu}
	\begin{aligned}
		\nu & := \Big(\sup_{ t\in(0,T)} (t+1)^{-\frac{1}{2}} \Ec^*(t)\Big)^{\frac{1}{2}}, \\
		M & := \max\Big\{1, \ \Big(\sup_{ t\in(0,T)} (t+1)^{-\frac{1}{2}} \Ec(t)\Big)^{\frac{1}{2}} \Big\}.
	\end{aligned}
\end{equation}

Now we are ready to state the a priori estimates. 

\begin{Prop}[A priori estimates]\label{Prop-apriori}
	Under the hypotheses of \cref{Thm}, suppose that for a fixed $T>0,$ $ ( \phi, \zeta) \in \mathbb{B}(0,T) $ is a solution to the problem \cref{equ-phipsi} on $\Omega\times[0,T].$
	Then there exist small positive constants $ \e_0, \chi_0 $ and $ \nu_0 $, independent of $t$ and $M,$ such that if 	
	\begin{equation}\label{assum}
		\e \leq \e_0, \quad \chi \leq \chi_0, \quad \nu \leq \nu_0,
	\end{equation}
	and 
	\begin{equation}\label{assum-M}
		(\e + \chi + \nu) M^{k_0} \leq 1,
	\end{equation}
	for some integer $k_0\geq 1$ that depends only on the dimension,
	then it holds that
	\begin{equation}\label{est-apriori}
		\begin{aligned}
			\sup_{t\in(0,T)} (t+1)^{-\frac{1}{2}} \Ec(t) & \lesssim \Ec(0), \\
			\sup_{t\in(0,T)} (t+1)^{-\frac{1}{2}} \Ec^*(t) & \lesssim \Ec^*(0) + \e \Ec(0).
		\end{aligned}
	\end{equation}
\end{Prop} 
The proof of \cref{Prop-apriori} is shown in \cref{Sec-pf-apriori}.

\subsection{List of some a priori bounds}\label{Sec-bound}
For readers' convenience, we list the a priori $L^2$- and $L^\infty$-bounds in terms of the constants in \cref{M&nu} for the following variables. 
The proofs of these a priori bounds are included in the appendix.
\vspace{.1cm}
~\\
\textit{A priori} $L^2$-bounds.
\begin{itemize}
	\item Anti-derivatives: 
	\begin{equation}\label{aL2-anti}
		\begin{aligned}
			\norm{\p_3^j (\Phi, \Psi_3, Z_3)}_{L^2} & \lesssim \nu (t+1)^{-\frac{j}{2} + \frac{1}{4}},\qquad j=0,1,2, \\
			\norm{\p_3^j (\Psi_\perp, \Zv_\perp)}_{L^2} & \lesssim M (t+1)^{-\frac{j}{2} + \frac{1}{4}}, \qquad j=0,1,2, \\
			\norm{\p_3^3 (\Phi, \Psi, \Zv)}_{H^1} & \lesssim M (t+1)^{-\frac{3}{4}}.
		\end{aligned} 
	\end{equation}
	
	\item Zero modes:
	\begin{equation}\label{aL2-od}
		\begin{aligned}
			\norm{\p_3^j (\phi^\od, \psi_3^\od, w_3^\od)}_{L^2} & \lesssim \nu (t+1)^{-\frac{j}{2} - \frac{1}{4}},  \qquad j=0,1, \\
			\norm{\p_3^j (\psi_\perp^\od, \wv_\perp^\od)}_{L^2} & \lesssim M (t+1)^{-\frac{j}{2} - \frac{1}{4}}, \qquad j=0,1, \\
			\norm{\p_3^2 (\phi^\od, \psi^\od, \wv^\od)}_{H^1} & \lesssim M (t+1)^{-\frac{3}{4}}.
		\end{aligned}
	\end{equation}
	
	\item Non-zero modes: 
	\begin{equation}\label{aL2-md}
		\begin{aligned}
			\norm{(\phi^\md, \psi^\md, \wv^\md)}_{H^1} & \lesssim \min\big\{ \nu (t+1)^{-\frac{1}{4}}, \ M (t+1)^{-\frac{3}{4}} \big\}, \\
			\norm{\nabla^2 (\phi^\md, \psi^\md, \wv^\md)}_{H^1} & \lesssim M (t+1)^{-\frac{3}{4}}.
		\end{aligned}
	\end{equation}
	
	\item Original perturbations:
	\begin{equation}\label{aL2-pert}
		\begin{aligned}
			\norm{(\phi, \psi_3, w_3)}_{L^2} & \lesssim  \nu (t+1)^{-\frac{1}{4}}, \\
			\norm{\nabla (\phi, \psi_3, w_3)}_{L^2} & \lesssim \min\big\{ \nu (t+1)^{-\frac{1}{4}}, \ M (t+1)^{-\frac{3}{4}} \big\}, \\
			\norm{\nabla^j (\psi_\perp, \wv_\perp)}_{L^2} & \lesssim M (t+1)^{-\frac{2j+1}{4}}, \qquad j=0,1, \\
			\norm{\nabla^2 (\phi, \psi, \wv)}_{H^1} & \lesssim M (t+1)^{-\frac{3}{4}}.
		\end{aligned}
	\end{equation}
	
\end{itemize}
~\\
\textit{A priori} $L^\infty$-bounds.
\begin{itemize}
	\item Anti-derivatives: 
	\begin{equation}\label{aLinf-anti}
		\begin{aligned}
			\norm{\p_3^j (\Phi, \Psi_3, Z_3)}_{L^\infty} & \lesssim \nu (t+1)^{-\frac{j}{2}}, \qquad j=0,1, \\
			\norm{\p_3^2 (\Phi, \Psi_3, Z_3)}_{W^{1, \infty}} &  \lesssim M^{\frac{3}{4}}\nu^{\frac{1}{4}} (t+1)^{-\frac{3}{4}}, \\
			\norm{\p_3^j (\Psi_\perp, \Zv_\perp)}_{L^\infty} & \lesssim M (t+1)^{-\frac{j}{2}}, \qquad j=0,1, \\
			\norm{\p_3^2 (\Psi_\perp, \Zv_\perp)}_{W^{1, \infty}} & \lesssim M (t+1)^{-\frac{3}{4}}.
		\end{aligned}
	\end{equation}
	
	\item Zero modes:
	\begin{equation}\label{aLinf-od}
		\begin{aligned}
			\norm{(\phi^\od, \psi_3^\od, w_3^\od)}_{L^\infty} & \lesssim \nu (t+1)^{-\frac{1}{2}}, \\
			\norm{\p_3 (\phi^\od, \psi_3^\od, w_3^\od)}_{W^{1, \infty}} & \lesssim M^{\frac{3}{4}}\nu^{\frac{1}{4}} (t+1)^{-\frac{3}{4}}, \\
			\norm{(\psi_\perp^\od, \wv_\perp^\od)}_{L^\infty} & \lesssim M (t+1)^{-\frac{1}{2}}, \\ 
			\norm{\p_3 (\psi_\perp^\od, \wv_\perp^\od)}_{W^{1, \infty}} & \lesssim M (t+1)^{-\frac{3}{4}}.
		\end{aligned}
	\end{equation}
	
	\item Non-zero modes: 
	\begin{equation}\label{aLinf-md}
		\begin{aligned}
			\norm{(\phi^\md, \psi^\md, \wv^\md)}_{W^{1,\infty}}
			& \lesssim \min\big\{ M^{\frac{3}{4}}\nu^{\frac{1}{4}} (t+1)^{-\frac{1}{2}}, M (t+1)^{-\frac{3}{4}} \big\}, 
		\end{aligned}
	\end{equation}
	
	\item original perturbations:
	\begin{equation}\label{aLinf-pert}
		\begin{aligned}
			\norm{(\phi, \psi_3, w_3)}_{L^\infty} & \lesssim M^{\frac{3}{4}}\nu^{\frac{1}{4}} (t+1)^{-\frac{1}{2}}, \\
			\norm{\nabla (\phi, \psi_3, w_3)}_{L^\infty} & \lesssim \min \big\{ M^{\frac{3}{4}}\nu^{\frac{1}{4}} (t+1)^{-\frac{1}{2}}, M (t+1)^{-\frac{3}{4}} \big\}, \\
			\norm{\nabla^j (\psi_\perp, \wv_\perp)}_{L^\infty} & \lesssim M (t+1)^{-\frac{j+2}{4}}, \qquad j=0,1.
		\end{aligned}
	\end{equation}
	
\end{itemize}

Moreover, using \cref{Lem-F,aLinf-pert}, if $\e>0$ is small and $\e^{\frac{1}{2}} M \leq 1 $, then 
\begin{equation}\label{bdd-rho}
	\begin{aligned}
		\frac{\rhob}{2} \leq \rhot(x,t) \leq 2 \rhob \andd
		\frac{\rhob}{4}  \leq \rho(x,t)  \leq 4 \rhob \qquad \forall x\in\Omega, \ t\in [0,T],
	\end{aligned}
\end{equation}
and
\begin{equation}\label{bdd-u}
	\sup_{t\in [0,T]}\norm{\uv}_{L^\infty} \lesssim \abs{\uvb} + \sup_{t\in [0,T]} \norm{\wv}_{L^\infty} \lesssim M.
\end{equation}

\begin{Lem}\label{Lem-rel}
    Under the assumptions of \cref{Prop-apriori}, if $\e^{\frac{1}{2}} M \leq 1,$ then the perturbation of the velocity, $\zeta,$ satisfies that
	\begin{itemize}
		\item for $i=1,2,3$ and $j=0,1,2,$
		\begin{equation}\label{rel-0-2}
		  \pm \norm{\p_3^j \zeta_i^\od}_{L^2} \lesssim \pm \big(\norm{\p_3^{j+1} Z_i}_{L^2}  + \norm{\nabla^j w_i^\md}_{L^2}\big) + \e^{\frac{1}{2}} M (t+1)^{-\frac{2j+1}{4}};
		\end{equation}
		
		\item for $j=0,1,2,$
		\begin{equation}\label{rel-zeta-w}
			\pm \norm{\nabla^j \zeta^\md}_{L^2}  \lesssim \pm \norm{\nabla^j \wv^\md}_{L^2} + \e^{\frac{1}{2}} \norm{\nabla^j \phi^\md}_{L^2};
		\end{equation}
		
		\item 
		\begin{equation}\label{rel-zw3}
			\pm \norm{\nabla^3 \zeta}_{L^2} \lesssim \pm \norm{\nabla^3 \wv}_{L^2} + \e^{\frac{1}{2}} M (t+1)^{-\frac{5}{4}}.
		\end{equation}
	\end{itemize}
\end{Lem}
The proof of \cref{Lem-rel} is placed in \cref{App-rel}.

\vspace{.3cm}



\section{Proof of a priori estimates}\label{Sec-pf-apriori}

In this section, we outline the main steps of the a priori estimates and finish the proof of \cref{Prop-apriori}. 
The proof of the each step is placed in
\cref{Sec-od,Sec-md,Sec-orn} separately.

\subsection{Main steps of a priori estimates}\label{Sec-steps}
~\\
\textit{Step 1.} $H^2$-estimates for the anti-derivative variables with smallness: 
~\\
If $ \big( \nu^{\frac{1}{2}} + \e^{\frac{1}{2}}\big) M \leq 1, $ then
\begin{equation}\label{ineq-od*}
	\begin{aligned}
		& \frac{d}{dt} \Ec^*_\od(t) + \Nc_{\od}^*(t) \\
		& \qquad \lesssim (\chi \nu + \nu^3 + \e M^2) (t+1)^{-\frac{1}{2}} + \nu^{\frac{1}{2}} (t+1) \norm{\nabla^2 (\phi, \wv^\md)}_{L^2}^2,
	\end{aligned}
\end{equation}
where $\Ec^*_\od$ and $ \Nc_\od^* $ are two energy functionals, satisfying that
\begin{equation}\label{ED-od*}
	\begin{aligned}
		\Ec_{\od}^*(t) & \sim \sum_{i=0}^{2} (t+1)^{i} \norm{\p_3^{i} (\Phi,  Z_3)}_{L^2}^2, \\
		\Nc_{\od}^*(t) & \sim \sum_{i=1}^{2} (t+1)^{i-1} \norm{\p_3^i (\Phi, Z_3)}_{L^2}^2 + (t+1)^2 \norm{\p_3^3 Z_3}_{L^2}^2.
	\end{aligned}
\end{equation}
~\\
\textit{Step 2.} 
$H^2$-estimates for full anti-derivative variables without smallness: 
~\\
If $ \big( \nu^{\frac{1}{2}} + \e^{\frac{1}{2}}\big) M \leq 1, $ then
\begin{equation}\label{ineq-od}
	\frac{d}{dt} \Ec_\od(t) + \Nc_\od(t) \lesssim (\chi + \nu^{\frac{1}{2}} + \e) M^2 (t+1)^{-\frac{1}{2}},
\end{equation}
where $ \Ec_\od $ and $ \Nc_\od $ are two energy functionals, satisfying that
\begin{equation}\label{ED-od}
	\begin{aligned}
		\Ec_\od(t) & \sim \sum_{i=0}^{2} (t+1)^{i} \norm{\p_3^{i} (\Phi,  \Zv)}_{L^2}^2, \\
		\Nc_\od(t) & \sim \sum_{i=0}^{1} (t+1)^{i} \norm{\p_3^{i+1} (\Phi, \Zv)}_{L^2}^2 + (t+1)^2 \norm{\p_3^3 \Zv}_{L^2}^2.
	\end{aligned}
\end{equation}
~\\
\textit{Step 3.} $H^1$-estimate for the non-zero modes: ~\\
If $ \big(\chi^{\frac{1}{2}} + \nu^{\frac{1}{6}} + \e^{\frac{1}{4}} \big) M \leq 1, $ then
\begin{equation}\label{ineq-md*}
	\frac{d}{dt} \Ec_\md(t) +\Nc_\md(t) \lesssim (t+1)^{-\frac{5}{4}} \Ec_\md(t) + \e^{\frac{1}{2}} (t+1) \norm{\nabla^2\phi}_{L^2}^2,
\end{equation}
where $ \Ec_\md $  and $ \Nc_\md $ are two energy functionals, satisfying that
\begin{equation}\label{Ec-md}
	\begin{aligned}
		\Ec_\md(t) & \sim \norm{(\phi^\md, \wv^\md)}_{L^2}^2 + (t+1) \norm{\nabla (\phi^\md, \wv^\md)}_{L^2}^2, \\
		\Nc_\md(t) & \sim \norm{\nabla (\phi^\md, \wv^\md)}_{L^2}^2 + (t+1) \norm{\nabla^2 \wv^\md}_{L^2}^2.
	\end{aligned}
\end{equation}
~\\
\textit{Step 4.} Dissipation of the second-order density gradient:   ~\\
If $ \e^{\frac{1}{2}} M \leq 1, $ then
	\begin{equation}\label{est8*}
		\begin{aligned}
			& \frac{d}{dt} \big[(t+1) \Ac_{2,\e\phi}\big] + (t+1) \norm{\nabla^2 \phi}_{L^2}^2 \\
			& \qquad \lesssim (t+1)^{-1} \Nc_\od^* + \Nc_\md + \e (t+1)^{-1} \Nc_\od + \e M^2 (t+1)^{-\frac{3}{2}},
		\end{aligned}
	\end{equation}
	where $\Ac_{2,\e\phi}$ is an energy functional satisfying that
	\begin{equation}\label{bdd-A2}
		\begin{aligned}
			\abs{\Ac_{2,\e\phi}(t)} \lesssim \e M^2 (t+1)^{-\frac{3}{2}}.
		\end{aligned}
	\end{equation}
~\\
\textit{Step 5.} $H^1$-estimate for $\nabla^2 (\phi,\zeta):$ ~\\
If $ \big(\chi^{\frac{1}{2}} + \nu^{\frac{1}{6}} + \e^{\frac{1}{4}} \big) M \leq 1, $ then
	\begin{equation}\label{ineq-ho}
		\begin{aligned}
			& \frac{d}{dt} \Ec_{ho}(t) +\Nc_{ho}(t) \lesssim \Nc_\od(t) + \Nc_\md(t) + (\chi + \e) M^2 (t+1)^{-\frac{1}{2}},
		\end{aligned}
	\end{equation}
	where $ \Ec_{ho} $ and $ \Nc_{ho} $ are two energy functionals, satisfying that
	\begin{equation}\label{Ec-ho}
		\begin{aligned}
			\pm \Ec_{ho}(t) & \lesssim \pm (t+1)^2 \norm{\nabla^2 (\phi, \zeta)}_{H^1}^2 + \e M^2 (t+1)^{\frac{1}{2}}, \\
			\Nc_{ho}(t) & \sim (t+1)^2 \norm{(\nabla^2 \phi, \nabla^3 \zeta)}_{H^1}^2.
		\end{aligned}
	\end{equation}

\vspace{.2cm}

\subsection{Proof of the a priori estimates} 
Now we apply the steps shown in \cref{Sec-steps} to prove \cref{Prop-apriori}.
\begin{proof}
1) We first show \cref{est-apriori}$_1$. 
It follows from \cref{ic-Ec} that
\begin{equation}\label{ic-bdd}
	\begin{aligned}
		\Ec_\od^*(0) + \Ec_\md(0) & \lesssim \Ec^*(0) \lesssim \chi, \\
		\Ec_\od(0) + \Ec_{ho}(0) & \lesssim \Ec(0) + \e M^2 \lesssim 1 + \e M^2.
	\end{aligned}
\end{equation}	
Collecting \cref{ineq-od,ineq-md*,ineq-ho}, one has that
\begin{equation}\label{ineq-Et}
	\begin{aligned}
		\frac{d}{dt} \widetilde{\Ec}(t) + \Nc(t) & \lesssim (t+1)^{-\frac{5}{4}} \Ec_\md(t) + (\chi + \nu^{\frac{1}{2}} + \e) M^2 (t+1)^{-\frac{1}{2}},
	\end{aligned}
\end{equation}
where $ \widetilde{\Ec} $ and $ \Nc $ are two energy functionals, satisfying that 
\begin{equation}\label{Ect&Nc}
	\widetilde{\Ec} \sim \Ec_\od + \Ec_\md + \Ec_{ho}, \quad \Nc \sim \Nc_\od + \Nc_\md + \Nc_{ho}.
\end{equation}
Note that 
$$ \Ec_\md(t) \lesssim \widetilde{\Ec}(t) + \e M^2 (t+1)^{\frac{1}{2}}. $$
Then \cref{ineq-Et} implies that
\begin{equation}\label{ineq-Et*}
	\begin{aligned}
		\frac{d}{dt} \widetilde{\Ec}(t) + \Nc(t) & \lesssim (t+1)^{-\frac{5}{4}} \widetilde{\Ec}(t) + (\chi + \nu^{\frac{1}{2}} + \e) M^2 (t+1)^{-\frac{1}{2}}.
	\end{aligned}
\end{equation}
Using the Gronwall inequality, one has that
\begin{equation}\label{ode-Ect}
	\widetilde{\Ec}(t) + \int_0^t \Nc(\tau) d\tau \lesssim \big[ \widetilde{\Ec}(0) + (\chi + \nu^{\frac{1}{2}} + \e) M^2 \big] (t+1)^{\frac{1}{2}}.
\end{equation}
Comparing \cref{Ec,Ect&Nc}, one has that
\begin{equation}\label{sim-Ec}
	\pm \Ec(t) \lesssim \pm \widetilde{\Ec}(t)+ \e M^2 (t+1)^{\frac{1}{2}}.
\end{equation} 
Thus, \cref{ode-Ect} gives that
$$
\Ec(t) + \int_0^t \Nc(\tau) d\tau \lesssim \big[ \Ec(0) + (\chi + \nu^{\frac{1}{2}} + \e) M^2 \big] (t+1)^{\frac{1}{2}},
$$
which implies that
$$
\sup_{ t\in(0,T)} (t+1)^{-\frac{1}{2}} \Ec(t) \lesssim \Ec(0) + (\nu^{\frac{1}{2}} + \e) M^2.
$$
Recall the definition of $M$ in \cref{M&nu}. It holds that
$$
M^2 \lesssim \Ec(0) + (\nu^{\frac{1}{2}} + \e) M^2,
$$
which yields that
\begin{equation}\label{bdd-M}
	\sup_{ t\in(0,T)} (t+1)^{-\frac{1}{2}} \Ec(t) \leq M^2 \lesssim \Ec(0).
\end{equation}
Then the proof of \cref{est-apriori}$_1$ is finished.
\vspace{.2cm}
~\\
2) Then we prove \cref{est-apriori}$_2$. Since we have shown that $M^2 \lesssim \Ec(0) \lesssim 1,$ then the combination of \cref{ineq-od*,ineq-md*} yields that
\begin{equation}\label{ineq-AB*}
	\begin{aligned}
		\frac{d}{dt} \big( \Ec_\od^* + \Ec_\md \big) + \Nc_\od^* + \Nc_\md & \lesssim (t+1)^{-\frac{5}{4}} \Ec_\md(t) + \big(\chi \nu + \nu^3 + \e \Ec(0) \big) (t+1)^{-\frac{1}{2}}  \\
		& \quad + \big( \nu^{\frac{1}{2}} + \e^{\frac{1}{2}} \big) (t+1) \norm{\nabla^2 \phi}_{L^2}^2.
	\end{aligned}
\end{equation}
If we use the high-order estimate \cref{ineq-ho} to cancel the large energy $ \nu^{\frac{1}{2}} (t+1) \norm{\nabla^2 \phi}_{L^2}^2, $ then we will finally get the trivial result that
$$ \nu^2 \lesssim \nu^{\frac{1}{2}} + \cdots, $$
which cannot result in the smallness of $\nu$ shown in \cref{est-apriori}$_2.$
Nevertheless, the effective way to cancel this large energy is the use of the density dissipation \cref{est8*}. 
In fact, combining \cref{ineq-AB*,est8*}, one can get that
\begin{align*}
	& \frac{d}{dt} \Big( 2 C_0 \big( \Ec_\od^* + \Ec_\md\big) + (t+1) \Ac_{2, \e\phi} \Big) + C_0 \big( \Nc_\od^* + \Nc_\md \big) + (t+1) \norm{\nabla^2 \phi}_{L^2}^2 \\
	& \qquad \lesssim (t+1)^{-\frac{5}{4}}   \Ec_\md  + \big(\chi \nu + \nu^3 + \e \Ec(0) \big) (t+1)^{-\frac{1}{2}},
\end{align*}
where $C_0 \sim 1$ is a suitably large constant and we have used the fact that
$$
\e (t+1)^{-1} \Nc_\od \lesssim \e M^2 (t+1)^{-\frac{1}{2}} \lesssim \e \Ec(0) (t+1)^{-\frac{1}{2}}.
$$
By denoting  
\begin{equation}\label{Ect*}
	\widetilde{\Ec}^{*} := 2C_0 \big( \Ec_\od^* + \Ec_\md\big) + (t+1) \Ac_{2, \e\phi},
\end{equation}
then one can use \cref{bdd-A2} to get that
$$
\Ec_\md \lesssim \widetilde{\Ec}^{*} + \e \Ec(0) (t+1)^{-\frac{1}{2}}.
$$
Thus, it holds that
$$
\frac{d}{dt} \widetilde{\Ec}^{*}(t) \lesssim (t+1)^{-\frac{5}{4}} \widetilde{\Ec}^{*}(t) + \big(\chi \nu + \nu^3 + \e \Ec(0) \big) (t+1)^{-\frac{1}{2}},
$$
which yields that
\begin{align}
	(t+1)^{-\frac{1}{2}} \widetilde{\Ec}^{*}(t) 
	& \lesssim \widetilde{\Ec}^*(0) + \chi \nu + \nu^3 + \e \Ec(0). \label{est-Ec**}
\end{align}
Comparing \cref{Ec*,Ect*}, it holds that
$$
\pm \Ec^*(t) \lesssim \pm \widetilde{\Ec}^*(t) + \e \Ec(0) (t+1)^{-\frac{1}{2}}, \qquad t\in [0,T].
$$
Then it holds that
$$ \nu^2 = \sup_{ t\in(0,T)} (t+1)^{-\frac{1}{2}} \Ec^*(t) \lesssim \chi \nu + \nu^3 + \e \Ec(0), $$
which yields that
\begin{equation}
	\nu^2 = \sup_{ t\in(0,T)} (t+1)^{-\frac{1}{2}} \Ec^*(t) \lesssim \chi + \e \Ec(0).
\end{equation}
The proof of \cref{Prop-apriori} is completed.
	
\end{proof}

\vspace{.2cm}

In the following sections, we show the proof of each step shown in \cref{Sec-steps}.
To simplify the statement, we always use $ \e_0, \chi_0 $ and $ \nu_0 $ to denote some small positive constants, all of which are independent of $t$ and $M$. We also omit the statement of the assumptions that $\e\leq \e_0, \chi \leq \chi_0$ and $\nu\leq \nu_0.$

\vspace{.1cm}

\section{Estimates for anti-derivative variables}\label{Sec-od}

In this section, we first establish the $L^2$-estimates for the anti-derivatives, $(\Phi, \Zv),$ and finish the proof of the Steps 1 and 2 in \cref{Sec-steps}. 
Furthermore, we shall use the Green function theory to achieve a more accurate estimate for the tangential anti-derivatives, $\Zv_\perp = (Z_1,Z_2).$ The refined estimate plays an essential role in the following non-zero mode estimates, which overcomes the difficulty arising from the nonlinear convection without smallness. 

\subsection{Anti-derivative system and the nonlinear estimates}
Integrating \cref{equ-phipsi} with respect to $ \xp \in \Torus^2 $ yields the zero-mode system of $(\phi^\od,\,\psi^\od)(x_3,t)$:
\begin{equation}\label{equ-pert-od}
	\begin{cases}
		\p_t \phi^\od + \e^{-1} \p_3 \psi_3^\od = -\p_3 F_0, \\
		\p_t \psi^\od + \p_3 \big[ \big(\frac{m_3 \mv}{\rho} - \frac{\mt_3 \mvt}{\rhot} \big)^\od \big] + \e^{-1} \p_3 \big( p'(\rhot) \phi^\od \big)\E_3 + \e^{-2} \p_3 \big(\varpi(\rho,\rhot)\big)^\od \E_3 \\
		\qquad - \mu \p_3^2 \big[ \big(\frac{\mv}{\rho} - \frac{\mvt}{\rhot} \big)^\od \big] - (\mu+\lambda) \p_3^2 \big[ \big(\frac{m_3}{\rho} - \frac{\mt_3}{\rhot} \big)^\od \big] \E_3 = - \p_3 \Fv.
	\end{cases}
\end{equation}
Here recall the notations, $\E_3 = (0,0,1)^t$ and
$$
\begin{aligned}
	\varpi(\rho, \rhot) & = p(\rhot + \e \phi) -p(\rhot) - \e p'(\rhot) \phi \\
	& = \e^2 \int_0^1\int_0^1 p''(\rhot + \e r_1 r_2 \phi) dr_2 r_1 dr_1 \, \phi^2.
\end{aligned}
$$
By integrating \cref{equ-pert-od} with respect to $x_3$ from $-\infty$ to $x_3$, one can arrive at the system of the anti-derivative variables, $(\Phi, \Zv), $ defined by \cref{def: anti-derivative Phi Psi,Zv},
\begin{equation}\label{equ-PhiZ}
	\begin{cases}
		\p_t \Phi + \e^{-1} \p_3 Z_3 + \p_3 (\ut_3 \Phi) = -F_0, \\
		\p_t \Zv  +\e^{-1} p'(\rhot) \p_3 \Phi \E_3
		+\ut_3 \p_3 \Zv
		- \mu \p_3 \big(\frac{1}{\rhot} \p_3 \Zv \big)
		-(\mu+\lambda) \p_3 \big(\frac{1}{\rhot} \p_3 Z_3\big)\E_3 \\
		\qquad 
		= \e \uvt F_0  - \Fv - \e \mathbf{d}_{0} \Phi - \e \db_{1} \p_3 \Phi - \Qv_{1}^\od + \e \p_3 \Qv_{2}^\od,
	\end{cases}
\end{equation}
where the variable coefficients in the linear terms are given by
\begin{equation}\label{d}
	\begin{aligned}
		\db_0 & := \p_t\uvt + \ut_3 \p_3 \uvt - \mu \p_3 \Big(\frac{\p_3\uvt}{\rhot}\Big) - (\mu+\lambda)\p_3 \Big(\frac{\p_3 \ut_3}{\rhot}\Big) \E_3, \\
		\db_1 & := - \frac{\mu \p_3\uvt}{\rhot} - \frac{(\mu+\lambda)\p_3 \ut_3}{\rhot} \E_3, 
	\end{aligned}
\end{equation}
and the nonlinear terms are 
\begin{equation}\label{Q}
	\begin{aligned}
		\Qv_1
		& := \frac{m_3 \mv}{\rho} - \frac{\mt_3 \mvt}{\rhot} - \frac{\mvt}{\rhot} \psi_3 - \frac{\mt_3}{\rhot} \psi + \frac{\e \mt_3\mvt}{\rhot^2} \phi  + \e^{-2} \varpi(\rho,\rhot) \E_3, \\
		\Qv_2 
		& := \frac{\mu}{\e} \Big(\frac{\mv}{\rho} - \frac{\mvt}{\rhot}-\frac{1}{\rhot} \psi + \frac{\e \mvt}{\rhot^2} \phi \Big) + \frac{\mu+\lambda}{\e} \Big(\frac{m_3}{\rho} + \frac{\mt_3}{\rhot} - \frac{1}{\rhot} \psi_3 + \frac{\e\mt_3}{\rhot^2} \phi \Big) \E_3.
	\end{aligned}
\end{equation}

\vspace{.1cm}

Using \cref{bdd-diff,bdd-F}, one can obtain the estimates of the linear coefficients in \cref{d} directly.

\begin{Lem}\label{Lem-d}
	For $j=0,1,2,\cdots, $ then the linear coefficients \cref{d} satisfy that  
	\begin{equation}\label{est-d}
		\begin{aligned}
			& \norm{\p_3^j d_{0,3}}_{L^\infty} \lesssim \chi (t+1)^{-\frac{3+j}{2}}, \qquad \norm{\p_3^j \db_{0,\perp}}_{L^\infty} \lesssim (t+\Lambda)^{-\frac{3+j}{2}}, \\
			& \norm{\p_3^j d_{1,3}}_{L^\infty} \lesssim \chi (t+1)^{-\frac{2+j}{2}}, \qquad \norm{\p_3^j \db_{1,\perp}}_{L^\infty} \lesssim (t+\Lambda)^{-\frac{1+j}{2}}.
		\end{aligned}
	\end{equation}
\end{Lem}

\vspace{.1cm}

To achieve the uniform smallness of the partial energy functional \cref{Ec*}, we need to estimate the nonlinear term, $\Qv_1,$ carefully. 
We first decompose the nonlinear terms into different orders of $\e$,
\begin{equation}\label{Qc-Qh}
	- \Qv_{1}^\od + \e \p_3 \Qv_{2}^\od = -\widecheck{\Qv}_1^\od + \e \big(\widehat{\Qv}_1^\od + \p_3 \Qv_{2}^\od\big),
\end{equation}
where 
\begin{equation}\label{Q-12}
	\begin{aligned}
	\widecheck{\Qv}_1 & = \frac{\psi_3 \psi}{\rhot} + \sigma(\rho,\rhot) \phi^2 \E_3 \quad \text{with } \  \sigma(\rho,\rhot) = \int_0^1 \int_0^1 p''(\rhot+\e r_1 r_2 \phi) dr_2 r_1 dr_1, \\
	\widehat{\Qv}_1 & = \frac{\psi_3 \psi}{\rhot\rho} + \frac{1}{\rho} \big( \uvt\phi \psi_3 + \ut_3 \phi\psi - \e \ut_3 \uvt \phi^2 \big), \\
	\Qv_{2} & = \frac{1}{\rhot\rho} \big[ \mu \big( \e \uvt \phi^2 - \phi \psi \big) + (\mu+\lambda) \big(\e \uvt_3 \phi^2 - \phi \psi_3 \big) \E_3\big].
	\end{aligned}
\end{equation}
Since $\widecheck{\Qv}_1^\od$ has no small coefficient $\e,$ the nonlinear interactions of non-zero modes, for example $\frac{1}{\rhot} (\psi_3^\md \psi^\md)^\od, $ result in difficulties

\begin{Lem}\label{Lem-Q-nor}
	If $\e M\leq 1, $ then the zero-order nonlinearity, $ \widecheck{\Qv}_1^\od,$  in \cref{Qc-Qh} satisfies that
	\begin{equation}\label{L2-Q3}
		\begin{aligned}
			\norm{\widecheck{Q}_{1,3}^\od}_{L^2}^2 & \lesssim \nu^2  (t+1)^{-\frac{1}{2}} \norm{\nabla (\phi^\md, \psi^\md_3)}_{L^2}^2 + \nu^4 (t+1)^{-\frac{3}{2}}, \\
			\norm{\p_3 \widecheck{Q}_{1,3}^\od}_{L^2}^2 & \lesssim \nu M (t+1)^{-1} \norm{\nabla^2 \big(\phi^\md, \psi_3^\md\big)}_{L^2}^2 + \nu^4 (t+1)^{-\frac{5}{2}},
		\end{aligned}
	\end{equation}
	and 
	\begin{equation}\label{L2-Q1p}
		\norm{\widecheck{\Qv}_{1,\perp}^\od}_{L^2}^2 \lesssim  \nu^2 M^2 (t+1)^{-\frac{3}{2}}, \quad \norm{\p_3 \widecheck{\Qv}_{1,\perp}^\od}_{L^2}^2 \lesssim \nu M^3 (t+1)^{-\frac{5}{2}}.
	\end{equation}
	
	\end{Lem}
	
\begin{proof}
	
	It follows from \cref{Q-12}$_1$ that 
	\begin{equation}\label{Q-c}
		\widecheck{Q}_{1,3} = \frac{1}{\rhot} \psi_3^2 + \sigma(\rho, \rhot) \phi^2, \quad \widecheck{\Qv}_{1,\perp} = \frac{1}{\rhot} \psi_3 \psi_\perp.
	\end{equation}
	1) To prove \cref{L2-Q3}$_1$, the combination of Lemmas \ref{Lem-GN} and \ref{Lem-poin} implies that 
		\begin{align*}
			\norm{\widecheck{Q}_{1,3}}_{L^2}^2 \lesssim \norm{(\phi, \psi_3)}_{L^4}^4  & \lesssim \norm{(\phi^\od, \psi_3^\od)}_{L^4}^4 + \norm{(\phi^\md, \psi_3^\md)}_{L^4}^4 \\
			& \lesssim \norm{\p_3 (\phi^\od, \psi_3^\od)}_{L^2} \norm{(\phi^\od, \psi_3^\od)}_{L^2}^3 + \norm{\nabla (\phi^\md, \psi^\md_3)}_{L^2}^{4}.
		\end{align*}
		Using \cref{aL2-od,aL2-md}, one has that 
		\begin{equation}\label{L4-0}
			\norm{(\phi, \psi_3)}_{L^4}^4 \lesssim \nu^4 (t+1)^{-\frac{3}{2}} + \nu^2 (t+1)^{-\frac{1}{2}} \norm{\nabla (\phi^\md, \psi^\md_3)}_{L^2}^2,
		\end{equation}
		which implies \cref{L2-Q3}$_1$. 

		2) Note that $\abs{\p_3 (\sigma(\rho, \rhot))} \lesssim \abs{\p_3 \rhot} + \e \abs{\p_3 \phi}.$
		Then using \cref{bdd-diff,aLinf-pert}, if $\e M \leq 1,$ one can get that
		\begin{equation}\label{ineq-Q3-1}
			\norm{\p_3 \widecheck{Q}_{1,3}}_{L^2}^2 \lesssim \norm{ (\phi, \psi_3)}_{L^\infty}^2  \norm{\p_3 (\phi, \psi_3)}_{L^2}^2 + \chi^2 (t+1)^{-2} \norm{(\phi, \psi_3)}_{L^4}^4.
		\end{equation}
		For the first term on the right-hand side of \cref{ineq-Q3-1}, it holds that
		\begin{align}
			\norm{(\phi, \psi_3)}_{L^\infty}^2 & \lesssim \norm{(\phi^\od, \psi^\od_3)}_{L^\infty}^2 + \norm{\big(\phi^\md, \psi_3^\md\big)}_{L^\infty}^2 \notag \\
			& \lesssim \nu^2 (t+1)^{-1} + \norm{\nabla^2 \big(\phi^\md,  \psi_3^\md\big)}_{L^2}^2, \label{L-inf}
		\end{align}
		where we have used \cref{aLinf-od} and the fact (by Lemmas \ref{Lem-GN} and \ref{Lem-poin}) that $$\norm{\big(\phi^\md, \psi_3^\md\big)}_{L^\infty} \lesssim \norm{\big(\phi^\md, \psi_3^\md\big)}_{H^2} \lesssim \norm{\nabla^2 \big(\phi^\md,  \psi_3^\md\big)}_{L^2}.$$
		Besides, it follows from \cref{aL2-od} that
		\begin{equation}\label{L-2}
			\begin{aligned}
				\norm{\p_3 (\phi, \psi_3)}_{L^2}^2 & \lesssim \norm{\p_3 (\phi^\od, \psi_3^\od)}_{L^2}^2 + \norm{\nabla(\phi^\md, \psi_3^\md)}_{L^2}^2 \\
				& \lesssim \nu^2 (t+1)^{-\frac{3}{2}} + \norm{\nabla (\phi^\md, \psi_3^\md)}_{L^2}^2.
			\end{aligned}
		\end{equation}
		Then plugging \cref{L4-0,L-2,L-inf} into \cref{ineq-Q3-1} yields that
		\begin{equation}\label{ineq-Q3-2}
			\begin{aligned}
			\norm{\p_3 \widecheck{Q}_{1,3}}_{L^2}^2 & \lesssim \nu^4 (t+1)^{-\frac{5}{2}} + \nu^2 (t+1)^{-1} \norm{\nabla^2 \big(\phi^\md, \psi_3^\md\big)}_{L^2}^2 \\
			& \quad + \norm{\nabla \big(\phi^\md,  \psi_3^\md\big)}_{L^2}^2 \norm{\nabla^2 \big(\phi^\md,  \psi_3^\md\big)}_{L^2}^2.
		\end{aligned}
		\end{equation}
		Using \cref{aL2-md,poincare}, one has that
		\begin{equation*}
			\begin{aligned}
				\norm{\nabla (\phi^\md, \psi_3^\md)}_{L^2}^2 & \lesssim \nu (t+1)^{-\frac{1}{4}} \norm{\nabla^2 (\phi^\md, \psi_3^\md)}_{L^2}, \\
				\norm{\nabla^2 (\phi^\md, \psi_3^\md)}_{L^2}^2 & \lesssim M (t+1)^{-\frac{3}{4}} \norm{\nabla^2 (\phi^\md, \psi_3^\md)}_{L^2},
			\end{aligned}
		\end{equation*}
		which yields that
		 \begin{equation*}
			\norm{\nabla (\phi^\md, \psi_3^\md)}_{L^2}^2 \norm{\nabla^2 (\phi^\md, \psi_3^\md)}_{L^2}^2 \lesssim \nu M (t+1)^{-1} \norm{\nabla^2 (\phi^\md, \psi_3^\md)}_{L^2}^2.
		\end{equation*}
		This, together with \cref{ineq-Q3-2}, yields \cref{L2-Q3}$_2$.
		
		3) To show \cref{L2-Q1p}, it follows from \cref{Q-c} that
		$$
		\norm{\widecheck{\Qv}_{1,\perp}^\od}_{L^2}^2  \lesssim \norm{\psi_3^\od}_{L^2}^2 \norm{\psi_\perp^\od}_{L^\infty}^2 + \norm{ \psi^\md}_{L^4}^4, 
		$$
		and
		\begin{equation}\label{p1-Qc}
			\begin{aligned}
			\norm{\p_3 \widecheck{\Qv}_{1,\perp}^\od}_{L^2}^2 & \lesssim \norm{\psi_3^\od}_{L^\infty}^2 \norm{\p_3 \psi_\perp^\od}_{L^2}^2 + \norm{\p_3 \psi_3^\od}_{L^2}^2 \norm{ \psi_\perp^\od}_{L^\infty}^2  \\
			& \quad + \norm{\psi^\md}_{L^\infty}^2 \norm{\p_3  \psi^\md}_{L^2}^2 + \e^2 \chi^2 (t+1)^{-2} \norm{\widecheck{\Qv}_{1,\perp}^\od}_{L^2}^2.
		\end{aligned}
		\end{equation}
		Using \cref{aL2-od,aL2-md,aLinf-od}, it holds that
		\begin{equation*}
			\begin{aligned}
				& \norm{\psi_3^\od}_{L^2}^2 \norm{\psi_\perp^\od}_{L^\infty}^2 \lesssim \nu^2 M^2 (t+1)^{-\frac{3}{2}}, \\
				& \norm{\psi^\md}_{L^4}^4 \lesssim \norm{\psi^\md}_{H^1}^4 \lesssim \norm{\nabla \psi^\md}_{L^2}^4 \\
				& \qquad \quad \ \lesssim \norm{\nabla  \psi^\md}_{L^2}^2 \norm{\nabla^2\psi^\md}_{L^2}^2 \lesssim \nu^2 M^2 (t+1)^{-2}, 
			\end{aligned}
		\end{equation*}
		which yields that $\norm{\widecheck{\Qv}_{1,\perp}^\od}_{L^2}^2 \lesssim  \nu^2 M^2 (t+1)^{-\frac{3}{2}}.$ 

		To estimate $\p_3 \widecheck{\Qv}_{1,\perp}^\od$, it follows from \cref{aLinf-md,aL2-md} that
		$$
		\norm{\psi^\md}_{L^\infty}^2 \norm{\p_3  \psi^\md}_{L^2}^2 \lesssim \nu M^3 (t+1)^{-\frac{5}{2}}.
		$$
		Besides, one can use \cref{aLinf-od,aL2-od}, together with the $L^2$-bound of $\widecheck{\Qv}_{1,\perp}^\od,$ to estimate the remaining terms in \cref{p1-Qc}, which can complete the proof of \cref{L2-Q1p}.
\end{proof}

\begin{Lem}\label{Lem-Qh2}
	If $\e M\leq 1, $ then the nonlinear terms, $\widehat{\Qv}_1^\od$ and $\Qv_{2}^\od,$ in \cref{Qc-Qh} satisfy that 
	\begin{equation}\label{L2-Qh2}
		\begin{aligned}
			\norm{\widehat{\Qv}_1}_{L^2}^2 + \norm{\Qv_{2}}_{L^2}^2 & \lesssim M^4 (t+1)^{-\frac{3}{2}},  \\
			\norm{\p_3 \widehat{\Qv}_1}_{L^2}^2 + \norm{\p_3 \Qv_{2}}_{H^1}^2 & \lesssim M^4 (t+1)^{-\frac{5}{2}}. 
		\end{aligned}
	\end{equation}
	
\end{Lem}

\begin{proof}

Using \cref{aL2-pert,aLinf-pert}, one can get that
\begin{equation}\label{ineq-1-per}
	\begin{aligned}
		\norm{\widehat{\Qv}_1}_{L^2} + \norm{\Qv_2}_{L^2} \lesssim \norm{\big(\phi, \psi\big)}_{L^\infty} \norm{\big(\phi, \psi\big)}_{L^2} \lesssim M^2 (t+1)^{-\frac{3}{4}},
	\end{aligned}
\end{equation}
which gives \cref{L2-Qh2}$_1$.
Similarly, one can get that if $\e M \leq 1, $
\begin{align*}
	& \norm{\p_3 \widehat{\Qv}_1}_{L^2} + \norm{\p_3 \Qv_2}_{L^2} \\
	& \quad \lesssim \norm{(\phi, \psi)}_{L^\infty} \norm{\p_3 (\phi, \psi)}_{L^2} + \big[(t+1)^{-\frac{1}{2}} + \e \norm{\p_3 \phi}_{L^\infty} \big] \norm{(\phi, \psi)}_{L^\infty} \norm{(\phi, \psi)}_{L^2}  \\
	& \quad \lesssim M^2 (t+1)^{-\frac{5}{4}}.
\end{align*}
Furthermore, note that $ \p_3^2 \Qv_2 = \frac{1}{\rho} \p_3^2 (\rho \Qv_2) -\frac{2}{\rho} \p_3 \rho \p_3 \Qv_2 - \frac{1}{\rho} \p_3^2 \rho \Qv_2 . $
It follows from \cref{Q-12}$_3$ that $\norm{\Qv_2}_{L^\infty} \lesssim \norm{(\phi,\psi)}_{L^\infty}^2 \lesssim M^2 (t+1)^{-1},$ and
\begin{align*}
	\norm{\p_3^2 (\rho \Qv_2)}_{L^2} & \lesssim \norm{(\phi, \psi)}_{L^\infty} \norm{\p_3^2 (\phi, \psi)}_{L^2} + \norm{\p_3(\phi,\psi)}_{L^\infty} \norm{\p_3(\phi,\psi)}_{L^2} \\
	& \quad +  \sum_{i=0}^{1} \e  (t+1)^{-1+\frac{i}{2}} \norm{(\phi, \psi)}_{L^\infty} \norm{ \p_3^i\big(\phi, \psi\big)}_{L^2} \\
	& \lesssim M^2 (t+1)^{-\frac{5}{4}}.
\end{align*}
Then one has that
\begin{align*}
	\norm{\p_3^2 \Qv_2}_{L^2} & \lesssim \norm{\p_3^2 (\rho \Qv_2)}_{L^2} + \big( \norm{\p_3 \rhot}_{L^\infty} + \e \norm{\p_3 \phi}_{L^\infty} \big) \norm{\p_3 \Qv_2}_{L^2} \\
	& \quad + \norm{\p_3^2 \rhot}_{L^\infty} \norm{\Qv_2}_{L^2} + \e \norm{\p_3^2 \phi}_{L^2}  \norm{\Qv_2}_{L^\infty} \\
	& \lesssim M^2 (t+1)^{-\frac{5}{4}} + \e M^3 (t+1)^{-\frac{7}{4}} \\
	& \lesssim M^2 (t+1)^{-\frac{5}{4}}, \qquad \text{if } \ \e M\leq 1.
\end{align*}

\end{proof}


\subsection{Proof of Step 1}

With the linear and nonlinear estimates obtained in Lemmas \ref{Lem-d}--\ref{Lem-Qh2}, we are able to prove the Step 1 in \cref{Sec-steps}. 
Extracting the equations of $(\Phi,\,Z_3)$ from \cref{equ-PhiZ} gives that
\begin{equation}\label{equ-PhiZ3}
	\begin{cases}
		\p_t \Phi + \e^{-1} \p_3 Z_3 + \p_3 (\ut_3 \Phi) = -F_0, \\
		\p_t Z_3 + \e^{-1} p'(\rhot) \p_3 \Phi + \ut_3 \p_3 Z_3  - \mut \p_3 \big(\frac{1}{\rhot} \p_3 Z_3\big) + \e d_{0,3} \Phi + \e d_{1,3} \p_3 \Phi \\
		\qquad = \e \ut_3 F_0 - F_3 - Q_{1,3}^\od  + \e \p_3 Q_{2,3}^\od,
	\end{cases}
\end{equation}
where $\mut= 2\mu+\lambda.$
The weighted $H^2$-estimates of the small anti-derivatives in \cref{Ec*} will be derived in the following three lemmas. 
\begin{Lem}\label{Lem-est1}
	If $\e^{\frac{1}{2}} M \leq 1$, then it holds that
	\begin{equation}\label{est-1}
		\begin{aligned}
			\frac{d}{dt} \Ac_{-1,\od} + \norm{\p_3 (\Phi, Z_3)}_{L^2}^2 &  \lesssim (\chi\nu + \nu^3 + \e M^2) (t+1)^{-\frac{1}{2}},
		\end{aligned}
	\end{equation}
	where $ \Ac_{-1,\od} $ is a functional satisfying that 
	\begin{equation}\label{AB-od-1}
		\Ac_{-1,\od} \sim \norm{(\Phi,Z_3)}_{L^2}^2 + \e^2 \norm{\p_3 \Phi}_{L^2}^2.
	\end{equation}
\end{Lem}

\begin{proof}
Multiplying $ p'(\rhot) \Phi $ and $ Z_3 $ on \cref{equ-PhiZ3}$ _1 $ and \cref{equ-PhiZ3}$ _2 $, respectively, and adding the resulting two equations together, one has that
\begin{equation}\label{eq1}
	\begin{aligned}
		& \p_t \Big[ \frac{1}{2} \big( p'(\rhot) \Phi^2 + \abs{Z_3}^2\big) \Big] + \frac{\mut}{\rhot} \abs{\p_3 Z_3}^2 = \p_3 (\cdots) + I_1 + I_2,
	\end{aligned}
\end{equation}
where  
\begin{align*}
	I_1 &= \frac{1}{2} p''(\rhot) \big( \p_t\rhot + \ut_3 \p_3 \rhot \big) \Phi^2 - p'(\rhot) F_0 \Phi + \frac{1}{2} \p_3 \ut_3 Z_3^2 \\
	&\quad + \big[ \e^{-1} p''(\rhot) \p_3 \rhot - \e d_{0,3} + \e \p_3 d_{1,3} \big] \Phi Z_3 + \e d_{1,3}  \Phi \p_3 Z_3 + ( \e \ut_3 F_0 - F_3) Z_3, \\
	I_2 &= - Q_{1,3}^\od Z_3 - \e Q_{2,3}^\od \p_3 Z_3.
\end{align*}
It follows from \cref{Lem-F}, \cref{aL2-anti,est-d} that
\begin{align}
	\int_\R I_1 dx_3 & \lesssim \chi (t+1)^{-1} \norm{(\Phi, Z_3)}_{L^2}^2 + \chi (t+1)^{-\frac{3}{4}} \norm{(\Phi, Z_3)}_{L^2} + \chi \norm{\p_3 Z_3}_{L^2}^2 \notag \\
	& \lesssim \chi \nu (t+1)^{-\frac{1}{2}}. \label{est-I1}
\end{align}
Besides, it follows from \cref{aL2-pert,aLinf-anti,Q-12} that
\begin{equation}\label{est-I2}
	\begin{aligned}
		\int_\R I_2 dx_3 \lesssim \norm{Z_3}_{W^{1,\infty}} \norm{\big(\phi, \psi_3\big)}_{L^2}^2 \lesssim \nu^3 (t+1)^{-\frac{1}{2}}.
	\end{aligned}
\end{equation}
Combining \cref{est-I1,est-I2} and using the fact that $ \rhot \leq 2 \rhob, $ one can obtain that
\begin{equation}\label{est1}
	\begin{aligned}
		& \frac{d}{dt} \Big[\int_\R \big( p'(\rhot)  \abs{\Phi}^2 + \abs{Z_3} ^2 \big) dx_3\Big] + \frac{\mut}{2\rhob}  \norm{\p_3 Z_3}^2_{L^2} \lesssim (\chi \nu + \nu^3) (t+1)^{-\frac{1}{2}}.
	\end{aligned}
\end{equation}

To achieve the dissipation of $\p_3\Phi,$ one can multiply $ \frac{\e^2 \mut}{\rhot} \p_3 \Phi $ and $ \e \p_3 \Phi $ on $ \p_3 $\cref{equ-PhiZ3}$ _1 $ and \cref{equ-PhiZ3}$ _2 $, respectively, to get that
\begin{align}
	& \frac{d}{dt} \int_\R \Big(\frac{\e^2 \mut}{2\rhot} \abs{\p_3 \Phi}^2 + \e Z_3 \p_3 \Phi \Big) dx_3 + \int_\R \frac{1}{2} p'(\rhot) \abs{\p_3 \Phi}^2 dx_3 \notag \\
	& \qquad \lesssim \norm{\p_3 Z_3}_{L^2}^2 + \e \chi \nu (t+1)^{-1} + \e^2 \norm{Q_{1,3}^\od}_{L^2}^2 + \e^2 \norm{\p_3 Q_{2,3}^\od}_{L^2}^2 \notag \\
	& \qquad \lesssim \norm{\p_3 Z_3}_{L^2}^2 + \e \chi \nu (t+1)^{-1} + \e^2 M^4 (t+1)^{-\frac{3}{2}}. \label{est1.5}
\end{align}
Here we have used the a priori bounds in \cref{Sec-bound}, together with the fact from \cref{L2-Qh2,L2-Q3} that
$$
\norm{\Qv_1}_{L^2}^2 \lesssim M^4 (t+1)^{-3/2}, \quad \norm{\p_2 \Qv_2}_{L^2}^2 \lesssim M^4 (t+1)^{-5/2}.
$$

Combining \cref{est1,est1.5}, one can complete the proof of \cref{est-1}.

\end{proof}


Similarly, we can derive the estimates of $\norm{\p_3 (\Phi,Z_3)}_{L^2}$.
\begin{Lem}\label{Lem-est2}
	If $ \e^{\frac{1}{2}} M \leq 1,$ then it holds that
	\begin{equation}\label{est-2}
		\begin{aligned}
			& \frac{d}{dt} \big[(t+1) \Ac_{0,\od}\big] + (t+1) \norm{\p_3^2 (\Phi, Z_3)}_{L^2}^2 \\
			& \qquad \lesssim \norm{\p_3 (\Phi, Z_3)}_{L^2}^2  + \nu^2 (t+1)^{\frac{1}{2}} \norm{\nabla \big(\phi^\md, \psi_3^\md\big)}_{L^2}^2 \\
			& \qquad\quad + (\chi \nu + \nu^4 + \e M^2) (t+1)^{-\frac{1}{2}},
		\end{aligned}
	\end{equation}
	where $ \Ac_{0,\od} $ is a functional satisfying that 
	\begin{equation}\label{AB-od-2}
		\Ac_{0,\od} \sim \norm{\p_3 (\Phi,Z_3)}_{L^2}^2 + \e^2 \norm{\p_3^2 \Phi}_{L^2}^2.
	\end{equation}
\end{Lem}

\begin{proof}

Multiplying $ p'(\rhot) \p_3\Phi$ and  $-\p_3^2 Z_3 $ on $\p_3$\cref{equ-PhiZ3}\textsubscript{1} and \cref{equ-PhiZ3}\textsubscript{2}, respectively, one can get that
\begin{equation}\label{est4}
	\begin{aligned}
		& \frac{d}{dt} \int_\R \big(p'(\rhot) \abs{\p_3 \Phi}^2 + \abs{\p_3 Z_3}^2\big) dx_3 + \int_\R \frac{\mut}{\rhot} \abs{\p_3^2 Z_3} ^2 dx_3 \\
		& \qquad \lesssim \chi \nu (t+1)^{-\frac{3}{2}} + \norm{\widecheck{Q}_{1,3}^\od}_{L^2}^2 + \e^2 \big( \norm{\widehat{Q}_{2,3}^\od}_{L^2}^2 + \norm{\p_3 Q_{2,3}^\od}_{L^2}^2 \big) \\
		& \qquad \lesssim \nu^2 (t+1)^{-\frac{1}{2}} \norm{\nabla (\phi^\md, \psi^\md_3)}_{L^2}^2 + \big(\chi \nu + \nu^4 + \e^2 M^4\big) (t+1)^{-\frac{3}{2}},
		\end{aligned}
\end{equation}
where we have used \cref{L2-Q3}$_1$ and \cref{L2-Qh2}, and the omitted proof of the linear estimates are similar to \cref{est-I1}.

By multiplying $ \frac{\e^2 \mut}{\rhot} \p_3^2 \Phi $ and $\e\p_3^2\Phi $ on $ \p_3^2 $\cref{equ-PhiZ3}$_1$ and $\p_3$\cref{equ-PhiZ3}$_2$, respectively, and using \cref{L2-Q3}$_1$ and \cref{L2-Qh2}, one can get that
\begin{equation}\label{est5}
	\begin{aligned}
		& \frac{d}{dt} \Big[ \int_\R \Big(\frac{\e^2 \mut}{2\rhot}\abs{\p_3^2\Phi}^2 + \e \p_3 Z_3 \p_3^2 \Phi\Big) dx_3 \Big] + \int_\R \frac{1}{2} p'(\rhot) \abs{\p_3^2\Phi}^2 dx_3 \\
		& \qquad \lesssim \norm{\p_3^2 Z_3}_{L^2}^2 +  \e \chi \nu (t+1)^{-2} + \e^2 \norm{\p_3 Q_{1,3}^\od}_{L^2}^2 + \e^2 \norm{\p_3^2 Q_{2,3}^\od}_{L^2}^2 \\
		& \qquad \lesssim \norm{\p_3^2 Z_3}_{L^2}^2 +  \big(\e \chi \nu + \e^2 M^4 \big) (t+1)^{-2}.
	\end{aligned}
\end{equation}
Combining \cref{est4,est5} can imply \cref{est-2}.

\end{proof}

For the high order estimates $\norm{\p_3^2 \big(\Phi, Z_3\big)}_{L^2}$, we can only get the dissipation of $\p_3^2Z_3$.

\begin{Lem}\label{Lem-est3}
	If $ \big( \e^{\frac{1}{2}} + \nu^{\frac{1}{2}} \big) M \leq 1,$ then it holds that
	\begin{equation}\label{est-3}
		\begin{aligned}
			& \frac{d}{dt} \big[(t+1)^2 \Ac_{1,\od}\big] + (t+1)^2 \norm{\p_3^3 Z_3}_{L^2}^2 \\
			& \qquad \lesssim (t+1) \norm{\p_3^2 (\Phi, Z_3)}_{L^2}^2  + \nu^{\frac{1}{2}} (t+1) \norm{\nabla^2 (\phi^\md, \psi_3^\md)}_{L^2}^2 \\
			& \qquad\quad +  (\chi\nu + \nu^4+ \e M^2) (t+1)^{-\frac{1}{2}},
		\end{aligned}
	\end{equation}
	where $ \Ac_{1,\od} $ is a functional satisfy that
	\begin{equation}\label{AB-3}
		\Ac_{1,\od} \sim \norm{\p_3^2 \big(\Phi, Z_3\big)}_{L^2}^2.
	\end{equation}
\end{Lem}

\begin{proof}
	
Multiplying $ (p'(\rhot) \p_3^2 \Phi,\, \p_3^2 Z_3) $ on $ \p_3^2 $\cref{equ-PhiZ3}, one can get that
\begin{equation}\label{est6}
	\begin{aligned}
		& \frac{d}{dt}  \Big( \int_\R p'(\rhot) \abs{\p^2_3\Phi}^2 dx_3 + \norm{\p^2_3 Z_3}_{L^2}^2 \Big) + \frac{\mut}{2 \rhob} \norm{\p_3^3 Z_3}_{L^2}^2 \\
		& \qquad \lesssim \chi \nu (t+1)^{-\frac{5}{2}} + \norm{\p_3 \widecheck{Q}_{1,3}^\od}_{L^2}^2 + \e^2 \big(\norm{\p_3 \widehat{Q}_{1,3}^\od}_{L^2}^2 + \norm{\p_3^2 Q_{2,3}^\od}_{L^2}^2\big).
	\end{aligned}
\end{equation}
This, together \cref{L2-Q3}$_2$ and \cref{L2-Qh2}$_2$, yields \cref{est-3}.
\end{proof}


It follows from \cref{poincare,wv-md} that
\begin{equation}\label{ineq-md}
	\begin{aligned}
	\norm{\nabla (\phi^\md, \psi^\md)}_{L^2} & \lesssim \norm{\nabla^2 (\phi^\md, \psi^\md)}_{L^2} \lesssim \norm{\nabla^2 (\phi^\md, \wv^\md)}_{L^2}.
\end{aligned}
\end{equation}
Then collecting Lemmas \ref{Lem-est1}--\ref{Lem-est3}, one can get \cref{ineq-od*}. 

\vspace{.3cm}

\subsection{Proof of Step 2}\label{Sec-perp}
In this section, we show the estimates for the anti-derivative variable, $\Zv_\perp$, which is associated with the tangential velocity. These estimates can complete the proof of Step 2 in \cref{Sec-steps}. 
By \cref{equ-PhiZ,Qc-Qh}, the system of $\Zv_\perp=(Z_1,\,Z_2)$ reads
\begin{equation}\label{equ-Zp}
	\begin{aligned}
		& \p_t \Zv_\perp  + \ut_3 \p_3 \Zv_\perp - \mu \p_3 \Big(\frac{1}{\rhot} \p_3 \Zv_\perp \Big)
		+ \e \mathbf{d}_{0,\perp} \Phi  + \e \db_{1,\perp} \p_3 \Phi \\
		&\qquad = \e \uvt_\perp F_0 - \Fv_\perp - \widecheck{\Qv}_{1,\perp}^\od + \e \widehat{\Qv}_{1,\perp}^\od + \e \p_3 \Qv_{2,\perp}^\od.
	\end{aligned}
\end{equation}

\begin{Lem}\label{Lem-est4}
	If $\big( \e^{\frac{1}{2}} + \nu^{\frac{1}{2}} \big) M \leq 1, $ then for $j=0,1,2,$ it holds that
	\begin{equation}\label{est-4}
		\begin{aligned}
			& \frac{d}{dt} \big( \norm{\p_3^j \Zv_\perp}_{L^2}^2 \big) + \frac{\mu}{8\rhob} \norm{\p_3^{j+1} \Zv_\perp}_{L^2}^2  \lesssim (\chi + \nu^{\frac{1}{2}} + \e) M^2 (t+1)^{-\frac{1+2j}{2}}.
		\end{aligned}
	\end{equation}
\end{Lem}

\begin{proof}
	
	Multiplying $\Zv_\perp$ on \cref{equ-Zp}, and using the bounds in \cref{Sec-bound}, \cref{L2-Q1p,L2-Qh2}, one can get that
	\begin{align*}
		& \frac{d }{dt } \big(\norm{\Zv_\perp}_{L^2}^2 \big) + \int_\R \frac{\mu}{2\rhot} \abs{\p_3 \Zv_\perp}^2 dx_3 \\
		& \qquad \lesssim \chi (t+1)^{-1} \norm{\Zv_\perp}_{L^2}^2 + \Big[\e \sum_{j=0}^{1} (t+1)^{-\frac{2-j}{2}} \norm{\p_3^j \Phi}_{L^2} + \chi (t+1)^{-\frac{3}{4}}\Big] \norm{\Zv_\perp}_{L^2} \\
		& \qquad \quad + \norm{\Zv_\perp}_{L^2} \big(\norm{\widecheck{\Qv}_{1,\perp}^\od}_{L^2} + \e^2 \norm{\widehat{\Qv}_{1,\perp}^\od}_{L^2}\big) + \e^2 \norm{\Qv_{2,\perp}^\od}_{L^2}^2 \\
		& \qquad \lesssim (\chi M^2 + \e \nu M + \chi M + \nu M^2 + \e^2 M^4) (t+1)^{-\frac{1}{2}}.
	\end{align*} 
	Similarly, for $j=1,2,$ multiplying $\p_3^j \Zv_\perp$ on $ \p_3^j $\cref{equ-Zp} yields that
	\begin{align*}
		& \frac{d}{dt} \big(\norm{\p_3^j \Zv_\perp}_{L^2}^2 \big) + \int_\R \frac{\mu}{2\rhot} \abs{\p_3^{j+1} \Zv_\perp}^2 dx_3 \\
		& \qquad \lesssim (\chi + \e) M^2 (t+1)^{-\frac{1+2j}{2}} + \norm{\p_3^{j-1}\widecheck{\Qv}_{1,\perp}^\od}_{L^2}^2 \\
		& \quad \qquad + \e^2 \big(\norm{\p_3^{j-1}\widehat{\Qv}_{1,\perp}^\od}_{L^2}^2 + \norm{\p_3^j \Qv_{2,\perp}^\od}_{L^2}^2\big) \\
		& \qquad \lesssim (\chi + \nu M + \e) M^2 (t+1)^{-\frac{1+2j}{2}}.
	\end{align*}
	Recall that we have assumed that $M\geq 1$ without loss of generality (see \cref{M&nu}). Then the proof is finished.
\end{proof}

It follows from \cref{est-4} that there exist some constant $c_j$ for $j=0,1,2,$ depending only on $\mu$ and $\rhob$, such that
\begin{align}
	& \frac{d}{dt} \big( \sum_{j=0}^{2} c_j (t+1)^j \norm{\p_3^j \Zv_\perp}_{L^2}^2 \big) + \sum_{j=0}^{2} (t+1)^j \norm{\p_3^{j+1} \Zv_\perp}_{L^2}^2 \\
	& \qquad \lesssim (\chi + \nu^{\frac{1}{2}} + \e) M^2 (t+1)^{-\frac{1}{2}}.
\end{align}
This, together with \cref{ineq-od*} and the fact that $\norm{\nabla^2 (\phi^\md, \wv^\md)}_{L^2}^2 \lesssim M^2 (t+1)^{-\frac{3}{2}},$ can yield \cref{ineq-od}.

\vspace{.1cm}

\subsection{Parabolic type equations of tangential zero modes}

We will see in \cref{Sec-md} that in the non-zero mode estimates, the nonlinear convection gives rise to an interaction of the zero mode and non-zero mode, $\p_3 \wv^\od w_3^\md,$ which satisfies that
$$ \p_3 \wv^\od w_3^\md= \p_3^2 \Zv \, w_3^\md + O(1) \e \norm{\Phi}_{W^{1,\infty}} \abs{w_3^\md}. $$ 
Note that the tangential part of the anti-derivative variable, $\p_3^2 \Zv_\perp,$ has no smallness. Thus, even with the aid of Poincar\'{e} inequality, the nonlinear term, $\p_3^2 \Zv_\perp \, w_3^\md,$ cannot be controlled by the dissipation.
Furthermore, the decay rate $\norm{\p_3^2 \Zv_\perp}_{L^\infty} \lesssim M (t+1)^{-3/4}$, which is obtained through the $L^2$-method (see \cref{aLinf-anti}), is not time integrable and makes it difficult to establish the estimate of $\norm{\wv^\md}_{L^2}$, either. 

To overcome this difficulty, we observe that for $\Zv_\perp$, the main structure of the system \cref{equ-Zp} is of the  following parabolic type. Namely, we rewrite \cref{equ-Zp} into
\begin{equation}\label{equ-parabolic}
    \p_t \Zv_\perp = \frac{\mu}{\rhob} \p_3^2 \Zv_\perp + \mathcal{S}_\perp,
\end{equation}
where the remainder $\mathcal{S} = (\mathcal{S}_\perp, \mathcal{S}_3)(x_3,t) \in \R^3$ is given by
\begin{align*}
    \mathcal{S} & = - \frac{\mu(\rhot-\rhob)}{\rhob \rhot} \p_3^2 \Zv -\Big(\ut_3 + \frac{\mu \p_3 \rhot}{\rhot^2}\Big) \p_3 \Zv - \e \big(\db_{0} \Phi + \db_{1} \p_3 \Phi \big) \\
    & \quad + \e \uvt F_0 - \Fv - \widecheck{\Qv}_1^\od + \e \widehat{\Qv}_1^\od + \e \p_3 \Qv_2^\od.
\end{align*}
By using \cref{bdd-diff,L2-Q1p,L2-Qh2}, for $j=0,1,$ it holds that
\begin{equation*}
	\norm{\p_3^j \mathcal{S}_\perp}_{L^2} \lesssim (\chi M + \nu^{\frac{1}{2}} M^{\frac{3}{2}} + \e M^2) (t+1)^{-\frac{3+2j}{4}}.
\end{equation*}
Thanks to  the small parameters, $\chi, \e$ and $\nu,$ if $\big(\chi^{\frac{1}{2}} + \nu^{\frac{1}{6}} + \e^{\frac{1}{4}} \big) M \leq 1, $ the source term $\mathcal{S}$ in \cref{equ-parabolic} has smallness,
\begin{equation}\label{est-S}
	\norm{\p_3^j \mathcal{S}_\perp}_{L^2} \lesssim (\chi^\frac{1}{2} + \nu^\frac{1}{4} + \e^{\frac{1}{2}}) (t+1)^{-\frac{3+2j}{4}}, \qquad j=0,1.
\end{equation}
Then through the standard Green function theory, we are able to achieve a more refined estimate for $\Zv_\perp$ than \cref{aLinf-anti}. 
It turns out that the large part of $\Zv_\perp $ comes totally from its initial data, and it achieves a faster decay rate than the small one arising from the source $\mathcal{S}.$

\begin{Lem}\label{Lem-Zp-inf}
    If $ \big(\chi^{\frac{1}{2}} + \nu^{\frac{1}{6}} + \e^{\frac{1}{4}} \big) M \leq 1,$ then for $ j =1,2,$ it holds that
	\begin{equation}\label{Zp-inf}
		\begin{aligned}
			\norm{\p_3^j \Zv_\perp}_{L^\infty} & \lesssim (t+1)^{-\frac{j}{2}-\frac{1}{4}} + (\chi^\frac{1}{2} + \nu^\frac{1}{4} + \e^{\frac{1}{2}}) (t+1)^{-\frac{j}{2}},
		\end{aligned}  \quad 0\leq t\leq T.
	\end{equation}
\end{Lem}
\begin{proof}

	With the use of \cref{Lem-ic} and \cref{est-S}, it follows from \cref{equ-parabolic} and the classical parabolic theory  that
	\begin{itemize}
		\item if $0\leq t\leq 1,$ then
		\begin{equation*}
			\begin{aligned}
				\norm{\p_3^{j} \Zv_\perp}_{L^\infty} 
				& \lesssim \norm{\p_3^{j} \Zv_\perp(\cdot,0)}_{L^\infty} + \int_0^{1} \tau^{-\frac{2j-1}{4}} \norm{\p_3 \mathcal{S}_\perp(\cdot,t-\tau)}_{L^2} d\tau \\
				& \lesssim \Ec(0) + \chi^\frac{1}{2} + \nu^\frac{1}{4} + \e^{\frac{1}{2}} \\
				& \lesssim 1.
			\end{aligned}
		\end{equation*}
		
		\item if $t\geq 1, $ then 
		\begin{equation*}
			\begin{aligned}
				\norm{\p_3^{j} \Zv_\perp}_{L^\infty} 
				& \lesssim t^{-\frac{2j+1}{4}} \norm{\Zv_\perp(\cdot,0)}_{L^2} + \int_0^{\frac{t}{2}} \tau^{-\frac{2j-1}{4}} \norm{\p_3 \mathcal{S}_\perp(\cdot,t-\tau)}_{L^2} d\tau \\
				& \quad + \int_{\frac{t}{2}}^t \tau^{-\frac{2j+1}{4}} \norm{\mathcal{S}_\perp(\cdot,t-\tau)}_{L^2}  d\tau \\
				& \lesssim (t+1)^{-\frac{2j+1}{4}} + (\chi^\frac{1}{2} + \nu^\frac{1}{4} + \e^{\frac{1}{2}}) (t+1)^{-\frac{j}{2}},
			\end{aligned}
		\end{equation*}
		where we have used the fact that
		\begin{align*}
			& \int_0^{\frac{t}{2}} \tau^{-\frac{2j-1}{4}} (t-\tau+1)^{-\frac{5}{4}} d\tau + \int_{\frac{t}{2}}^t \tau^{-\frac{2j+1}{4}} (t-\tau+1)^{-\frac{3}{4}} d\tau \\
			& \qquad \lesssim (t+1)^{-\frac{5}{4}} \int_0^{\frac{t}{2}} \tau^{-\frac{2j-1}{4}} d\tau + t^{-\frac{2j+1}{4}} \int_{\frac{t}{2}}^t  (t-\tau+1)^{-\frac{3}{4}} d\tau \\
			& \qquad \lesssim (t+1)^{-\frac{j}{2}}.
		\end{align*}
	\end{itemize}
	
\end{proof}

\vspace{.2cm}

\section{Esitmates for non-zero modes}\label{Sec-md}

We now establish the non-zero mode estimate, \cref{ineq-md*}, to finish the proof of Step 3 in \cref{Sec-steps}.
Recall the effective variable, $\wv=\psi - \e \uvt \phi,$ defined in \eqref{wv}, and note that, for $i=1,2,3,$
\begin{align*}
	\frac{m_i \mv}{\rho} - \frac{\mt_i \mvt}{\rhot} & = \ut_i \psi + \zeta_i \mv = \ut_i \wv + \e \ut_i \uvt \phi + \frac{w_i}{\rho}  \big(\mvt + \wv + \e \uvt \phi\big) \\
	& = \ut_i \wv  + \uvt w_i + \frac{1}{\rho} w_i \wv + \e \ut_i \uvt \phi.
\end{align*}
Then using \cref{equ-phipsi,wv}, one can get the system of $(\phi,\wv)$:
\begin{equation}\label{equ-phiw}
	\begin{cases}
		\p_t \phi + \uvt \dnab \phi + \e^{-1} \dv \wv + \p_3 \ut_3 \phi = - \p_3 F_0, \\
		\p_t \wv + \uvt\dnab \wv + \e^{-1} p'(\rhot) \nabla \phi + \Lv + \Qv_3 + \e \Qv_4 \\ 
		\qquad\qquad 
		= \mu \lap \zeta + (\mu+\lambda) \nabla\dv \zeta - \p_3 \Fv + \e \uvt \, \p_3 F_0,
	\end{cases}
\end{equation}
where $\Lv $ and $\Qv_i$ for $i=3,4$ are the linear and nonlinear terms, given by
\begin{equation}\label{Q-4}
	\begin{aligned}
		\Lv & := \big( \e \p_t \uvt + \e \ut_3 \p_3 \uvt + \e^{-1} p''(\rhot) \nabla\rhot \big) \phi +  \p_3 \ut_3 \wv + \p_3 \uvt \, w_3, \\
		\Qv_3 & := \frac{1}{\rhot} \big(\wv \dv \wv + \wv \dnab \wv\big) + \e^{-2} \nabla\varpi(\rho, \rhot), \\
		\Qv_4 & := - \frac{\phi}{\rhot \rho} \big(\wv \dv \wv + \wv \dnab \wv\big) - \frac{1}{\rho^2} \big( \e^{-1} \nabla\rhot + \nabla \phi \big) \cdot 
		\wv \wv.
	\end{aligned}
\end{equation}
Correspondingly, the system of the non-zero mode, $ (\phi^\md, \wv^\md), $ is given as follows,
\begin{equation}\label{equ-md}
	\begin{cases}
		\p_t \phi^\md+ \uvt\dnab \phi^\md + \e^{-1} \dv \wv^\md  + \p_3 \ut_3 \phi^\md = 0, \\
		\p_t \wv^\md + \uvt \dnab \wv^\md + \e^{-1}   p'(\rhot) \nabla \phi^\md + \Lv^\md + \Qv_3^\md + \e \Qv_4^\md \\
		\qquad\qquad = \mu \lap \zeta^\md + (\mu+\lambda) \nabla \dv \zeta^\md.
	\end{cases}
\end{equation}

To establish the $L^2$-estimate for the non-zero mode system \cref{equ-md}, we encounter difficulties arising from the convection in both the linear and nonlinear levels.

\subsection{Difficulty in the linear level}

It follows from \cref{Q-4}$_1$ that 
\begin{equation*} 
	\Lv^\md = \big( \e \p_t \uvt + \e \ut_3 \p_3 \uvt + \e^{-1} p''(\rhot) \nabla\rhot \big) \phi^\md +  \p_3 \ut_3 \wv^\md + \p_3 \uvt \, w_3^\md.
\end{equation*}
The linear difficulty comes from the last term, $\p_3 \uvt \, w_3^\md,$ which is not able to be controlled, if the tangential part, $\abs{\p_3 \uvt_\perp} \lesssim \abs{\uvb} (t+\Lambda)^{-\frac{1}{2}},$ has neither smallness nor a sufficiently fast decay rate. 
To overcome this difficulty, we choose the constant $\Lambda$ in \cref{profile-t} to be suitably large such that this bad term can be controlled by the dissipation.
In fact, one can use \cref{profile-t,bdd-diff,poincare} to obtain the following linear estimate.

\begin{Lem}\label{Lem-L-md}
	It holds that
	\begin{equation}\label{large-Lam}
		\begin{aligned}
			\norm{\Lv^\md}_{H^1} & \lesssim (\e+\chi) (t+1)^{-1} \norm{\nabla \phi^\md}_{L^2} + (t+\Lambda)^{-\frac{1}{2}} \norm{\nabla \wv^\md}_{L^2}.
		\end{aligned}
	\end{equation}
\end{Lem}

\vspace{.1cm}

\subsection{Difficulty in the nonlinear level}

To deal with the nonlinear terms in \cref{equ-md}, we first decompose $\Qv_3^\md$ into $\Qv_3^\md = \widecheck{\Qv}_3^\md + \widehat{\Qv}_3^\md, $ where
\begin{equation}\label{Q-md}
	\begin{aligned}
		\widecheck{\Qv}_3^\md & =  \frac{1}{\rhot} \big[ \wv^\od \dv \wv^\md + w_3^\md \p_3 \wv^\od + \wv^\od \dnab \wv^\md  \big], \\
		\widehat{\Qv}_3^\md & = \frac{1}{\rhot} \big[ \p_3 w_3^\od \wv^\md  + \big( \wv^\md \dv \wv^\md + \wv^\md \dnab \wv^\md \big)^\md \big] + \e^{-2} \nabla \big(\varpi(\rho, \rhot)\big)^\md.
	\end{aligned}
\end{equation}
In fact, it is noted that the nonlinear terms, $\widehat{\Qv}_3^\md$ and $\e \Qv_4^\md$ in \cref{equ-md}, contain small parameters $\nu$ and $\e,$ respectively, while $\widecheck{\Qv}_3^\md$ is indeed large due to the zero mode associated with the tangential velocity, $\wv_\perp^\od$.

\begin{Lem}[Nonlinearity with smallness]\label{Lem-Q3h4}
	If $ \big(\nu^{\frac{1}{6}} + \e^{\frac{1}{4}} \big) M \leq 1, $ then it holds that 
	\begin{equation}\label{est-Q3h4}
		\begin{aligned}
			\norm{\widehat{\Qv}_3^\md}_{L^2} + \norm{\e \Qv_4^\md}_{L^2} & \lesssim (\nu^{\frac{1}{8}}+ \e^{\frac{1}{2}}) (t+1)^{-\frac{1}{2}} \norm{\nabla \big(\phi^\md,\wv^\md \big)}_{L^2},
		\end{aligned}
	\end{equation}
\end{Lem}
\begin{proof}
	Since both $\rhot$ and $\phi^\od$ are independent of $\xp,$ then it holds that
	\begin{align}
		\big( \varpi(\rho, \rhot) \big)^\md & = \big( p(\rhot + \e \phi^\od + \e \phi^\md) - p(\rhot + \e \phi^\od) \big)^\md - \e p'(\rhot) \phi^\md \notag \\
		& = \e \Big\{ \int_0^1 \big[ p'(\rhot+\e \phi^\od + \e r\phi^\md) - p'(\rhot) \big] dr \, \phi^\md \Big\}^\md. \label{varpi-md}
	\end{align}
	Then the remaining proof is based on the a priori bounds shown in \cref{Sec-bound}.
	In fact ,it follows \cref{aLinf-od,aLinf-md} that  
	\begin{align*}
		\norm{\widehat{\Qv}_3^\md}_{L^2} & \lesssim \big(\norm{\p_3 w_3^\od}_{L^{\infty}} + \norm{\wv^\md}_{L^\infty} \big) \norm{ \wv^\md}_{H^1} + \big(\norm{\phi^\od}_{W^{1,\infty}} + \norm{\phi^\md}_{L^\infty} \big) \norm{ \phi^\md}_{H^1} \\
		&  \lesssim M^{\frac{3}{4}} \nu^{\frac{1}{4}} (t+1)^{-\frac{1}{2}} \norm{ \nabla (\wv^\md, \phi^\md)}_{L^2}.
	\end{align*}
	Similarly, one can use \cref{aLinf-pert} to obtain that
	\begin{align*}
		\norm{\e \Qv_4^\md}_{L^2} & \lesssim \e M^2 (t+1)^{-\frac{1}{2}} \norm{(\phi^\md, \wv^\md)}_{H^1} \\
		& \lesssim 
		\e M^2 (t+1)^{-\frac{1}{2}} \norm{\nabla \big(\phi^\md, \wv^\md \big)}_{L^2}.
	\end{align*}
	Thus, \cref{est-Q3h4} holds true if $ M^{\frac{3}{4}} \nu^{\frac{1}{8}} + \e^{\frac{1}{2}} M^2 \leq 1.$
\end{proof}

\vspace{.1cm}

\begin{Lem}[Nonlinearity without smallness]\label{Lem-Q3c}
	If $\big(\chi^{\frac{1}{2}} + \nu^{\frac{1}{6}} + \e^{\frac{1}{4}} \big) M \leq 1, $ then it holds that 
	\begin{equation}\label{est-Q3c}
		\begin{aligned}
			\norm{\widecheck{\Qv}_3^\md}_{L^2} & \lesssim \big[ (t+1)^{-\frac{3}{4}} + \big(\chi^{\frac{1}{2}} + \e^{\frac{1}{2}} + \nu^{\frac{1}{8}}\big) (t+1)^{-\frac{1}{2}} \big] \norm{\nabla \wv^\md}_{L^2}.
		\end{aligned}
	\end{equation}
\end{Lem}
\begin{proof}
	It follows from \cref{poincare} and \cref{Q-md}$_1$ that
	$$
	\norm{\widecheck{\Qv}_3^\md}_{L^2} \lesssim \norm{\wv^\od}_{W^{1,\infty}} \norm{\wv^\md}_{H^1} \lesssim \norm{\wv^\od}_{W^{1,\infty}} \norm{\nabla \wv^\md}_{L^2}.
	$$
	Recall that $\wv^\od = \p_3 \Zv + \e \p_3 \uvt \Phi.$ 
	Using \cref{Zp-inf,aLinf-anti}, if $ \big(\chi^{\frac{1}{2}} + \nu^{\frac{1}{6}} + \e^{\frac{1}{4}} \big) M \leq 1,$ then
	\begin{align*}
		\norm{\wv^\od}_{W^{1,\infty}} & \lesssim \norm{\p_3 \Zv_\perp}_{W^{1,\infty}} + \norm{\p_3 Z_3}_{W^{1,\infty}} + \e (t+1)^{-\frac{1}{2}} \norm{\Phi}_{W^{1,\infty}} \\
		& \lesssim (t+1)^{-\frac{3}{4}} + \big(\chi^{\frac{1}{2}} + \e^{\frac{1}{2}} + \nu^{\frac{1}{8}}\big) (t+1)^{-\frac{1}{2}},
	\end{align*}
	which yields \cref{est-Q3c}. 
	
\end{proof}

Unfortunately, the estimate \cref{est-Q3c} is still not sufficient to establish the zero-order $L^2$-estimate for the non-zero mode system \cref{equ-md}, owing to the large size of initial perturbation as well as the slow rate $(t+1)^{-\frac{3}{4}}.$
Thus, we shall take an additional care of the large nonlinearity, $\widecheck{\Qv}_3^\md,$ and the following estimate plays the most important role in the non-zero mode estimates.

\begin{Lem}[Key estimate]\label{Lem-Q-key}
	If $\big(\chi^{\frac{1}{2}} + \nu^{\frac{1}{6}} + \e^{\frac{1}{4}} \big) M \leq 1, $ then it holds that	
	\begin{equation}\label{Q-34m}
		\begin{aligned}
			& \abs{\int_\Omega \wv^\md \cdot \widecheck{\Qv}_3^\md dx + \frac{d}{dt} \Big( \int_\Omega \frac{\e}{\rhot} \wv^\od \cdot \wv^\md \phi^\md dx \Big)} \\
			& \qquad \lesssim (t+1)^{-\frac{5}{4}} \norm{\wv^\md}_{L^2}^2 + (\chi^\frac{1}{2} +\e^\frac{1}{2}+\nu^\frac{1}{8}) \norm{\nabla\big(\phi^\md, \wv^\md\big)}_{L^2}^2.
		\end{aligned}
	\end{equation}
\end{Lem}

\begin{proof}

It follows from \cref{Q-md}$_1$ that
\begin{align*}
	\wv^\md \cdot \widecheck{\Qv}_3^\md & = \dv\Big( \frac{\wv^\od}{2\rhot} \abs{\wv^\md}^2 \Big) + \sum_{i=1}^{3} I_i,
\end{align*}
where 
\begin{equation}\label{I-Q}
	\begin{aligned}
		I_1 & = \frac{1}{\rhot} \wv^\od \cdot \wv^\md \dv \wv^\md, \\
		I_2 & = \frac{1}{\rhot} \p_3 w_3^\od \abs{w_3^\md}^2 - \p_3 \Big( \frac{w_3^\od}{2\rhot} \Big) \abs{\wv^\md}^2, \\
		I_3 & =  \frac{1}{\rhot} \p_3 \wv_\perp^\od \cdot \wv_\perp^\md w_3^\md.
	\end{aligned}
\end{equation}
Now we estimate each term in \cref{I-Q}.
~\\
1) We first claim that
\begin{equation}\label{I-Q-1}
	\left\lvert\int_\Omega I_1 dx + \frac{d}{dt} \Big( \int_\Omega \frac{\e}{\rhot} \wv^\od \cdot \wv^\md \phi^\md dx \Big)\right\rvert \lesssim (\nu^{\frac{1}{8}} + \e^{\frac{1}{2}}) \norm{\nabla (\phi^\md, \wv^\md)}_{L^2}^2.
\end{equation}
In fact, it follows from \cref{equ-md}$_1$ that
\begin{equation}\label{I1}
	\begin{aligned}
		& I_1 + \p_t \Big( \frac{\e}{\rhot} \wv^\od \cdot \wv^\md \phi^\md \Big) \\
		& \qquad = I_{1,1} + I_{1,2} + O(1) \e \norm{\wv^\od}_{L^\infty} \abs{\wv^\md} \big[\chi (t+1)^{-1} \abs{\phi^\md} + \abs{\nabla\phi^\md} \big].
	\end{aligned}
\end{equation}
where 
$$
I_{1,1} = \frac{\e}{\rhot}  \p_t \wv^\od \cdot \wv^\md \phi^\md, \qquad I_{1,2} = \frac{\e}{\rhot} \wv^\od \cdot \p_t \wv^\md  \phi^\md.
$$
Since $\wv^\od = \psi^\od - \e \uvt \phi^\od $, then it holds that
\begin{equation}\label{ptw}
	\begin{aligned}
		\e \norm{\p_t \wv^\od}_{L^\infty} 
		& \lesssim  \e \norm{\p_t \psi^\od + \e^{-1} \p_3 \big(p'(\rhot) \phi^\od\big) \E_3}_{L^\infty} + \norm{\p_3 \big(p'(\rhot) \phi^\od\big)}_{L^\infty} \\
		&\quad  + \e^2 (t+1)^{-1} \norm{\phi^\od}_{L^\infty} + \e^2 \norm{ \p_t \phi^\od}_{L^\infty}.
	\end{aligned}
\end{equation}
It follows from \cref{aLinf-od} that
\begin{equation}\label{ptw-1}
	\norm{\p_3 \big(p'(\rhot) \phi^\od\big)}_{L^\infty} + \e^2 (t+1)^{-1} \norm{\phi^\od}_{L^\infty} \lesssim (\nu^{\frac{1}{4}} M^{\frac{3}{4}} + \e \nu) (t+1)^{-\frac{3}{4}}.
\end{equation}
Besides, the zero-mode system \cref{equ-pert-od}, together with the bounds shown in \cref{Sec-bound}, yields that
\begin{align}
	& \e \norm{\p_t \psi^\od + \e^{-1} \p_3 \big(p'(\rhot) \phi^\od\big) \E_3}_{L^\infty}  + \e^2 \norm{\p_t \phi^\od}_{L^\infty} \notag \\ 
	& \qquad \lesssim \e  \norm{(\phi, \psi, \zeta)}_{L^\infty} \big( 1+ \norm{\p_3 (\phi, \psi, \zeta)}_{L^\infty} \big) + \e \norm{\p_3^2 \zeta^\od}_{H^1} + \e \chi \lesssim \e M^2. \label{ptw-2}
\end{align}
Thus, collecting \cref{ptw,ptw-1,ptw-2}, if $ \nu^{\frac{1}{8}} M^{\frac{3}{4}} + \e^{\frac{1}{2}} M^2 \leq 1,$ it holds that
\begin{equation}\label{bdd-ptw}
	\e \norm{\p_t \wv^\od}_{L^\infty}  \lesssim \nu^{\frac{1}{8}} + \e^{\frac{1}{2}},
\end{equation}
which, together with \cref{poincare}, yields that
\begin{equation}\label{small-1}
	\abs{\int_\Omega I_{1,1} dx} \lesssim \big(\nu^{\frac{1}{8}} + \e^{\frac{1}{2}}\big) \norm{\nabla(\phi^\md, \wv^\md)}_{L^2}.
\end{equation}

For the second nonlinear term, $I_{1,2}, $ in \cref{I1}, it follows from \cref{equ-md}$_2$ that
\begin{align*}
	I_{1,2} & = \frac{\e}{\rhot} \phi^\md\, \wv^\od \cdot \big(\p_t \wv^\md + \e^{-1} p'(\rhot) \nabla \phi^\md\big) - \frac{p'(\rhot)}{\rhot} \phi^\md\,  \wv^\od \cdot \nabla \phi^\md \\
	& = \dv (\cdots) + \p_3 \Big( \frac{p'(\rhot)}{2 \rhot} w_3^\od \Big) \abs{\phi^\md}^2 \\
    & \quad + O(1) \e \norm{\wv^\od}_{L^\infty} \abs{\phi^\md} \big(\abs{\nabla\wv^\md} + \abs{\Lv^\md} + \abs{\Qv_3^\md} + \e \abs{\Qv_4^\md} \big)  \\
	& \quad + O(1) \e \norm{\wv^\od}_{W^{1,\infty}}\abs{\nabla\zeta^\md} \big(\abs{\nabla\phi^\md} + \abs{\phi^\md} \big).
\end{align*}
Using Lemmas \ref{Lem-L-md}--\ref{Lem-Q3c} and the bounds shown in \cref{Sec-bound}, one has that if $\nu^{\frac{1}{8}} M^{\frac{3}{4}} + \e^{\frac{1}{2}} M^2 \leq 1, $ then
\begin{align}
	\abs{\int_\Omega I_{1,2} dx} & \lesssim M^{\frac{3}{4}} \nu^{\frac{1}{4}} \norm{\phi^\md}_{L^2}^2 + \e M^2 \norm{\nabla \big(\phi^\md, \wv^\md\big)}_{L^2}^2 \notag \\
	& \lesssim (\nu^{\frac{1}{8}} + \e^{\frac{1}{2}}) \norm{\nabla (\phi^\md, \wv^\md)}_{L^2}^2. \label{small-2}
\end{align}
The last term on the right-hand side of \cref{I1} is similar to estimate. Then combining \cref{small-1,small-2}, one can obtain \cref{I-Q-1}.

2) For $I_2$ in \cref{I-Q}, using \cref{aLinf-od}, if $M^{\frac{3}{4}} \nu^{\frac{1}{8}} \leq 1, $ it holds that
\begin{equation}\label{I-Q-2}
	\abs{\int_\Omega I_2 dx} \lesssim \norm{ w_3^\od}_{W^{1,\infty}} \norm{ \wv^\md}_{L^2}^2 \lesssim \nu^{\frac{1}{8}} \norm{\nabla  \wv^\md}_{L^2}^2.
\end{equation}

3) Note that $ \p_3 \wv_\perp^\od = \p_3^2 \Zv_\perp + \e \p_3 \big( \uvt_\perp \p_3 \Phi\big). $ Then using \cref{aLinf-anti} and \cref{Lem-Zp-inf}, the $I_3$  in \cref{I-Q} satisfies that
\begin{equation}\label{I-Q-3}
	\abs{\int_\Omega I_3 dx} \lesssim (t+1)^{-\frac{5}{4}} \norm{\wv^\md}_{L^2}^2 + (\chi^\frac{1}{2}+\nu^\frac{1}{4} +\e^\frac{1}{2}) \norm{\nabla \wv^\md}_{L^2}^2.
\end{equation}

Collecting \cref{I-Q-1,I-Q-2,I-Q-3}, one can finish the proof.

\end{proof}

\vspace{.1cm}

\subsection{Proof of Step 3}
With the lemmas shown above, we are able to establish the $L^2$-estimates for the non-zero system \cref{equ-md}. These estimates complete the proof of Step 3 in \cref{Sec-steps}.

\begin{Lem}\label{Lem-md-1}
	If $\big(\chi^{\frac{1}{2}} + \nu^{\frac{1}{6}} + \e^{\frac{1}{4}} \big) M \leq 1, $  then
	\begin{equation}\label{est-md-1*}
		\begin{aligned}
			\frac{d}{dt} \Ac_{0,\sharp}  + \Bc_{0,\md} & \lesssim (t+1)^{-\frac{5}{4}} \norm{\wv^\md}_{L^2}^2 + \e \norm{\nabla^2 \phi^\md}^2_{L^2},
		\end{aligned}
	\end{equation}
	where
	\begin{align*}
		\Ac_{0,\sharp} \sim \norm{\big(\phi^\md, \wv^\md\big)}_{L^2}^2 + \e^2 \norm{\nabla \phi^\md}_{L^2}^2 \andd \Bc_{0,\md} \sim \norm{\nabla (\phi^\md, \wv^\md)}_{L^2}^2.
	\end{align*}
\end{Lem}
\begin{proof}
	Multiplying $ p'(\rhot) \phi^\md $ and $ \wv^\md $ on \cref{equ-md}$ _1 $ and \cref{equ-md}$ _2 $, respectively, one can get that
	\begin{align*}
		& \p_t \Big( \frac{p'(\rhot)}{2} \abs{\phi^\md}^2 + \frac{1}{2} \abs{\wv^\md}^2 \Big) + \frac{\mu }{\rhot} \abs{\nabla\wv^\md}^2 + \frac{\mu+\lambda}{\rhot} \abs{\dv\wv^\md}^2 = \dv(\cdots) + \sum_{i=4}^{6} I_i,
	\end{align*}
	where 
	\begin{equation*}
		\begin{aligned}
			I_4 & = \frac{1}{2} \big[ p''(\rhot) \big(\p_t \rhot + \ut_3 \p_3 \rhot \big) - p'(\rhot) \p_3 \ut_3 \big] \abs{\phi^\md}^2 + \e^{-1} p''(\rhot) \p_3 \rhot \phi^\md w_3^\md + \frac{1}{2} \p_3 \ut_3 \abs{\wv^\md}^2, \\
			I_5 & = - \wv^\md \cdot \big( \Lv^\md + \widecheck{\Qv}_3^\md + \widehat{\Qv}_3^\md + \e \Qv^\md_4 \big), \\
			I_6 & = \mu \sum_{i=1}^{3} \p_i \wv^\md \cdot \Big( \frac{1}{\rhot} \p_i \wv^\md - \p_i \zeta^\md \Big) + (\mu+\lambda) \dv \wv^\md \Big( \frac{1}{\rhot} \dv \wv^\md - \dv \zeta^\md \Big).
		\end{aligned}
	\end{equation*} 
	First, one can use \cref{bdd-diff} to get that
	\begin{equation}\label{I-4}
		\abs{\int_\Omega I_4 dx} \lesssim (\chi + \e) \norm{ (\phi^\md, \wv^\md)}_{L^2}^2 \lesssim (\chi + \e) \norm{\nabla (\phi^\md, \wv^\md)}_{L^2}^2.
	\end{equation}
	The combination of Lemmas \ref{Lem-L-md}, \ref{Lem-Q3h4} and \ref{Lem-Q-key} yields that
	\begin{equation}\label{I-5}
		\begin{aligned}
			& \abs{\int_\Omega I_5 dx - \frac{d}{dt} \Big( \int_\Omega \frac{\e}{\rhot} \wv^\od \cdot \wv^\md \phi^\md dx \Big)} \\
			& \qquad \lesssim (t+1)^{-\frac{5}{4}} \norm{\wv^\md}_{L^2}^2 + \big(\Lambda^{-\frac{1}{2}} + \chi^{\frac{1}{2}} +\nu^{\frac{1}{8}} + \e^{\frac{1}{2}}\big) \norm{\nabla (\phi^\md, \wv^\md)}_{L^2}^2.
		\end{aligned}
	\end{equation}
	To estimate $I_6$, it follows from \cref{rel-zeta-w} that if $\e^{\frac{1}{2}} M \leq 1, $ then
	$$
	\norm{\zeta^\md}_{L^2} \lesssim \norm{\wv^\md}_{L^2} + \e^{\frac{1}{2}} \norm{\phi^\md}_{L^2}.
	$$
	Then using \cref{wv-md}, one can get that
	\begin{equation*}
		\norm{\nabla \wv^\md - \rhot\, \nabla \zeta^\md}_{L^2}  \lesssim \e M \norm{(\phi^\md, \zeta^\md)}_{H^1} 
		\lesssim \e^{\frac{1}{2}} \norm{\nabla (\phi^\md, \wv^\md)}_{L^2},
	\end{equation*}
	which yields that
	\begin{equation}\label{I-6}
		\abs{\int_\Omega I_6 dx} \lesssim \e^{\frac{1}{2}} \norm{\nabla (\phi^\md, \wv^\md)}_{L^2}^2.
	\end{equation}
	
	Collecting \cref{I-4,I-5,I-6} and choosing the constant $\Lambda$ to be suitably large (see \cref{Lambda}), one can get that
	\begin{equation}\label{ineq3}
		\begin{aligned}
			& \frac{d}{dt} \Big[ \int_\Omega \big( p'(\rhot) 	\abs{\phi^\md}^2 + \abs{\wv^\md}^2 \big) dx - 2\e \int_\Omega \frac{1}{\rhot} \wv^\od \cdot \wv^\md  \phi^\md dx \Big] + \frac{\mu}{4\rhob}  \norm{\nabla\wv^\md}_{L^2}^2 \\
			& \qquad \lesssim (t+1)^{-\frac{5}{4}} \norm{\wv^\md}_{L^2}^2 + (\chi^{\frac{1}{2}} +\nu^{\frac{1}{8}} + \e^{\frac{1}{2}}) \norm{\nabla \phi^\md}_{L^2}^2.
		\end{aligned}
	\end{equation}

\vspace{.1cm}

On the other hand, to achieve the density dissipation, one can multiply $ \frac{\e^2 \mut}{\rhot} \nabla \phi^\md $ and $ \e \nabla \phi^\md $ on $\nabla$\cref{equ-md}$ _1 $ and \cref{equ-md}$ _2 $, respectively, to get that
\begin{equation}
	\begin{aligned}
		& \p_t \Big( \frac{\e^2 \mut}{2\rhot} \abs{\nabla\phi^\md}^2 + \e \wv^\md \dnab \phi^\md \Big) + p'(\rhot) \abs{\nabla\phi^\md}^2 \\
	& \qquad = \dv(\cdots) + \abs{\dv\wv^\md}^2 + O(1) \e M^2 \abs{\nabla\phi^\md} \Big[ \sum_{i=0}^{1 } \big( \abs{\nabla^i \phi^\md} + \abs{\nabla^i \wv^\md} \big) \Big] \\
	& \qquad \quad + O(1) \e \abs{\nabla^2 \phi^\md} \big(\abs{\nabla\wv^\md} + \abs{\nabla\zeta^\md} \big).
	\end{aligned}
\end{equation}
This, together with \cref{rel-zeta-w}, yields that if $\e^{\frac{1}{2}} M^2 \leq 1,$ then 
\begin{equation}\label{ineq4}
	\begin{aligned}
		& \frac{d}{dt} \Big(\int_\Omega \frac{\e^2 \mut}{\rhot} \abs{\nabla\phi^\md}^2 dx + \int_\Omega 2 \e \wv^\md \dnab \phi^\md dx \Big) + \int_\Omega p'(\rhot) \abs{\nabla\phi^\md}^2 dx \\
		& \qquad \lesssim \norm{\nabla \wv^\md}_{L^2}^2 + \e \norm{ \nabla^2 \phi^\md}_{L^2}^2. 
	\end{aligned}
\end{equation}

At last, the combination of \cref{ineq3,ineq4} can imply \cref{est-md-1*}.

\end{proof}


\begin{Lem}\label{Lem-md-2}
	If $\big(\chi^{\frac{1}{2}} + \nu^{\frac{1}{6}} + \e^{\frac{1}{4}} \big) M \leq 1, $ then
	\begin{equation}\label{est-md-2}
		\begin{aligned}
			& \frac{d}{dt} \big[(t+1) \Ac_{1,\sharp}\big] + (t+1) \Bc_{1,\sharp}  \lesssim \Bc_{0,\md} + \e^{\frac{1}{2}} (t+1) \norm{\nabla^2 \phi^\md}_{L^2}^2,
		\end{aligned}
	\end{equation}
	where
	\begin{align*}
		\Ac_{1,\sharp} \sim \norm{\nabla (\phi^\md, \wv^\md)}^2_{L^2}  \andd \Bc_{1,\sharp} \sim \norm{\nabla^2  \wv^\md}_{L^2}^2
	\end{align*}
\end{Lem}
\begin{proof}
	Multiplying $ p'(\rhot) \nabla \phi^\md $ and $ -\lap \wv^\md $ on $ \nabla $\cref{equ-md}\textsubscript{1} and \cref{equ-md}$ _2 $, respectively, one can get that
	\begin{equation}\label{eq-1}
		\begin{aligned}
			& \p_t \Big( \frac{p'(\rhot)}{2} \abs{\nabla\phi^\md}^2 \Big) + \e^{-1} p'(\rhot) \nabla \dv \wv^\md \dnab \phi^\md \\
			& \qquad = \dv(\cdots) + O(1) \abs{ \nabla\phi^\md} \big( \abs{ \nabla\phi^\md} + \abs{ \phi^\md} \big),
		\end{aligned}
	\end{equation}
	and
	\begin{equation}\label{eq-2}
		\begin{aligned}
			& \p_t \Big( \frac{1}{2} \abs{\nabla\wv^\md}^2 \Big) - \e^{-1} p'(\rhot) \nabla \dv \wv^\md \dnab \phi^\md +  \frac{\mu}{\rhot} \abs{\lap \wv^\md}^2 + \frac{\mu+\lambda}{\rhot} \abs{\nabla\dv \wv^\md}^2 \\
			& \qquad = \dv(\cdots) + \sum_{i=7}^{9} I_i + O(1) \Big[ \sum_{i=0}^{1} \big(\abs{\nabla^i \phi^\md}^2 + \abs{\nabla^i \wv^\md}^2\big) \Big],
		\end{aligned}
	\end{equation}
	where 
	\begin{equation}\label{I7-9}
		\begin{aligned}
			I_7 & = \e^{-1} p'(\rhot) \nabla  \phi^\md \cdot \big(\lap \wv^\md - \nabla \dv \wv^\md \big) -  (\mu+\lambda) \nabla\dv \wv^\md \cdot \big( \lap \wv^\md - \nabla \dv \wv^\md \big), \\
			I_8 & = \lap \wv^\md \cdot \big(\widecheck{\Qv}_3^\md + \widehat{\Qv}_3^\md + \e \Qv_4^\md\big), \\
			I_9 & = \mu \lap\wv^\md \cdot \Big(\frac{1}{\rhot}\lap\wv^\md - \lap \zeta^\md\Big) + (\mu+\lambda) \lap \wv^\md \cdot \Big(\frac{1}{\rhot} \nabla\dv\wv^\md - \nabla\dv\zeta^\md\Big).
		\end{aligned}
	\end{equation}
	
	1) Note that for any functions $ h_0\in H^1, \mathbf{h}=(h_1,h_2,h_3) \in H^2, $ it holds that
	\begin{equation}\label{fact}
		\nabla h_0 \cdot (\lap \mathbf{h} - \nabla\dv \mathbf{h}) = \dv (\nabla h_0 \times \mathsf{curl}\, \mathbf{h}).
	\end{equation} 
	Using this fact on \cref{I7-9}$_1$ yields that
	$$I_7 = \dv(\cdots) - \e^{-1} p''(\rhot) \nabla\rhot \cdot \big( \nabla \phi^\md \times \text{curl}\, \wv^\md \big). $$
	Then it holds that
	\begin{equation}\label{I-7}
		\abs{\int_\Omega I_7 dx} \lesssim \chi (t+1)^{-1} \norm{\nabla (\phi^\md, \wv^\md)}_{L^2}^2.
	\end{equation}
	
	2) Applying Lemmas \ref{Lem-Q3h4} and \ref{Lem-Q3c} on \cref{I7-9}$_2$, one can get that
	\begin{align}
		\abs{\int_\Omega I_8 dx} & \lesssim (t+1)^{-\frac{1}{2}} \norm{\nabla^2 \wv^\md}_{L^2} \norm{\nabla (\phi^\md, \wv^\md)}_{L^2} \notag \\
		& \leq \frac{\mu}{100 \rhob} \norm{\nabla^2 \wv^\md}_{L^2}^2 + C (t+1)^{-1} \norm{\nabla (\phi^\md, \wv^\md)}_{L^2}^2. \label{I-8}
	\end{align}

	3) To estimate $I_9,$ it follows from \cref{wv-md,rel-zeta-w} that
		\begin{equation*}
			\begin{aligned}
				& \norm{\lap \wv^\md - \rhot \lap \zeta^\md}_{L^2} + \norm{\nabla\dv \wv^\md - \rhot\,\nabla\dv \zeta^\md}_{L^2} \\
				& \qquad \lesssim \norm{\p_3 \rhot}_{W^{1,\infty}} \norm{\zeta^\md}_{H^1} + \e \big( \norm{(\phi^\od, \zeta^\od)}_{W^{2,\infty}} + \norm{(\phi^\md,\zeta^\md)}_{W^{1,\infty}} \big) \norm{ (\phi^\md, \zeta^\md)}_{H^2} \\
				& \qquad \lesssim (\e \chi + \e M) \norm{\nabla^2 (\phi^\md, \zeta^\md)}_{L^2} \\
				& \qquad \lesssim \e^{\frac{1}{2}} \norm{\nabla^2 (\phi^\md, \wv^\md)}_{L^2}, \qquad \text{if } \ \e^{\frac{1}{2}} M \leq 1.
			\end{aligned}
	\end{equation*}
	Then it holds that
	\begin{equation}\label{I-9}
		\abs{\int_\Omega I_9 dx} \leq \frac{\mu}{100 \rhob} \norm{\nabla^2 \wv^\md}_{L^2}^2 + C \e \norm{\nabla^2 \phi^\md}_{L^2}^2.
	\end{equation}
 	
 	Note that $\norm{\lap \wv^\md}_{L^2}^2 = \norm{\nabla^2 \wv^\md}_{L^2}^2 $, which follows from integration by parts. Then adding \cref{eq-1,eq-2} up, and using \cref{I-7,I-8,I-9}, one can get \cref{est-md-2}.

\end{proof}


Combining Lemmas \ref{Lem-md-1} and \ref{Lem-md-2}, one can complete the proof of \cref{ineq-md*}.

\vspace{.2cm}

\section{Higher order estimates}\label{Sec-orn}

Now we show the proof of Step 4 in \cref{Sec-steps}. 
We consider the system associated with the density and velocity, that is, the one satisfied by the perturbations,
\begin{equation*}
	(\phi, \zeta) = \big( \e^{-1} (\rho-\rhot), \uv-\uvt \big).
\end{equation*}
The benefit of the consideration is to avoid the complex nonlinearity in the system of $(\phi,\,\psi)$; see \cref{equ-phipsi}. 
In fact, the system of $(\phi,\,\zeta)$, derived from \cref{NS,equ-ansatz}, is given by
\begin{equation}\label{equ-phizeta}
	\begin{cases}
		\p_t  \phi + \uv \cdot \nabla \phi + \e^{-1} \rho \dv \zeta + \widecheck{L}_0 = 0, \\
		\p_t \zeta + \uv \dnab \zeta + \e^{-1} \frac{p'(\rho)}{ \rho}  \nabla \phi + \widecheck{\Lv} = \frac{\mu}{\rho} \lap \zeta + \frac{\mu+\lambda}{\rho} \nabla \dv \zeta,
	\end{cases}
\end{equation}
where
\begin{equation}\label{Lc}
\begin{aligned}
	\widecheck{L}_0 & := \e^{-1} \p_3 \rhot \zeta_3 + \p_3 \ut_3 \phi + \p_3 F_0, \\
	\widecheck{\Lv} & = (\widecheck{L}_1, \widecheck{L}_2, \widecheck{L}_3) := \p_3 \uvt \zeta_3 + \db_2 \phi + \frac{1}{\rhot }\big(\p_3 \Fv - \e \p_3 F_0 \uvt\big),
\end{aligned}
\end{equation}
with 
\begin{equation*}
	\db_2 = \e^{-1} \p_3\rhot \int_0^1 \frac{d}{ds} \Big(\frac{p'(s)}{s}\Big)|_{s=\rhot + \e r \phi} dr \, \E_3
	+ \frac{\e}{\rhot\rho} \big[ \mu \p_3^2 \uvt + (\mu+\lambda) \p_3^2 \ut_3 \E_3 \big].
\end{equation*}

\begin{Lem}
	If $\e M \leq 1, $ for $ k=1,2,3, $ it holds that
	\begin{equation}\label{est-Lc}
		\begin{aligned}
			\norm{\nabla^k \widecheck{L}_0}_{L^2} & \lesssim \chi \sum_{j=0}^{k} (t+1)^{-\frac{k-j+2}{2}} \norm{\nabla^j (\phi, \zeta)}_{L^2} + \chi (t+1)^{-\frac{2k+5}{4}},  \\
			\norm{\nabla^k \widecheck{\Lv}}_{L^2} & \lesssim \sum_{j=0}^{k} \big[(t+\Lambda)^{-\frac{k-j+1}{2}} \norm{\nabla^j \zeta}_{L^2} + (\chi+\e)(t+1)^{-\frac{k-j+2}{2}} \norm{\nabla^j \phi}_{L^2} \big] \\
			& \quad + \chi (t+1)^{-\frac{2k+5}{4}}.
		\end{aligned}
	\end{equation}
\end{Lem}
\begin{proof}
	We only estimate the nonlinear term $\db_2 \phi$, since the other terms are linear and direct to be estimated. 
	Recall that for $ j=0, 1,2,\cdots, $
	\begin{equation*}
		\e^{-1}\abs{\p_3^j \rhot} + \abs{\p_3^j \ut_3} \lesssim \chi (t+1)^{-\frac{1+j}{2}} \andd \abs{\p_3^j \uvt_\perp} \lesssim (t+\Lambda)^{-\frac{j}{2}}.
	\end{equation*}
	Then it holds that
	\begin{align*}
		\abs{\db_2} & \lesssim (\chi+\e)(t+1)^{-1}, \\
		\abs{\nabla \db_2} & \lesssim (\chi+\e) \big[ (t+1)^{-\frac{3}{2}} + \e (t+1)^{-1} \abs{\nabla \phi} \big].
	\end{align*}
	Note that \cref{aLinf-pert} gives $ \norm{\phi}_{W^{1, \infty}} \lesssim M (t+1)^{-\frac{1}{2}}. $
	Thus, if $\e M \leq 1, $ one has that
	\begin{align*}
		\abs{\nabla^2 \db_2} & \lesssim (\chi+\e) \big[ (t+1)^{-2} + \e \sum_{j=1}^{2} (t+1)^{-\frac{4-j}{2}} \abs{\nabla^j \phi} \big], \\
		\abs{\nabla^3 \db_2} & \lesssim (\chi+\e) \big[ (t+1)^{-\frac{5}{2}} + \e \sum_{j=1}^{3} (t+1)^{-\frac{5-j}{2}} \abs{\nabla^j \phi} \big].
	\end{align*}	
	Collecting the estimates above, one can obtain that if $\e M \leq 1, $ 
	\begin{align*}
		\norm{\nabla^k (\db_2 \phi)}_{L^2} 
		& \lesssim (\chi+\e) \Big[ \sum_{j=0}^{k}   (t+1)^{-\frac{k-j+2}{2}} \norm{\nabla^j \phi}_{L^2} \Big] \qquad \text{for }\, k=1,2,3.
	\end{align*}
	
\end{proof}


The rest of the section is devoted to the proof of the high order estimates in \eqref{ineq-ho}. 

\subsection{Proof of Step 4}
This section is devoted to derive the dissipation of the second-order density gradient, which can complete the proof of Step 4 in \cref{Sec-steps}.

The proof of \cref{est-apriori}$_2$, that is, the uniform (in time) smallness of $\nu,$ is based on the estimates \cref{ineq-od*,ineq-md*}. However, on the right-hand side of \cref{ineq-od*}, the second-order term, $\nu^{\frac{1}{2}} (t+1) \norm{\nabla^2 \phi}_{L^2}^2,$ is a large energy with the order $\nu^{\frac{1}{2}}$, while the left-hand side is of the order $\nu^2.$ Thus, we fail to achieve the smallness of $\nu$.  
To overcome this difficulty, we find that in the $L^2$-estimate of the second-order density dissipation, $\nabla^2\phi$, the small Mach number can help control the large energy.


\begin{Lem}\label{Lem-est8}
	If $\e^{\frac{1}{2}} M \leq 1,$ then it holds that 
	\begin{equation}\label{est8}
		\begin{aligned}
			& \frac{d}{dt} \Ac_{2,\e\phi} + \norm{\nabla^2 \phi}_{L^2}^2 \\
			& \qquad \lesssim \norm{\p_3^3 Z_3}_{L^2}^2 +  \norm{\nabla^2 \wv^\md}_{L^2}^2 + \e \norm{\p_3^3 \Zv_\perp}_{L^2}^2 + \e M^2 (t+1)^{-\frac{5}{2}},
		\end{aligned}
	\end{equation}
	where $\Ac_{2,\e\phi}$ is a functional, satisfying that
	\begin{equation}\label{A2}
		\begin{aligned}
			\pm \Ac_{2,\e\phi} \lesssim \pm \e^2 \norm{\nabla^2 \phi}_{L^2}^2 + \e \norm{\nabla^2 \phi}_{L^2} \norm{\nabla \zeta}_{L^2}.
		\end{aligned}
	\end{equation}
\end{Lem}
It is noted that the first two $L^2$-norms on the right-hand side of \cref{est8} are indeed the dissipative terms obtained in \cref{ineq-od*,ineq-md*}.

\begin{proof}[Proof of \cref{Lem-est8}]
	Let $i\in \{1,2,3\} $ be fixed. It follows from \cref{equ-phizeta} that
	\begin{equation}\label{equ-1}
		\begin{aligned}
			& \p_t \nabla\p_i \phi + (\uv \cdot \nabla) \nabla\p_i \phi + \e^{-1} \rho \nabla \dv \p_i \zeta + \nabla \p_i \widecheck{L}_0 \\
			& \qquad = \sum_{j=1}^{2} O(1) \big( \abs{\nabla^j \uv} \abs{\nabla^{3-j}\phi} + \e^{-1} \abs{\nabla^j \rho} \abs{\nabla^{3-j} \zeta} \big),
		\end{aligned}
	\end{equation}
	and
	\begin{equation}\label{equ-2}
		\begin{aligned}
			& \p_t \p_i \zeta +  \frac{p'(\rho)}{\e \rho} \nabla \p_i \phi - \frac{\mut}{\rho} \nabla\dv \p_i\zeta + \p_i \widecheck{\Lv} = \frac{\mu}{\rho} \big(\lap \p_i \zeta - \nabla \dv \p_i \zeta\big) \\
			& \qquad + O(1) \Big[  \sum_{j=0}^{1}  \abs{\nabla^j \uv} \abs{\nabla^{2-j} \zeta} + \e^{-1} \abs{\nabla\rho} \abs{\nabla\phi} + \abs{\nabla\rho} \abs{\nabla^2\zeta} \Big].
		\end{aligned}
	\end{equation}
	Using \cref{equ-phizeta}$_1$, one has that
	\begin{equation}\label{fact1}
		\begin{aligned}
			\e \nabla\p_i \phi \cdot \p_t \p_i \zeta & = \dv(\cdots) + \p_t (\e \nabla\p_i \phi \cdot \p_i \zeta) - \rho \abs{\dv\p_i\zeta}^2 \\
			& \quad - \e \dv \p_i \zeta \big[ \p_i (\uv\dnab \phi) + \e^{-1} \p_i\rho \dv\zeta + \p_i \widecheck{L}_0 \big].
		\end{aligned}
	\end{equation}
	Besides, it follows from \cref{fact} that
	\begin{equation}\label{fact2}
		\nabla\p_i \phi \cdot (\lap \p_i\zeta - \nabla\dv\p_i \zeta) = \dv \big( \nabla\p_i \phi \times \text{curl } \p_i \zeta \big).
	\end{equation}
	Then by multiplying $\frac{\e^2 \mut}{\rho^2} \nabla \p_i\phi $ and $ \e \nabla \p_i \phi $ on \cref{equ-1,equ-2}, respectively, and using \cref{fact1,fact2}, one can get that
	\begin{equation}\label{eq7a}
		\begin{aligned}
			& \p_t\Big( \frac{\e^2 \mut}{2\rho^2} \abs{\nabla\p_i \phi}^2 + \e \nabla \p_i\phi  \cdot \p_i \zeta \Big) + \frac{p'(\rho)}{\rho} \abs{\nabla\p_i \phi}^2 + \dv (\cdots) \\
			& \qquad = \rho \abs{\dv\p_i\zeta}^2 + \e \nabla \p_i \phi \cdot I_1 + \e \dv \p_i \zeta \, I_2,
		\end{aligned}
	\end{equation}
	where we use $ I_1$ and $I_2$ to denote the sums of the remaining terms, which satisfy
	\begin{equation}\label{I-1-2}
		\begin{aligned}
		\abs{I_1} & \lesssim \e \abs{\p_t\rho} \abs{\nabla\p_i \phi} +  \sum_{j=1}^{2} \big(\e \abs{\nabla^j \uv} \abs{\nabla^{3-j}\phi} + \abs{\nabla^j \rho} \abs{\nabla^{3-j} \zeta} \big) \\
		 & \quad + \e \abs{\uv} \abs{\nabla\rho} \abs{\nabla^2 \phi} + \e \abs{\nabla^2 \widecheck{L}_0} + \sum_{j=0}^{1}  \abs{\nabla^j \uv} \abs{\nabla^{2-j} \zeta} \\
		 & \quad  + \e^{-1} \abs{\nabla\rho} \abs{\nabla\phi} + \abs{\nabla \widecheck{\Lv}}, \\
		 \abs{I_2} & \lesssim \abs{\nabla \uv} \abs{\nabla\phi} + \e^{-1} \abs{\nabla\rho} \abs{\nabla\zeta} + \abs{\nabla \widecheck{L}_0}.
	\end{aligned}
	\end{equation}
	Now we estimate the terms on the right-hand side of \cref{eq7a}.
	~\\
	1) It is noted that the divergence of the velocity perturbation satisfies that $\dv \zeta = \p_{3} \zeta_3^\od + \dv \zeta^\md,$ which excludes the tangential zero mode, $\zeta_\perp^\od$.
	Then using \cref{rel-0-2,rel-zeta-w}, one can get that if $\e^{\frac{1}{2}} M \leq 1,$ then
	\begin{equation}\label{L2-dv-zeta}
		\norm{\dv\p_i \zeta}_{L^2}^2 \lesssim \norm{\p_3^3 Z_3}_{L^2}^2 + \norm{\nabla^2 \wv^\md}_{L^2}^2 + \e \norm{\nabla^2 \phi^\md}_{L^2}^2 + \e M^2 (t+1)^{-\frac{5}{2}}.
	\end{equation}
	2) Then we estimate $\norm{I_i}_{L^2} $ for $i=1,2$ in \cref{I-1-2}. First, it follows from \cref{aL2-pert} and \cref{Lem-rel} that 
	\begin{equation}\label{L2-zeta}
		\begin{aligned}
			\norm{\nabla^j (\phi, \zeta)}_{L^2} & \lesssim M (t+1)^{ -\frac{2j+1}{4} }, \qquad j=0,1, \\
			\norm{\nabla^2 \phi}_{L^2} & \lesssim M (t+1)^{-\frac{3}{4}}, \\
			\norm{\nabla^2 \zeta}_{L^2} & \lesssim \norm{\p_3^3 \Zv}_{L^2} + \norm{\nabla^2 (\wv^\md, \phi^\md)}_{L^2} + \e^{\frac{1}{2}} M (t+1)^{-\frac{5}{4}} \\
			& \lesssim M (t+1)^{-\frac{3}{4}}.
		\end{aligned}
	\end{equation}
	Plugging these estimates into \cref{est-Lc} implies that
	\begin{equation}\label{L2-Lc}
		\begin{aligned}
		\norm{\nabla^k \widecheck{L}_0}_{L^2} & \lesssim \chi M (t+1)^{-\frac{2k+5}{4}}, \quad \text{for } k=1,2, \\
		\norm{\nabla \widecheck{\Lv}}_{L^2} & \lesssim M (t+1)^{-\frac{5}{4}}.
	\end{aligned}
	\end{equation}
	Then using \cref{bdd-diff,bdd-u,L2-Lc}, one can estimate the linear terms in \cref{I-1-2} and then obtain that 
	\begin{align*}
		\norm{I_1}_{L^2} & \lesssim \e \norm{\p_t \phi}_{L^\infty} \norm{\nabla^2 \phi}_{L^2} + \e \norm{\nabla \zeta}_{L^\infty} \norm{\nabla^2 \phi}_{L^2} + \e \norm{\nabla \phi}_{L^\infty}\norm{\nabla^2 \zeta}_{L^2} \\
		 & \quad + \e M \norm{\nabla \phi}_{L^\infty} \norm{\nabla^2 \phi}_{L^2} + M \norm{\nabla^2 \zeta}_{L^2}  + \norm{\nabla \zeta}_{L^\infty}\norm{\nabla \zeta}_{L^2} \\
		 & \quad + \norm{\nabla \phi}_{L^\infty} \norm{\nabla \phi}_{L^2} +  M (t+1)^{-\frac{5}{4}}, \\
		 \norm{I_2}_{L^2} & \lesssim \norm{\nabla \zeta}_{L^2} \norm{\nabla \phi}_{L^\infty} + M (t+1)^{-\frac{5}{4}}.
	\end{align*}
	Using the $L^\infty$-norms in \cref{aLinf-pert}, one has that if $\e^{\frac{1}{2}} M \leq 1,$ then
	\begin{align*}
		\norm{\zeta}_{L^\infty} & \lesssim \norm{\wv}_{L^\infty} \lesssim M (t+1)^{-\frac{1}{2}}, \\
		\norm{\nabla\zeta}_{L^\infty} & \lesssim \norm{\nabla\wv}_{L^\infty} + \norm{\nabla\rho}_{L^\infty} \norm{\wv}_{L^\infty} \lesssim M (t+1)^{-\frac{3}{4}}.
	\end{align*}
	In addition, it follows from \cref{equ-phizeta}$_1$ that
	\begin{align*}
		\e \norm{\p_t \phi}_{L^\infty} & \lesssim \e \norm{\uv}_{L^\infty} \norm{\nabla\phi}_{L^\infty} + \norm{\dv \zeta}_{L^\infty} + \e \norm{ \widecheck{L}_0}_{L^\infty} \\
		& \lesssim M (t+1)^{-\frac{3}{4}}.
	\end{align*}
	Collecting these $L^\infty$-bounds and using \cref{L2-zeta}, one can get that if $\e^{\frac{1}{2}} M \leq 1,$ then
	\begin{align*}
		\norm{I_1}_{L^2} & \lesssim  M^2 (t+1)^{-\frac{5}{4}} + M \norm{\nabla^2 \zeta}_{L^2}, \qquad
		\norm{I_2}_{L^2} \lesssim  M^2 (t+1)^{-\frac{5}{4}}.
	\end{align*}
	Hence, with the use of \cref{L2-zeta}$_3$, the last two terms on the right-hand side of \cref{eq7a} satisfy that
	\begin{equation}\label{I1-1}
		\begin{aligned}
			& \int_\Omega \e \nabla^2 \phi \cdot I_1 dx \lesssim \e \norm{\nabla^2 \phi}_{L^2} \big[ M^2 (t+1)^{-\frac{5}{4}} + M \norm{\nabla^2 \zeta}_{L^2} \big] \\
			& \qquad \leq \int_\Omega \frac{p'(\rho)}{100 \rho} \abs{\nabla^2 \phi}^2 dx + C \e^2 M^2 \big(\norm{\p_3^3 \Zv}_{L^2}^2 + \norm{\nabla^2 (\phi^\md, \wv^\md)}_{L^2}^2 \big) \\
			& \qquad\quad + C \e^2 M^4 (t+1)^{-\frac{5}{2}},
		\end{aligned}
	\end{equation}
	and 
	\begin{equation}\label{ineq-I2}
		\int_\Omega \e \dv \p_i\zeta I_2 dx \lesssim \norm{\dv \p_i\zeta}_{L^2}^2 + \e^2 M^4 (t+1)^{-\frac{5}{2}}.
	\end{equation}
	By integrating \cref{eq7a} over $\Omega$, and using \cref{I1-1,L2-dv-zeta,ineq-I2}, one can get \cref{est8}.
	
\end{proof}

\begin{Lem}
	If $\e^{\frac{1}{2}} M \leq 1, $ then the estimate \cref{est8} implies that
	\begin{equation}\label{est8*1}
		\begin{aligned}
			& \frac{d}{dt} \big[(t+1) \Ac_{2,\e\phi}\big] + (t+1) \norm{\nabla^2 \phi}_{L^2}^2 \\
			& \qquad \lesssim (t+1)^{-1} \Nc_\od^* + \Nc_\md + \e (t+1)^{-1} \Nc_\od + \e M^2 (t+1)^{-\frac{3}{2}},
		\end{aligned}
	\end{equation}
	and 
	\begin{equation}\label{est8**}
	\begin{aligned}
	& \frac{d}{dt} \big[(t+1)^2 \Ac_{2,\e\phi}\big] + (t+1)^2 \norm{\nabla^2 \phi}_{L^2}^2 \\
	& \qquad \lesssim \Nc_{\od}^* + (t+1)^2 \norm{\nabla^3 \zeta}_{L^2}^2 + \e \Nc_\od + \e M^2 (t+1)^{-\frac{1}{2}}. 
\end{aligned}
\end{equation}
Here $\Nc_\od^*, \Nc_\od$ and $\Nc_\md$ are the dissipative terms defined in \cref{ED-od*,ED-od,Ec-md}, respectively.
\end{Lem}

\begin{proof}
Applying \cref{L2-zeta} onto \cref{A2}, one has that
$$\abs{\Ac_{2,\e\phi}} \lesssim \e M^2 (t+1)^{-\frac{3}{2}}.$$
Then multiplying $(t+1)$ and $(t+1)^2$ on \cref{est8}, one can obtain \cref{est8*1} and 
\begin{align*}
	& \frac{d}{dt} \big[(t+1)^2 \Ac_{2,\e\phi}\big] + (t+1)^2 \norm{\nabla^2 \phi}_{L^2}^2 \\
	& \qquad \lesssim \Nc_{\od}^* + (t+1)^2 \norm{\nabla^2 \wv^\md}_{L^2}^2 + \e \Nc_\od  + \e M^2 (t+1)^{-\frac{1}{2}},
\end{align*}
respectively.
It follows from \cref{rel-zw3} that
\begin{align*}
	(t+1)^2 \norm{\nabla^2 \wv^\md}_{L^2}^2 & \lesssim (t+1)^2 \norm{\nabla^3 \wv^\md}_{L^2}^2 \lesssim (t+1)^2 \norm{\nabla^3 \zeta}_{L^2}^2 + \e M^2 (t+1)^{-\frac{1}{2}}.
\end{align*}
Then one can get \cref{est8**}.
\end{proof}

\vspace{.2cm}

\subsection{Proof of Step 5}
It remains to establish the $H^1$-estimate for $\nabla^2 (\phi, \zeta),$ which completes the proof of Step 5 in \cref{Sec-steps}.

\begin{Lem}\label{Lem-est7}
	If $ \big(\chi^{\frac{1}{2}} + \nu^{\frac{1}{6}} + \e^{\frac{1}{4}} \big) M \leq 1,$ then it holds that
	\begin{equation}\label{est7}
		\begin{aligned}
			& \frac{d}{dt} \big[(t+1)^2 \Ac_2 \big] +  (t+1)^2 \norm{\nabla^3 \zeta}_{L^2}^2  \\
			& \qquad \lesssim \Nc_\od(t) + \Nc_\md(t) + \big(\Lambda^{-\frac{1}{2}} + \chi^{\frac{1}{2}} + \nu^{\frac{1}{8}} + \e^{\frac{1}{2}} \big) (t+1)^2 \norm{\nabla^2 \phi}_{L^2}^2 \\
			& \qquad \quad  + (t+1) \norm{\nabla^2 \phi}_{L^2}^2 + (\chi+ \e) M^2 (t+1)^{-\frac{1}{2}},
		\end{aligned}
	\end{equation}
	where $\Ac_2$ is a functional satisfying 
	\begin{equation}\label{A1}
		\Ac_2 \sim \norm{\nabla^2(\phi, \zeta)}_{L^2}^2,
	\end{equation}
	and $\Nc_\od$ and $\Nc_\md$ are given by \cref{ED-od} and \cref{Ec-md}, respectively.
\end{Lem}

\begin{proof}	
	Let $ i \in \{1,2,3\} $ be fixed.
	~\\
	1) Multiplying $ \p_i $\cref{equ-phizeta}\textsubscript{2} by $-\lap \p_i \zeta $, and using the facts derived from \cref{fact} that
	\begin{align*}
		-\frac{p'(\rho)}{\e \rho} \nabla\p_i \phi \cdot \big(\lap \p_i \zeta - \nabla \dv \p_i\zeta \big) & = -\frac{p'(\rho)}{\e \rho} \dv \big( \nabla\p_i \phi \times \text{curl } \p_i \zeta \big), \\
		-\frac{\mu+\lambda}{\rho} \nabla \dv \p_i \zeta \cdot \big(\lap \p_i \zeta - \nabla \dv \p_i\zeta \big) & = -\frac{\mu+\lambda}{\rho} \dv \big( \nabla \dv \p_i \zeta \times \text{curl } \p_i \zeta \big),
	\end{align*}
	one can get that
	\begin{equation}\label{eq7-1}
		\begin{aligned}
		& \p_t \Big( \frac{1}{2} \abs{\nabla\p_i  \zeta}^2 \Big) - \frac{p'(\rho)}{\e \rho} \nabla \p_i \phi \cdot \nabla \p_i \dv \zeta + \frac{\mu}{\rho} \abs{ \lap\p_i\zeta}^2 + \frac{\mu+\lambda}{\rho} \abs{\nabla\dv \p_i \zeta}^2 \\
		& \qquad = \dv (\cdots) + I_3  + I_4,
		\end{aligned}
	\end{equation}
	where $$I_3 = \lap \p_i \zeta \cdot \p_i (\uv \dnab \zeta),$$ and $I_4$ denotes the remaining terms, satisfying that
	\begin{equation}\label{I4}
		\begin{aligned}
			\abs{I_4} & \lesssim \abs{\lap \p_i\zeta} \big( \e^{-1} \norm{\p_i \rho}_{L^\infty} \abs{\nabla\phi} + \abs{\p_i \widecheck{\Lv}} \big) \\
			& \quad + \e^{-1} \norm{\nabla\rho}_{L^\infty} \abs{\text{curl}\, \p_i\zeta} \big( \abs{\nabla\p_i\phi} + \abs{\nabla\dv\p_i\zeta} \big).
		\end{aligned}
	\end{equation}
	Due to the largeness and slow decay rate of the velocity field, we need to estimate $I_3$ through the equality,
	\begin{equation*}
		\begin{aligned}
			I_3 & = \dv(\cdots) - \sum_{j=1}^{3} \p_{ij} \zeta \cdot \big( \p_{ij} \uv \cdot \nabla \zeta + 2 \p_i \uv \dnab \p_j \zeta \big) + \frac{1}{2} \dv \uv \abs{\nabla\p_i \zeta}^2.
		\end{aligned}
	\end{equation*}
	This, together with \cref{bdd-diff}, yields that
	\begin{equation}\label{est-I3}
		\begin{aligned}
		\abs{ \int_\Omega I_3 dx } & \lesssim \Lambda^{-\frac{1}{2}} \sum_{j=1}^2 (t+1)^{-2+j} \norm{\nabla^j \zeta}_{L^2}^2 \\
		& \quad + \norm{\nabla\zeta}_{L^\infty} \norm{\nabla^2 \zeta}_{L^2}^2 + \chi M^2 (t+1)^{-\frac{5}{2}}.
	\end{aligned}
	\end{equation}	
	To estimate $I_4$ in \cref{I4},  it is noted that 
	\begin{equation}\label{Linf-rho-1}
		\e^{-1}\norm{\nabla\rho}_{L^\infty} \lesssim \chi (t+1)^{-1} + M^{\frac{3}{4}} \nu^{\frac{1}{4}} (t+1)^{-\frac{1}{2}}.
	\end{equation}
	This, together with \cref{est-Lc}$_2,$ yields that if $M^{\frac{3}{4}} \nu^{\frac{1}{8}} \leq 1,$ then
	\begin{equation}\label{est-I4}
		\begin{aligned}
		\int_\Omega \abs{I_4} dx
		& \leq \frac{\mu}{100\rhob} \norm{\nabla^3 \zeta}_{L^2}^2 + C \big[  \big(\chi+ \nu^{\frac{1}{8}}\big) \norm{\nabla^2 (\phi, \zeta)}_{L^2}^2 \\
		& \qquad + \sum_{j=0}^{1} (t+1)^{-2+j} \norm{\nabla^j (\phi, \zeta)}_{L^2}^2 + \chi^2 (t+1)^{-\frac{7}{2}}\big].
	\end{aligned}
	\end{equation}
	~\\
	2) On the other hand, multiplying $ \nabla \p_i  $\cref{equ-phizeta}\textsubscript{1} by $ \frac{p'(\rho)}{\rho^2} \nabla \p_i \phi $ yields that
	\begin{equation}\label{eq7-2}
		\p_t \Big(\frac{p'(\rho) }{2\rho^2} \abs{\nabla \p_i \phi}^2 \Big)  + \frac{p'(\rho)}{\e \rho} \nabla \p_i \phi \cdot \nabla \p_i \dv \zeta = \dv (\cdots) + I_5,
	\end{equation}
	where $I_5$ denotes the sum of several lower-order terms, satisfying that
	\begin{equation*}
		\begin{aligned}
			\abs{I_5} & \lesssim \big( \abs{\nabla\uv} + \abs{\uv} \abs{\nabla\rho} \big) \abs{\nabla^2 \phi}^2 + \abs{\nabla^2 \uv} \abs{\nabla \phi} \abs{\nabla^2 \phi} \\
			& \quad + \abs{\nabla^2 \phi} \big( \e^{-1} \abs{\nabla\rho} \abs{\nabla^2 \zeta} + \e^{-1} \abs{\nabla^2 \rho} \abs{\nabla \zeta}  + \abs{\nabla^2 \widecheck{L}_0} \big).
		\end{aligned}
	\end{equation*}
	Using \cref{bdd-diff,bdd-u,L2-Lc,Linf-rho-1}, one has that if $\e^{\frac{1}{2}} M \leq 1, $ then
	\begin{equation}
		\begin{aligned}
			\int_\Omega \abs{I_5} dx & \lesssim \big(\Lambda^{-\frac{1}{2}} + \e^{\frac{1}{2}} + \chi + \nu^{\frac{1}{8}}\big) \norm{\nabla^2 (\phi, \zeta)}_{L^2}^2 + \norm{\nabla\zeta}_{L^\infty} \norm{\nabla^2 \phi}_{L^2}^2  \\
			& \quad + \sum_{j=0}^{1} (t+1)^{-2+j} \norm{\nabla^j (\phi, \zeta)}_{L^2}^2 + \chi M^2 (t+1)^{-\frac{9}{2}}.
		\end{aligned}
	\end{equation}
	It follows from the $L^\infty$-bounds in \cref{Sec-bound} and \cref{Zp-inf} that if $M^{\frac{3}{4}} \nu^{\frac{1}{8}} \leq 1,$ then
	\begin{align}
		\norm{\nabla\zeta}_{L^\infty} & \lesssim \norm{\p_3^2 \Zv_\perp}_{L^\infty} + \norm{\p_3^2 Z_3}_{L^\infty} + \norm{\nabla \wv^\md}_{L^\infty} \notag \\
		& \quad + \e \big[ \chi (t+1)^{-1} + \norm{\nabla\phi}_{L^\infty} \big] \norm{\wv}_{L^\infty} \notag \\
		& \lesssim (t+1)^{-\frac{5}{4}} + \chi^{\frac{1}{2}} + \nu^{\frac{1}{8}} + \e^{\frac{1}{2}}. \label{Linf-zeta-1}
	\end{align}
	This, together with \cref{L2-zeta}, \cref{L2-Lc}$_1$ and \cref{Linf-rho-1}, implies that if $M^{\frac{3}{4}} \nu^{\frac{1}{8}} \leq 1,$ then
	\begin{equation}\label{est-I5}
		\begin{aligned}
			\int_\Omega \abs{I_5} dx & \lesssim \big(\Lambda^{-\frac{1}{2}} + \chi^{\frac{1}{2}} + \e^{\frac{1}{2}} + \nu^{\frac{1}{8}}\big) \norm{\nabla^2 (\phi,\zeta)}_{L^2}^2 + (t+1)^{-\frac{5}{4}} \norm{\nabla^2 \phi}_{L^2}^2 \\
			& \quad + \sum_{j=0}^{1} (t+1)^{-2+j} \norm{\nabla^j (\phi, \zeta)}_{L^2}^2 + \chi M^2 (t+1)^{-\frac{9}{2}}.
		\end{aligned}
	\end{equation}
	~\\
	3) Note that $\norm{\lap \p_i\zeta}_{L^2} = \norm{\nabla^2 \p_i \zeta}_{L^2}. $ 
	Then by adding \cref{eq7-1,eq7-2} together and using \cref{est-I3,est-I4,est-I5}, one can get that
	\begin{equation}\label{ineq-1}
		\begin{aligned}
			& \frac{d}{dt} \Big[ (t+1)^2 \Big( \norm{\nabla^2 \zeta}_{L^2}^2 + \int_\Omega \frac{p'(\rho)}{\rho^2} \abs{\nabla^2 \phi}^2 dx \Big) \Big] + \frac{\mu}{4\rhob} (t+1)^2 \norm{\nabla^3 \zeta}_{L^2}^2 \\
			& \qquad \lesssim \big(\Lambda^{-\frac{1}{2}} + \chi^{\frac{1}{2}} + \nu^{\frac{1}{8}} + \e^{\frac{1}{2}} \big) (t+1)^2 \norm{\nabla^2 (\phi, \zeta)}_{L^2}^2 + (t+1) \norm{\nabla^2 (\phi, \zeta)}_{L^2}^2  \\
			& \qquad \quad  + \sum_{j=0}^{1} (t+1)^{j} \norm{\nabla^j (\phi, \zeta)}_{L^2}^2 + \chi M^2 (t+1)^{-\frac{1}{2}}.
		\end{aligned}
	\end{equation}
	
	Next we estimate the terms on the right-hand side of \cref{ineq-1}.
	\begin{itemize}
		\item It follows from \cref{rel-0-2,rel-zeta-w} that
		\begin{equation}\label{rel-zeta2}
			\norm{\nabla^2 \zeta}_{L^2}^2 \lesssim \norm{\p_3^3 \Zv}_{L^2}^2 + \norm{\nabla^2 \wv^\md}_{L^2}^2 + \e M^2 (t+1)^{-\frac{5}{2}},
		\end{equation}
		which yields that
		\begin{align}\label{est-p2z1}
			(t+1) \norm{\nabla^2 \zeta}_{L^2}^2 \lesssim \Nc_\od + \Nc_\md + \e M^2 (t+1)^{-\frac{3}{2}}.
		\end{align}

		\item Using \cref{rel-zeta2,rel-zw3}, it also holds that
		\begin{align*}
			\norm{\nabla^2 \zeta}_{L^2}^2 & \lesssim \norm{\p_3^3 \Zv}_{L^2}^2 + \norm{\nabla^3 \wv^\md}_{L^2}^2 + \e M^2 (t+1)^{-\frac{5}{2}} \\
			& \lesssim \norm{\p_3^3 \Zv}_{L^2}^2 + \norm{\nabla^3 \zeta}_{L^2}^2 + \e M^2 (t+1)^{-\frac{5}{2}}.
		\end{align*}
		Then the term on the third line of \cref{ineq-1} satisfies that 
		\begin{equation}\label{est-p2z2}
			(t+1)^2 \norm{\nabla^2 \zeta}_{L^2}^2 \lesssim \Nc_\od + (t+1)^2 \norm{\nabla^3 \zeta}_{L^2}^2 + \e M^2 (t+1)^{-\frac{1}{2}}.
		\end{equation}
		
		\item It follows from \cref{rel-0-2} that
		\begin{align}
			& \sum_{j=0}^{1} (t+1)^{j} \norm{\nabla^j (\phi, \zeta)}_{L^2}^2 \notag \\
			& \quad \lesssim \sum_{j=1}^{2} (t+1)^{-1+j} \norm{\p_3^j (\Phi, \Zv)}_{L^2}^2 + (t+1) \norm{\nabla \big( \phi^\md, \wv^\md\big)}_{L^2}^2 + \e M^2 (t+1)^{-\frac{1}{2}} \notag \\
			& \quad \lesssim \Nc_\od + \Nc_\md +  \e M^2 (t+1)^{-\frac{1}{2}}. \label{est-1pz}
		\end{align}

	\end{itemize}
		
	Then applying \cref{est-1pz,est-p2z1,est-p2z2} in \cref{ineq-1}, one can finish the proof.
	
\end{proof}


\begin{Lem}[Second-order estimate]\label{Lem-2est}
	If $ \big(\chi^{\frac{1}{2}} + \nu^{\frac{1}{6}} + \e^{\frac{1}{4}} \big) M \leq 1,$ then it holds that
	\begin{equation}\label{ineq1}
		\frac{d}{dt} \Ec^{(1)}_{ho}(t) + \Nc^{(1)}_{ho}(t) \lesssim \Nc_\od(t) + \Nc_\md(t) + (\chi+\e ) M^2 (t+1)^{-\frac{1}{2}},
	\end{equation}
where $\Ec^{(1)}_{ho}$ and $ \Nc^{(1)}_{ho} $ are two functionals satisfying that
\begin{equation}\label{ED-1}
	\begin{aligned}
		\pm \Ec^{(1)}_{ho}(t) & \lesssim \pm (t+1)^2 \norm{\nabla^2 (\phi, \zeta)}_{L^2}^2 + \e M^2 (t+1)^{\frac{1}{2}},  \\
		\Nc^{(1)}_{ho}(t) & \sim (t+1)^2 \norm{\big(\nabla^2 \phi, \nabla^3 \zeta\big)}_{L^2}^2,
	\end{aligned}
\end{equation}
and $\Nc_\od$ and $\Nc_\md$ are given by \cref{ED-od} and \cref{Ec-md}, respectively.
\end{Lem}

\begin{proof}
First, we can combine \cref{est8**,est7} to get that there exists a suitable large constant $ C_1 \sim 1, $ such that
\begin{equation}\label{ineq-4}
	\begin{aligned}
		& \frac{d}{dt} \big[ (t+1)^2 \big(2 C_1 \Ac_2 + \Ac_{2,\e\phi}\big) \big] + (t+1)^2 \big(C_1 \norm{\nabla^3 \zeta}_{L^2}^2 + \norm{\nabla^2 \phi}_{L^2}^2 \big) \\
		& \qquad \lesssim \Nc_\od + \Nc_\md + (t+1) \norm{\nabla^2 \phi}_{L^2}^2 + (\chi+\e) M^2 (t+1)^{-\frac{1}{2}}.
	\end{aligned}
\end{equation}
Here we have let $\Lambda \geq C_2 $ for a suitably large constant $C_2=C_2(\rhob, \mu, |\uvb|, M_0).$

Secondly, the combination of \cref{ineq-4,est8*1} yields that one can choose a suitable large constant $ C_3 = C_3(\rhob, \mu, |\uvb|)>0 $ such that
\begin{equation}\label{ineq5}
	\begin{aligned}
	& \frac{d}{dt} \big[ (t+1)^2 \big(2 C_1 \Ac_2 + \Ac_{2,\e\phi} \big) + 2 C_3 (t+1) \Ac_{2,\e\phi} \big] \\
	& \qquad \qquad + (t+1)^2 \big(C_1 \norm{\nabla^3 \zeta}_{L^2}^2 + \norm{\nabla^2 \phi}_{L^2}^2 \big) + C_3 (t+1) \norm{\nabla^2 \phi}_{L^2}^2 \\
	& \qquad \lesssim \Nc_\od +  \Nc_\md + (\chi+\e) M^2 (t+1)^{-\frac{1}{2}}.
	\end{aligned}
\end{equation}
According to \cref{ineq5}, we denote
\begin{equation}\label{EN-ho-1}
	\begin{aligned}
		\Ec^{(1)}_{ho}(t) & := (t+1)^2 \big(2C_1 \Ac_2 + \Ac_{2,\e\phi} \big) + 2 C_3  (t+1) \Ac_{2,\e\phi}, \\
		\Nc^{(1)}_{ho}(t) & := (t+1)^2 \big(C_1 \norm{\nabla^3 \zeta}_{L^2}^2 + \norm{\nabla^2 \phi}_{L^2}^2 \big) + C_3 (t+1) \norm{\nabla^2 \phi}_{L^2}^2.
	\end{aligned}
\end{equation}
Then it remains to verify \cref{ED-1}$_1$. In fact, it follows from \cref{A1,A2} that
$$
\pm \Ec^{(1)}_{ho} \lesssim \pm (t+1)^2 \norm{\nabla^2 (\phi, \zeta)}_{L^2}^2 + \e (t+1)^2 \norm{\nabla\zeta}_{L^2}^2.
$$
Using \cref{rel-0-2,rel-zeta-w}, one has that
\begin{align*}\label{ineq6}
	\norm{\nabla\zeta}_{L^2} & \lesssim \norm{\p_3 \zeta^\od}_{L^2} + \norm{\nabla\zeta^\md}_{L^2} \\
	& \lesssim \norm{\p_3^2 \Zv}_{L^2} + \norm{\nabla^2 \big(\phi^\md, \wv^\md\big)}_{L^2}  + M (t+1)^{-\frac{3}{4}} \lesssim M (t+1)^{-\frac{3}{4}},
\end{align*}
where we have used the fact that $\norm{\nabla \big(\phi^\md, \wv^\md\big)}_{L^2}  \lesssim \norm{\nabla^2 \big(\phi^\md, \wv^\md\big)}_{L^2}.$ 
The proof is completed.

\end{proof}



Similarly, we can prove the following third-order estimates. It is noted that the key \cref{Lem-Zp-inf} is also essential in the proof of the following \cref{Lem-est9}.

\begin{Lem}\label{Lem-est10}
	If $\e^{\frac{1}{2}} M \leq 1,$ then it holds that
	\begin{equation}\label{est10}
		\begin{aligned}
			\frac{d}{dt} \Ac_{3,\e \phi} + \norm{\nabla^3 \phi}_{L^2}^2 & \lesssim \norm{\nabla^3 \zeta}_{L^2}^2 + \e M^2 (t+1)^{-\frac{5}{2}},
		\end{aligned}
	\end{equation}
	where $\Ac_{3,\e \phi}$ is a functional satisfying that
	\begin{equation}\label{AB4}
		\pm \Ac_{3,\e \phi} \lesssim  \pm \e^2 \norm{\nabla^3 \phi}_{L^2}^2 + \e \norm{\nabla^3 \phi}_{L^2} \norm{\nabla^2 \zeta}_{L^2}.
	\end{equation}
\end{Lem}


\begin{Lem}\label{Lem-est9}
	If $ \big(\chi^{\frac{1}{2}} + \nu^{\frac{1}{6}} + \e^{\frac{1}{4}} \big) M \leq 1,$ then it holds that
	\begin{equation}\label{est9}
		\begin{aligned}
			\frac{d}{dt} \Ac_3 +  \norm{\nabla^4 \zeta}_{L^2}^2 & \lesssim (t+1)^{-2} \big(\Nc_\od + \Nc_\md + \Nc_{ho}^{(1)}\big) \\
			& \quad + \norm{\nabla^3 \phi}_{L^2}^2 + (\chi+\e) M^2 (t+1)^{-\frac{5}{2}},
		\end{aligned}
	\end{equation}
	where $\Ac_3$ is a functional satisfying that
	\begin{align*}
		\Ac_3 \sim \norm{\nabla^3 \big(\phi,  \zeta\big)}_{L^2}^2,
	\end{align*}
	and $\Nc_\od$, $\Nc_\md$ and $ \Nc_{ho}^{(1)} $ are given by \cref{ED-od,Ec-md,ED-1}, respectively.
\end{Lem}

We show only the sketch of the proof of \cref{Lem-est9}.

\begin{proof}[Proof of \cref{Lem-est9}]
In the lower-order estimates, we have shown that if $ \big(\chi^{\frac{1}{2}} + \nu^{\frac{1}{6}} + \e^{\frac{1}{4}} \big) M \leq 1,$ 
\begin{equation}\label{bdd-gradient}
	\begin{aligned}
	& \norm{\nabla\phi}_{L^\infty} \lesssim \nu^{\frac{1}{8}}, \qquad \norm{\dv \zeta}_{L^\infty} \lesssim \nu^{\frac{1}{8}} + \e^{\frac{1}{2}}, \\
	& \norm{\nabla\zeta}_{L^\infty} \lesssim (t+1)^{-\frac{5}{4}} + \chi^{\frac{1}{2}} + \nu^{\frac{1}{8}} + \e^{\frac{1}{2}}.
\end{aligned}
\end{equation}
Then similar to the proof of \cref{Lem-est7}, one can get that
\begin{equation*} 
	\begin{aligned}
	& \frac{d }{dt} \Big( \norm{\nabla^3 \zeta}_{L^2}^2 + \int_\Omega \frac{p'(\rho)}{2 \rho^2} \abs{\nabla^3 \phi}^2 dx \Big) + c_0 \norm{\nabla^4 \zeta}_{L^2}^2 \\
	& \qquad \lesssim \norm{\nabla^3 (\phi, \zeta)}_{L^2}^2 + (\chi +\e) M^2 (t+1)^{-\frac{5}{2}} + \norm{\nabla^2 \widecheck{\Lv}}_{L^2}^2 + I_6 + I_7,
\end{aligned}
\end{equation*}
where 
\begin{equation*}
	\begin{aligned}
	I_6 & := \norm{\nabla^2 \phi}_{L^4}^2 \norm{\nabla^2 \zeta}_{L^4}^2 + \norm{\nabla^2 \phi}_{L^4}^2 \norm{\nabla^3 \zeta}_{L^2}, \\
	I_7 & = \sum_{i,j=1}^{3} \abs{\int_\Omega \lap \p_{ij} \zeta \cdot \p_{ij} \big(\uv \dnab \zeta\big) dx}.
\end{aligned}
\end{equation*}
Using \cref{est-Lc}$_2$, \cref{est-p2z2}, \cref{est-1pz}, one has that
\begin{equation}\label{Lv-2}
	\norm{\nabla^2 \widecheck{\Lv}}_{L^2}^2 \lesssim (t+1)^{-3} (\Nc_\od + \Nc_\md) + \norm{\nabla^3 \zeta}_{L^2}^2 + \e M^2 (t+1)^{-\frac{7}{2}}.
\end{equation}

On the other hand, it follows from \cref{Lem-GN} that for any $h\in H^3(\Omega), $
\begin{equation}\label{GN-1}
	\norm{\nabla^2 h}_{L^4} \lesssim \norm{\nabla^3 h}_{L^2}^{\frac{1}{2}} \norm{\nabla h}_{L^\infty}^{\frac{1}{2}}.
\end{equation}
Then with the use of \cref{bdd-gradient}, one can get that
\begin{align}
	I_6 & \lesssim \norm{\nabla\phi}_{L^\infty} \norm{\nabla\zeta}_{L^\infty} \norm{\nabla^3 \phi}_{L^2} \norm{\nabla^3 \zeta}_{L^2} + \norm{\nabla\phi}_{L^\infty} \norm{\nabla^3 \phi}_{L^2} \norm{\nabla^3 \zeta}_{L^2} \notag \\
	& \lesssim \norm{\nabla^3 (\phi, \zeta)}_{L^2}^2. \label{est-I6}
\end{align}
To estimate $I_7,$ using integration by parts and the estimates, \cref{est-p2z2}, \cref{est-1pz} and \cref{GN-1}, one has
\begin{align*}
	& \abs{\int_\Omega \lap \p_{ij} \zeta \cdot \p_{ij} \big(\uv \dnab \zeta\big) dx} \\
	& \qquad \lesssim \sum_{k=1}^{3} \int_\Omega \abs{\nabla^3 \zeta} \abs{\nabla^k \uv} \abs{\nabla^{4-k}\zeta} dx \\
	& \qquad \lesssim \sum_{k=1}^{3} (t+\Lambda)^{-\frac{2k-1}{2}} \norm{\nabla^{4-k} \zeta}_{L^2}^2 + \norm{\nabla \zeta}_{L^\infty} \norm{\nabla^3 \zeta}_{L^2}^2 \\
	& \qquad \lesssim \norm{\nabla^3 \zeta}_{L^2}^2 + (t+1)^{-\frac{7}{2}} (\Nc_\od + \Nc_\md) + \e M^2 (t+1)^{-4}.
\end{align*}
This, together with \cref{Lv-2,est-I6}, can complete the proof.
\end{proof}

\vspace{.2cm}
One can use \cref{est9,est10} directly to obtain that there exist two functionals $ \Ec_{ho}^{(2)} $ and $ \Nc_{ho}^{(2)} $ satisfying that 
\begin{equation}\label{ineq2}
	\frac{d}{dt} \Ec_{ho}^{(2)}(t) + \Nc^{(2)}_{ho} \lesssim \Nc_\od + \Nc_\md + \Nc^{(1)}_{ho} + (\chi+\e) M^2 (t+1)^{-\frac{1}{2}},
\end{equation}
and
\begin{equation}\label{ED-2}
	\begin{aligned}
		\pm \Ec_{ho}^{(2)}(t) & \lesssim \pm (t+1)^2 \norm{\nabla^3 (\phi, \zeta)}_{L^2}^2 + \e M^2 (t+1)^{\frac{1}{2}},  \\
		\Nc^{(2)}_{ho}(t) & \sim (t+1)^2 \norm{\big(\nabla^3 \phi, \nabla^4 \zeta\big)}_{L^2}^2.
	\end{aligned}
\end{equation}
This, together with \cref{Lem-2est}, can complete the proof of \cref{ineq-ho}.

%
%

\vspace{.2cm}

\section{Incompressible limit}\label{Sec-lim}





The low Mach number limit $\e\to 0$ is used to simply the fluid dynamics of highly subsonic flows, and the mathematical theory of this approximation has been widely studied; see \cite{Ebin1977,KM1981,KM1982,Metivier2001} and \cite{DG1999,Danchin2002,Alazard06} for instance.

In this section, we prove the low Mach number limit around the vortex sheets for the Navier-Stokes equations.
It is noted that given any fixed $t_0>0,$ the vortex layer \cref{profile-t0},  
\begin{equation}\label{vs-incom}
    \uv^\vs(x_3,t) = \Theta\Big( \frac{x_3}{\sqrt{t+t_0}} \Big) \uvb,
\end{equation}
is independent of the Mach number and is also a smooth solution to the incompressible Navier-Stokes equations, 
\begin{equation}\label{InNS}
\begin{cases}
    \rhob \, \p_t \uv + \rhob \, \uv \dnab \uv + \nabla \mathcal{P} = \mu \lap \uv, \\
    \dv \uv = 0.
\end{cases}
\end{equation}

As a byproduct of \cref{Thm}, we can utilize the classical analysis in \cite{Metivier2001,Alazard06} to establish the incompressible limit of the solutions $\{(\rho^\e, \uv^\e)\}_{\e>0}$ to the Cauchy problem, \cref{NS,ic}.

\begin{Thm}[Incompressible limit]\label{thm: lim}
    Under the assumptions of \cref{Thm}, assume further that the perturbation $(b_0, \vv_0) \in H^4_{3/2}(\Omega)$. Denote that
    \begin{align}
        M_0 & := \norm{(b_0, \vv_0)^\od}_{L^2_{3/2}} + \norm{(b_0, \vv_0)}_{H^4},\\
        \chi & := \norm{(b_0, \vv_{03})^\od}_{L^2_{3/2}} + \norm{(b_0, \vv_0)^\md}_{H^1}.
    \end{align}
    If $\chi\leq \chi_0, 0< \e \leq \e_0$ and $(\chi+\e) M_0^{k_0} \leq 1$, where $\chi_0, \e_0$ and $k_0$ are the constants in \cref{Thm}, then the sequence $\{(\rho^\e, \uv^\e)\}_{0<\e\leq \e_0}$ converges to a limit $(\rhob,\uv^0) \in C(0,\infty; H^4(\Omega))$ weakly in $ L^\infty(0,\infty; H^4(\Omega)), $ and strongly in $L^2(0,\infty; H^{4-}_{loc}(\Omega))$. 
    In addition, 
	the limit $\uv^0$ is a classical solution to the incompressible Navier-Stokes equations \cref{InNS} with the initial data
    \begin{equation}\label{ic-incom}
        \uv(x,t=0) = \uv^\vs(x_3,t=0) + \Pi\zetab_0(x), \quad x\in \Omega,
    \end{equation}
    where $\Pi \zetab_0 = (\mathrm{Id}- \nabla \lap^{-1} \dv) \zetab_0 \in H^4(\Omega)$. 
\end{Thm}

As a corollary of \cref{thm: lim}, we can achieve the nonlinear asymptotic stability of the vortex layer for the incompressible Navier-Stokes equations.

\begin{Thm}[Stability in incompressible flows]
    Given any fixed $\rhob>0,$ $\uvb = (\ub_1, \ub_2, 0) \in \R^3$ and $t_0>0,$ let $\uv^\vs=\uv^\vs(x_3,t)$ be the vortex layer given by \cref{vs-incom}.
    Suppose that $\vv_0=(v_{01}, v_{02}, v_{03})(x) \in H^4(\Omega) $ is solenoidal.
    Denote $$ M_0 := \norm{\vv_0}_{H^4_{3/2}} \andd \chi := \norm{\vv_0^\md}_{H^1}. $$
    Then there exist 
    \begin{itemize}
        \item a constant $ \chi_0>0 $, depending on $\mu, \rhob, t_0$ and $\max\{\abs{\uvb}, 1\} $,

        \item and an integer $k_0>0,$ depending only on the space dimension,
    \end{itemize}
    such that given any $M_0>0,$ if $\chi \leq \chi_0$ and $\chi M_0^{k_0} \leq 1,$ then the Cauchy problem for the incompressible Navier-Stokes equations \cref{InNS} with the initial data 
    \begin{equation}
        \uv(x,t=0) = \uv^\vs(x_3,t=0) + \vv_0(x), \quad x\in\Omega,
    \end{equation}
    admits a classical bounded solution $\uv$ globally in time, satisfying that
    \begin{equation}
        \begin{aligned}
            & \sup_{t\geq 0} \norm{\uv -\uv^\vs}_{H^4}^2 + \int_0^\infty \norm{ \nabla (\uv -\uv^\vs) }_{H^{4}}^2 dt \leq C, \\
			& \norm{\uv - \uv^\vs}_{L^\infty} \leq C (t+1)^{-\frac{1}{2}},
        \end{aligned}
        \end{equation}
    where $ C>0 $ is a constant, independent of $ t. $
\end{Thm}

In the rest of this section, we give a sketchy proof of Theorem \ref{thm: lim}.
\vspace{.2cm}
~\\
\textit{Proof of \cref{thm: lim}.}
In the proof, we put back the upper index $\e$ and use the notations,
\begin{itemize}
	\item $(\rho^\e, \uv^\e)$ is the solution to the Cauchy problem \cref{NS}, \cref{ic}, 
	\item $(\rhot^\e, \uvt^\e)$ is the ansatz defined in \cref{ansatz}, 
	\item $(\phi^\e, \zeta^\e) $ is the perturbation defined by $ (\phi^\e, \zeta^\e) = (\e^{-1}(\rho^\e - \rhot^\e), \uv^\e - \uvt^\e). $
\end{itemize}
Recall that both the background flow, $(\rhob, \uv^\vs)$ in \cref{profile-t0}, and the auxiliary one, $(\rhob, \uvt^\vs)$ in \cref{profile-t}, are independent of the Mach number $\e.$

Set
\begin{equation}
	b^\e = \e^{-1} (\rho^\e - \rhob), \qquad \vv^\e := \uv^\e - \uv^\vs,
\end{equation}
and
\begin{equation}
    q^\e := \e^{-1} (p(\rho^\e) - p(\rhob)), \qquad \eta^\e := \frac{1}{p'(\rho^\e) \rho^\e}=\frac{1}{\gamma p(\rho^\e)}.
\end{equation}
It follows from \cref{NS,profile-t} that
\begin{equation}\label{equ-Ue}
    \begin{cases}
        \eta^\e(\p_t q^\e+\uv^\e\cdot\nabla q^\e)
        +\frac{1}{\e}\dv \vv^\e=0, \\
        \rho^\e(\p_t\vv^\e+\uv^\e\cdot\nabla \vv^\e ) + \frac{1}{\e}\nabla q^\e \\
		\qquad = \mu\Delta \vv^\e+(\mu+\lambda)\nabla\dv \vv^\e - \e b^\e \p_t \uv^\vs - \rho^\e v^\e_3\p_3\uvs.
    \end{cases}
\end{equation}
Using \cref{Thm-L2} (the fourth-order estimates can be obtained similarly), one can get that
\begin{equation}\label{uni-bdd}
	\sup_{t\geq 0} \norm{(q^\e, \vv^\e)}_{H^4(\Omega)}^2  + \int_0^\infty \Big(\norm{\nabla q^\e}_{H^3(\Omega)}^2 + \norm{ \nabla \vv^\e }_{H^4(\Omega)}^2 \Big) dt \leq C,
\end{equation}
where the constant $C>0$ is independent of $\e.$ Thus, we can extract a subsequence of $\{(q^\e, \vv^\e)\}_{0<\e\leq \e_0}$ such that 
$$
(q^\e, \vv^\e) \to (q^0,\vv^0) \quad \text{weakly} * \text{ in } L^\infty(0,\infty; H^4(\Omega)).
$$
Then we use the analysis in \cite{Metivier2001,Alazard06} to prove the strong convergence.

\begin{Lem}\label{lem: lim of pressure}
	For all $T>0$, it holds that
	\begin{equation}\label{lim of pressure}
		\begin{aligned}
			q^\e & \to q^0 = 0 \quad \text{strongly in } \ L^2((0,T); H_{loc}^{4-}), \\
			\dv \vv^\e & \to \dv \vv^0 = 0 \quad \text{strongly in } \ L^2((0,T); H_{loc}^{4-}).
		\end{aligned}
	\end{equation}
\end{Lem}
\begin{proof}
	It follows from \cref{equ-Ue} that
	\begin{equation}\label{equ-wave}
		\e^2 \p_t (\eta^\e \p_t q^\e) 
        -\dv\big(\frac{1}{\rho^\e}\nabla q^\e\big)  = \e h^\e,
	\end{equation}
	where
	\begin{align*}
		h^\e & = \dv \big[ (\uv^\vs + \vv^\e) \dnab \vv^\e + \p_3 \uv^\vs v_3^\e + \frac{\e b^\e}{\rho^\e} \p_t \uv^\vs \\
		& \qquad\quad - \frac{1}{\rho^\e}\big( \mu \lap \vv^\e + (\mu+\lambda) \nabla \dv \vv^\e \big) \big] - \e \p_t \big[\eta^\e(\uv^\vs + \vv^\e) \dnab q^\e \big].
	\end{align*}
	Note that $\dv \uv^\vs = 0.$ Then one has that
	\begin{equation}
		\int_0^{+\infty} \norm{h^\e}_{L^2}^2 dt \lesssim \int_0^{+\infty} \big(\norm{\nabla b^\e}_{H^1}^2 + \norm{\nabla \vv^\e}_{H^2}^2 \big) dt \lesssim 1.
	\end{equation}
	Then it follows from \cite{Metivier2001} or \cite[Theorem 8.3]{Alazard06} that 
	\begin{equation}
		q^\e \to 0  \quad \text{strongly in } \ L^2((0,T); L_{loc}^2).
	\end{equation}

	As for $\vv^\e,$ it follows from \cref{equ-Ue}$_1$ that
	\begin{equation}
		\dv \vv^\e = - \eta^\e (\e \p_t q^\e + \e \uv^\e \dnab q^\e).
	\end{equation}
	Thus, it suffices to prove that
	\begin{equation}\label{lim-pt}
		\e \p_t q^\e \to 0 \quad \text{strongly in } \ L^2((0,T); L_{loc}^2).
	\end{equation}
	By denoting that $\widetilde{q}^\e = \e \p_t q^\e, $  then $$ \sup_{t\geq 0} \norm{\widetilde{q}^\e}_{H^2} \lesssim \sup_{t\geq 0} \norm{\nabla (q^\e, \vv^\e)}_{H^2} \lesssim 1. $$
	It follows from \cref{equ-wave} that
	\begin{equation}
		\e^2 \p_t (\eta^\e \p_t \widetilde{q}^\e ) 
        - \dv \big(\frac{1}{\rho^\e} \nabla\widetilde{q}^\e \big) = \e \widetilde{h}^\e,
	\end{equation}
	where
	\begin{align*}
		\widetilde{h}^\e & = \e \p_t h^\e 
        -\dv\Big(\frac{\p_t \rho^\e}{\abs{\rho^\e}^2}\nabla q^\e\Big) - \e^2 \p_t\big( \p_t \eta^\e \p_t q^\e \big).
	\end{align*}
	Using \cref{equ-Ue,uni-bdd}, one has
	\begin{align*}
		\int_0^{+\infty} \norm{\widetilde{h}^\e}_{L^2}^2 dt &  \lesssim \int_0^{+\infty} \big( \norm{\nabla q^\e}_{H^2}^2 + \norm{\nabla \vv^\e}_{H^4}^2\big) dt \lesssim 1.
	\end{align*}
	Thus, using \cite[Theorem 8.3]{Alazard06} again, one can obtain \cref{lim-pt}.

	With the convergence in $L^2(0,T; L^2_{loc})$ and the uniform boundedness \cref{uni-bdd}, one can use the interpolation to achieve \cref{lim of pressure}.
	
\end{proof}

\begin{Lem}\label{lem: lim of shear}
    For all $T>0$, there exits a subsequence of  $\vv^\e$ which converges strongly in $L^2((0,T);\,H^{4-}_{loc})$ to the limit $\vv^0$. 
\end{Lem}
\begin{proof}
Denote the projection operator $ \Pi:= \mathrm{Id}-\nabla\Delta^{-1}\dv $ and set $\wv^\e := \Pi \vv^\e$. 
Then $ \wv^\e$ is uniformly bounded in $C^0(0,+\infty;\,H^4)$ and \cref{equ-Ue}$_2$ implies that
\begin{equation}\label{eq: wv}
    \rho^\e(\p_t+\uv^\e\cdot\nabla)\wv^\e 
    -\mu\lap\wv^\e
    =g^\e,
\end{equation}
where
\begin{equation*}
    g^\e=[\Pi,\, \rho^\e(\p_t+\uv^\e\cdot\nabla)]\vv^\e 
    -\Pi(\e b^\e \p_t \uvs + \rho^\e v^\e_3\p_3\uvs).
\end{equation*}
Using \cref{uni-bdd}, one has that $\p_t\wv^\e$ is uniformly bounded in $L^2(0,+\infty;\,H^{3})$. Then the Aubin-Lions lemma tells that a subsequence satisfies
\begin{equation}
    \text{$\wv^\e\to \Pi \vv^0 = \vv^0$ strongly in $C^0(0,+\infty;\,H^{4-})$},
\end{equation}
where we have used the result in \eqref{lim of pressure}$_2$ that $\dv \vv^0 = 0$. Thus, $\vv^\e\to\vv^0$ strongly in $L^2((0,T);\,H^{4-})$. 
\end{proof}

With the key steps in Lemmas \ref{lem: lim of pressure}--\ref{lem: lim of shear}, Theorem \ref{thm: lim} follows from \cite{Metivier2001,Alazard06} for a subsequence of $\{ (\rho^\e, \uv^\e) \}$.
At last, for the Cauchy problem, \cref{InNS,ic-incom}, the classical solutions belonging to
$$
\mathbb{B} := \{ \uv = \uv^\vs + \vv: \vv \in C(0,+\infty; H^4), \nabla \vv \in L^2(0,+\infty; H^4) \},
$$
are unique. Thus, the strong convergences above are actually valid for the full sequence.


\vspace{.3cm}

\appendix

\section{A priori bounds}

The appendix is used to prove the \textit{a priori} bounds in \cref{Sec-bound}.

\subsection{Proof of some a priori bounds}\label{App-bound}

The estimates \cref{aL2-anti,aL2-od,aL2-md,aL2-pert} follow directly from \cref{M&nu,poincare,bdd-diff}.
To show \cref{aLinf-anti}, one can use the Gagliardo--Nirenberg inequality to get that for $j=0,1,$
\begin{equation*}
	\begin{aligned}
		\norm{\p_3^j (\Phi, \Psi_3, Z_3)}_{L^\infty} & \lesssim \norm{\p_3^{j+1} (\Phi, \Psi_3, Z_3)}_{L^2}^{\frac{1}{2}} \norm{\p_3^j (\Phi, \Psi_3, Z_3)}_{L^2}^{\frac{1}{2}} \\
		&  \lesssim \nu (t+1)^{-\frac{j}{2}},
	\end{aligned}
\end{equation*}
and for $j=2,3,$ 
\begin{align*}
	\norm{\p_3^j (\Phi, \Psi_3, Z_3)}_{L^\infty} & \lesssim \norm{\p_3^{j+1} (\Phi, \Psi_3, Z_3)}_{L^2}^{1-\frac{1}{2(j-1)}} \norm{\p_3^2 (\Phi, \Psi_3, Z_3)}_{L^2}^{\frac{1}{2(j-1)}} \\
	& \lesssim M^{\frac{3}{4}} \nu^{\frac{1}{4}} (t+1)^{-\frac{3}{4}}.
\end{align*}
The other two estimates in \cref{aLinf-anti} can be obtained similarly, and \cref{aLinf-od} follows from \cref{aLinf-anti} directly.

It follows from \cref{Lem-GN} and the Poincar\'{e} inequality that for $j=0,1,$
\begin{align*}
	\norm{\nabla^j (\phi^\md, \psi^\md, \wv^\md)}_{L^\infty} & \lesssim \sum_{k=2}^{3} \norm{\nabla^{j+2} (\phi^\md, \psi^\md, \wv^\md)}_{L^2}^{\frac{k}{4}} \norm{\nabla^j (\phi^\md, \psi^\md, \wv^\md)}_{L^2}^{1-\frac{k}{4}}, \\
	& \lesssim \sum_{k=2}^{3} \norm{\nabla^3 (\phi^\md, \psi^\md, \wv^\md)}_{L^2}^{\frac{k}{4}} \norm{\nabla (\phi^\md, \psi^\md, \wv^\md)}_{L^2}^{1-\frac{k}{4}} \\
	& \lesssim M^{\frac{3}{4}} \nu^{\frac{1}{4}} (t+1)^{-\frac{1}{2}}.
\end{align*}
Meanwhile, the second inequality above yields that
\begin{align*}
	\norm{(\phi^\md, \psi^\md, \wv^\md)}_{W^{1, \infty}} & \lesssim  \norm{\nabla^3 (\phi^\md, \psi^\md, \wv^\md)}_{L^2} \lesssim M (t+1)^{-\frac{3}{4}}.
\end{align*}
The proof of \cref{aLinf-md} is completed.
At last, \cref{aLinf-pert} can follow from a combination of \cref{aL2-od,aL2-md}.


\subsection{Proof of Lemma \ref{Lem-rel}}\label{App-rel}
Let $i \in \{1,2,3\}$ be fixed.
~\\
i) It follows from \cref{wv-od,aL2-anti} that 
\begin{equation}\label{ineq-5}
	\begin{aligned}
		\pm \norm{\p_3^j w_i^\od}_{L^2} & \lesssim \pm \norm{\p_3^{j+1} Z_i}_{L^2} + \e \nu (t+1)^{-\frac{2j+1}{4}}, \qquad j=0,1,2.
	\end{aligned}
\end{equation}
Note that $\zeta = \frac{\wv}{\rho}.$ Then one has that
\begin{equation}\label{zeta-0}
	\pm \norm{\zeta_i}_{L^2}  \lesssim \pm \norm{w_i}_{L^2} \lesssim \pm \big(\norm{\p_3 Z_i}_{L^2} + \norm{w_i^\md}_{L^2}\big) + \e \nu (t+1)^{-\frac{1}{4}}.
\end{equation}
Similarly, it also follows from \cref{ineq-5,aL2-pert,bdd-diff,aLinf-pert} that
\begin{equation}\label{zeta-1}
	\begin{aligned}
		\pm \norm{\nabla \zeta_i}_{L^2} & \lesssim \pm \norm{\nabla w_i}_{L^2} + \e \chi (t+1)^{-1} \norm{w_i}_{L^2} + \e\norm{w_i}_{L^\infty} \norm{\nabla \phi}_{L^2}, \\
		& \lesssim \pm \big(\norm{\p_3^2 Z_i}_{L^2} + \norm{\nabla  w_i^\md}_{L^2}\big) + \e (\chi+\nu) M (t+1)^{-\frac{3}{4}},
	\end{aligned}
\end{equation}
and
\begin{equation}\label{zeta-2}
	\begin{aligned}
		\pm \norm{\nabla^2 \zeta_i}_{L^2} &  \lesssim \pm \norm{\nabla^2 w_i}_{L^2} + \e \chi M (t+1)^{-\frac{7}{4}} +  \e  \norm{\nabla\phi}_{L^\infty} \norm{\nabla\wv}_{L^2} \\
		& \quad + \e  \norm{\wv}_{L^\infty} \big(\norm{\nabla^2 \phi}_{L^2} + \norm{\nabla \phi}_{L^\infty} \norm{\nabla \phi}_{L^2} \big) \\
		& \lesssim \pm \big(\norm{\p_3^3 Z_i}_{L^2} + \norm{\nabla^2 w_i^\md}_{L^2}\big) + \e M^2 (t+1)^{-\frac{5}{4}}.
	\end{aligned}
\end{equation}
Thus, collecting \cref{zeta-0,zeta-1,zeta-2} yields \cref{rel-0-2}.

This implies that
\begin{equation}\label{app-L2-zeta}
	\norm{\zeta}_{L^2} \lesssim M (t+1)^{-\frac{1}{4}}, \qquad \norm{\nabla \zeta}_{H^1} \lesssim M (t+1)^{-\frac{3}{4}}.
\end{equation}
In addition, one can get from \cref{bdd-diff,app-L2-zeta} that
\begin{equation}\label{ineq-3}
	\begin{aligned}
		\pm \norm{\nabla^3 \wv}_{L^2} & \lesssim \pm \norm{\nabla^3 \zeta}_{L^2} + \e \chi M (t+1)^{-\frac{9}{4}} + \e \Big[ \norm{\nabla \phi}_{L^\infty} \norm{\nabla^2 \zeta}_{L^2} \\
		& \quad + \norm{\nabla^2 \phi}_{L^4} \norm{\nabla \zeta}_{L^4} + \norm{\zeta}_{L^\infty} \norm{\nabla^3 \phi}_{L^2} \Big].
	\end{aligned}
\end{equation}
It follows from \cref{G-N-type-1,app-L2-zeta,aL2-pert} that 
\begin{equation}\label{ineq-9}
	\begin{aligned}
		\norm{\nabla \zeta}_{L^4} & \lesssim \sum_{k=1}^{3} \norm{\nabla^2 \zeta}_{L^2}^{\frac{k}{4}} \norm{\nabla \zeta}_{L^2}^{1-\frac{k}{4}}  \lesssim M (t+1)^{-\frac{3}{4}}, \\ \norm{\nabla^2 \phi}_{L^4} & \lesssim \norm{\nabla^2 \phi}_{H^1} \lesssim M(t+1)^{-\frac{3}{4}}, \\
		\norm{\zeta}_{L^\infty} & \lesssim \norm{\nabla \zeta}_{L^2}^{\frac{1}{2}} \norm{ \zeta}_{L^2}^{\frac{1}{2}} + \sum_{k=2}^{3} \norm{\nabla^2 \zeta}_{L^2}^{\frac{k}{4}} \norm{ \zeta}_{L^2}^{1-\frac{k}{4}} \lesssim M (t+1)^{-\frac{1}{2}}.
	\end{aligned}
\end{equation}
This, together with \cref{aL2-pert,aLinf-pert}, yields that
\begin{equation}\label{zeta-3}
	\pm \norm{\nabla^3 \wv}_{L^2} \lesssim \pm \norm{\nabla^3 \zeta}_{L^2} + \e M^2 (t+1)^{-\frac{5}{4}}.
\end{equation}
Thus, when $\e^{\frac{1}{2}} M^2 \leq 1 $, one has
\begin{equation}\label{rel-3}
	\pm \norm{\nabla^3 \zeta}_{L^2} \lesssim \pm \norm{\nabla^3 \wv}_{L^2} + \e^{\frac{1}{2}} M (t+1)^{-\frac{5}{4}}.
\end{equation}

\vspace{.1cm}
~\\
ii) 
It follows from \cref{wv-md,poincare} that for $j=0,1,2,$
\begin{align*}
	\pm \norm{\nabla^j \wv^\md}_{L^2} & \lesssim \pm \norm{\nabla^j \zeta^\md}_{L^2} + \norm{\p_3 \rhot}_{W^{1,\infty}} \norm{\nabla^j \zeta^\md}_{L^2} \\
	& \quad + \e \big(\norm{(\phi^\od, \zeta^\od)}_{W^{2,\infty}} + \norm{(\phi^\md, \zeta^\md)}_{W^{1,\infty}}\big) \norm{\nabla^j (\phi^\md, \zeta^\md)}_{L^2} \\
	& \lesssim \pm \norm{\nabla^j \zeta^\md}_{L^2} + \e (\chi + M) \norm{\nabla^j (\phi^\md, \zeta^\md)}_{L^2}.
\end{align*}
Then if $\e^{\frac{1}{2}} M \leq 1 $ and $\e$ is suitably small, then \cref{rel-zeta-w} is true.  The proof of \cref{thm: lim} is finished.

\vspace{.3cm}

\textbf{Acknowledgments.} 
The authors want to thank Prof. Qiangchang Ju and Dr. Jiawei Wang for helpful discussions in the context of the incompressible limit.

Qian Yuan is partially supported by the National Natural Science Foundation of China (No.12201614), Youth Innovation Promotion Association of CAS (No.2022003) and CAS Project for Young Scientists in Basic Research (No.YSBR-031). 
Wenbin Zhao is partially supported by the National Natural Science Foundation of China (No.12401286), the Fundamental Research Funds for the Central Universities and the Research Funds of Renmin University of China (No.24XNKJ09).

\vspace{.1cm}

\textbf{Data Availability Statement}. Data sharing is not applicable to this article, as no data sets were generated or analyzed during the current study.

\textbf{Conflict of interest}. The authors declare that there is no conflict of interest.


\providecommand{\bysame}{\leavevmode\hbox to3em{\hrulefill}\thinspace}
\providecommand{\MR}{\relax\ifhmode\unskip\space\fi MR }
\providecommand{\MRhref}[2]{%
	\href{http://www.ams.org/mathscinet-getitem?mr=#1}{#2}
}
\providecommand{\href}[2]{#2}




\end{document}